%% file: main.tex
\documentclass[10pt]{article}
\usepackage[T1]{fontenc}
\usepackage{textalpha}

\usepackage{amsmath}
\usepackage{amsthm}
\usepackage{amsfonts}
\usepackage{amssymb}
\usepackage{bbm}
\usepackage{mathrsfs}

\usepackage{hyperref}
\usepackage[nameinlink,capitalize]{cleveref}
\usepackage{chngcntr}
\usepackage{cite}
\usepackage[nottoc,numbib]{tocbibind}

\usepackage{pgf,tikz,pgfplots}
\usepgfplotslibrary{groupplots}
\usetikzlibrary{arrows}
\usetikzlibrary{arrows.meta}

\usepackage{subcaption}

\newtheorem{observation}{Observation}

\counterwithin{theorem}{section}
\counterwithin{lemma}{section}
\counterwithin{remark}{section}
\counterwithin{observation}{section}
\counterwithin{definition}{section}

\newcommand{\vb}[1]{\mathbf{#1}}
\newcommand{\bm}[1]{\boldsymbol{#1}}

\newcommand{\wrt}{\text{w.r.t.}}

\newcommand{\one}{\bm{\mathbbm{1}}}

\newcommand{\Vol}{V}
\newcommand{\Area}{A}
\newcommand{\curv}{s}

\newcommand{\norm}[1]{\| {#1} \|}
\newcommand{\con}[2]{\langle {#1} , \, {#2} \rangle}
\newcommand{\jump}[1]{\ensuremath{[\![#1]\!]} }
\DeclareMathOperator{\at}{\bigg\vert}

\DeclareMathOperator*{\argmin}{arg\,min}
\DeclareMathOperator*{\argmax}{arg\,max}
\DeclareMathOperator{\polar}{\mathrm{polar}}
\DeclareMathOperator{\sym}{\mathrm{sym}}
\DeclareMathOperator{\dev}{\mathrm{dev}}
\DeclareMathOperator{\sph}{\mathrm{sph}}
\DeclareMathOperator{\skw}{\mathrm{skew}}
\DeclareMathOperator{\Anti}{\mathrm{Anti}}
\DeclareMathOperator{\axl}{\mathrm{axl}}
\DeclareMathOperator{\tr}{\mathrm{tr}}
\DeclareMathOperator{\id}{\mathrm{id}}

\DeclareMathOperator{\spa}{\mathrm{span}}
\DeclareMathOperator{\ad}{\mathrm{ad}}

\newcommand{\dd}{\mathrm{d}}
\newcommand{\Grad}{\mathrm{D}}
\DeclareMathOperator{\curl}{\mathrm{curl}}
\DeclareMathOperator{\Curl}{\mathrm{Curl}}
\DeclareMathOperator{\di}{\mathrm{div}}
\DeclareMathOperator{\Di}{\mathrm{Div}}

\newcommand{\R}{\mathbb{R}}
\newcommand{\SO}{\mathrm{SO}}
\newcommand{\so}{\mathfrak{so}}
\renewcommand{\sl}{\mathfrak{sl}}
\newcommand{\Sym}{\mathrm{Sym}}

\newcommand{\GL}{\mathrm{GL}}

\newcommand{\Diag}{\mathrm{Diag}}
\renewcommand{\P}{\mathrm{P}}

\newcommand{\C}{\mathit{C}}
\newcommand{\Le}{\mathit{L}^2}
\newcommand{\Hone}{\mathit{H}^1}

\newcommand{\Hcurl}[1]{\mathit{H}(\mathrm{curl}#1)}

\newcommand{\CG}{\mathcal{CG}}

\newcommand{\Ned}{\mathcal{N}}
\newcommand{\GEO}{\mathcal{GE}}

\newcommand{\defmap}{\bm{\varphi}}

\newcommand{\Double}{\bm{\mathfrak{P}}}
\newcommand{\Distress}{\bm{\mathfrak{D}}}
\newcommand{\Rdistress}{\bm{\mathfrak{T}}}
\newcommand{\Piola}{\bm{P}}
\newcommand{\Rpiola}{\bm{T}}
\newcommand{\defgrad}{\bm{F}}
\newcommand{\Coss}{\overline{\bm{R}}}
\newcommand{\Stretch}{\overline{\bm{U}}}
\newcommand{\Curva}{\overline{\bm{\mathfrak{A}}}}
\newcommand{\Wry}{\overline{\bm{\mathfrak{K}}}}

\newcommand{\Cm}{\mathbb{C}}
\newcommand{\J}{\mathbb{J}}

\newcommand{\muc}{\mu_{\mathrm{c}}}
\newcommand{\Lc}{L_\mathrm{c}}
\newcommand{\Lm}{\mathbb{L}}

\newcommand{\Energy}{\mathcal{I}}

\newcommand{\Action}{\mathcal{A}}

\newcommand{\AS}[1]{{\color{black} #1}}

\makeatletter
\let\@fnsymbol\@arabic
\makeatother

\title{A structure-preserving discretisation of SO(3)-rotation fields for finite Cosserat micropolar elasticity}

\author{
\normalsize{Lucca Schek}\thanks{Lucca Schek,
Fachgebiet für Kontinuumsmechanik und Materialtheorie, Institut für Mechanik, 
Fak. V,
Technische Universität Berlin,
Einsteinufer 5
D-10587 Berlin, Germany, email: l.schek@tu-berlin.de}
, \quad
    \normalsize{Peter Lewintan}\thanks{Peter Lewintan, Karlsruhe Institute of Technology, Englerstrasse 2, 76131 Karlsruhe, Germany, email: peter.lewintan@kit.edu}
, \quad
\normalsize{Wolfgang M\"uller}\thanks{Wolfgang M\"uller,
Fachgebiet für Kontinuumsmechanik und Materialtheorie, Institut für Mechanik, 
Fak. V,
Technische Universität Berlin,
Einsteinufer 5
D-10587 Berlin, Germany, email: wolfgang.h.mueller@tu-berlin.de}
, \quad
\normalsize{Ingo Muench}\thanks{Ingo Muench, Institute of Structural Mechanics, Statics and Dynamics, Technische Universit\"at Dortmund, August-Schmidt-Str. 8, 44227 Dortmund, Germany, email: ingo.muench@tu-dortmund.de}
 , \quad
    \normalsize{Andreas Zilian}\thanks{Andreas Zilian, Institute of Computational Engineering and Sciences, Department of Engineering, Faculty of Science, Technology and Medicine, University of Luxembourg, 6 Avenue de la Fonte, L-4362 Esch-sur-Alzette, Luxembourg, email: andreas.zilian@uni.lu}
    , \\
    \normalsize{St\'ephane P. A. Bordas}\thanks{St\'ephane P. A. Bordas, Institute of Computational Engineering and Sciences, Department of Engineering, Faculty of Science, Technology and Medicine, University of Luxembourg, 6 Avenue de la Fonte, L-4362 Esch-sur-Alzette, Luxembourg, email: stephane.bordas@alum.northwestern.edu}\,
    \thanks{Computational Modelling and Data Science, Department of Engineering, University of Exeter, United Kingdom, email: s.bordas@exeter.ac.uk}
    , \quad
	\normalsize{Patrizio Neff}\thanks{Patrizio Neff, Chair for Nonlinear 
		Analysis and Modelling, Faculty of Mathematics, Universit\"{a}t Duisburg-Essen,
		Thea-Leymann Str. 9, 45127 Essen, Germany, email: patrizio.neff@uni-due.de}
    \quad
    \normalsize{and} \quad
    \normalsize{Adam Sky$^*$}\thanks{$^*$Corresponding author: Adam Sky, Institute of Computational Engineering and Sciences, Department of Engineering, Faculty of Science, Technology and Medicine, University of Luxembourg, 6 Avenue de la Fonte, L-4362 Esch-sur-Alzette, Luxembourg, email: adam.sky@uni.lu}
}

\begin{document}

\maketitle

\begin{abstract}
We introduce a new method, dubbed \textbf{\textit{Geometric Structure-Preserving Interpolation (Γ-SPIN)}} to preserve physics-constraints inherent in the material parameter limits of the finite-strain Cosserat micropolar model.
The method advocates to interpolate the Cosserat rotation tensor using geodesic elements, which maintain objectivity and correctly represent curvature measures. At the same time, it proposes relaxing the interaction between the rotation tensor and the deformation tensor to alleviate locking effects. This relaxation is achieved in two steps. First, the regularity of the Cosserat rotation tensor is reduced by interpolating it into the N\'ed\'elec space. Second, the resulting field is projected back onto the Lie-group of rotations. Together, these steps define a lower-regularity projection-based interpolation. The construction allows the discrete Cosserat rotation tensor to match the polar part of the discrete deformation tensor. This ensures stable behaviour in the asymptotic regime as the Cosserat couple modulus tends to infinity, which constrains the model towards its couple-stress limit. We establish the consistency, stability, and optimality of the proposed method through several benchmark problems. The study culminates in a demonstration of its efficacy on a more intricate curved domain, contrasted with outcomes obtained from conventional interpolation techniques.  
\\
\vspace*{0.25cm}
\\

\end{abstract}

\section{Introduction}

Classical continuum mechanics models a material as a continuous collection of points, each possessing three translational degrees of freedom \cite{Truesdell1992}. By the balance of angular momentum, the Cauchy stress tensor is symmetric $\bm{\sigma} = \bm{\sigma}^T$ and the material response depends solely on the deformation map $\defmap:\Vol \subset \R^3 \to \Vol_\varphi \subset \R^3$ through its first gradient $\defgrad = \Grad \defmap$ \cite{Truesdell1984}. While this framework accurately describes many macroscopic phenomena, it fails in situations where the strain field varies rapidly \cite{Sansour2009}. In regions with high strain gradients or pronounced microstructural interactions, such as in foams, composites, granular media and metamaterials, internal rotations and couple-stresses become relevant, and the classical theory is no longer adequate \cite{Diebels2002}.

The underlying rotations and couple-stresses of material points are naturally accounted for in the Cosserat micropolar model \cite{Abreu2021}. The theory enriches classical continuum mechanics by introducing independent rotational degrees of freedom for each material point, thereby capturing the influence of microstructure on the macroscopic response \cite{Diepolder1991}. In this framework, every point of the continuum is described by a position vector $\defmap$ and an associated rotation tensor $\Coss:\Vol \subset \R^3 \to \SO(3)$ that represents the local microrotation \cite{Neff2006}, which is independent of the macro-rotation $\bm{R}=\polar \bm{F}$. 
The original motivation of François and Eugène Cosserat, formulated in their 1909 work \emph{Théorie des corps déformables}, was to construct a unified field theory of deformable bodies in which the balance of forces and balance of moments follow from a single variational principle invariant under Euclidean transformations \cite{Cosserat}. As a result, the Cosserat micropolar model can be conveniently described through the minimisation of an action-functional $\Action(\defmap,\Coss) \to \min \; \wrt \{\defmap,\Coss\}$ \cite{Neff2006}. 

Nowadays, the three-dimensional Cosserat micropolar continuum is widely employed in problems where microrotations, couple-stresses, and intrinsic length scales govern material behaviour \cite{Schek2026Pamm}. It effectively describes architected and chiral metamaterials with bending-dominant microstructures, capturing size-dependent stiffness and twist-coupling effects \cite{Liu2012,Frenzel2021}, and it is used in geosciences to model granular assemblies and fault-block interactions in seismic processes \cite{Lin2015}. At the nanoscale, the same framework accounts for magneto-mechanical couplings such as flexomagnetism in centrosymmetric magnetic materials \cite{SkyFlex}.
For thin structures, Cosserat shell theory reduces the three-dimensional model to a rotationally enriched two-dimensional formulation that enables natural mixed-dimensional couplings \cite{Sky2024}. In this context, recent $\Gamma$-convergence results clarify the relation between three-dimensional and reduced energies and provide explicit curvature expressions for efficient modelling and parameter identification \cite{MohammadiSaem2025,NEFF2010}, while further analysis establishes applicability to non-orientable manifolds via domain reconstruction \cite{Nebel2023}. For dynamic Cosserat plate models, appropriate initial and boundary conditions have been identified to ensure well-posedness under transient waves in thin microstructured media \cite{Varygina2018}.
Notwithstanding, the most widely used specialisation of the theory is the one-dimensional Cosserat rod, which offers a geometrically exact description of slender beams. In fact, the Simo--Reissner beam \cite{Simo1986,Simo1985,Simo1991} arises as a particular case of it, incorporating shear deformation and finite rotations under the same kinematic assumptions. Over the past decades, Cosserat rods have become central in engineering and applied physics, ranging from real-time simulation and control in soft robotics \cite{Till2019,Nguyen2023}, modelling biofilaments such as DNA and collagen in biomechanics and biophysics \cite{Gazzola2018,PalanthandalamMadapusi2011,Amuasi}, to capturing the nonlinear dynamics of drill strings, risers, and mooring lines in offshore and energy applications \cite{Tengesdal2023,Taran2025,Kim2019FlexibleRisers}, among many other developments and applications, see for example \cite{Kumar2026,Rajagopal2023,Bartholdt2021,Mochiyama2020}.
By virtue of its popularity, the numerical treatment of the Cosserat rod has received significant attention in the computational community. The main challenge here lies in maintaining objectivity and geometric consistency under large rotations while avoiding artefacts such as membrane- or shear-locking. To overcome this, geometrically exact formulations based on Lie group methods and rotation interpolation on $\SO(3)$ have been introduced to ensure frame invariance \cite{Sonneville2014,Ghosh2009,Crisfield1999}, while locking effects have been mitigated through mixed formulations, selective integration \cite{BaierSaip2020}, and Petrov--Galerkin schemes employing distinct trial and test spaces for selected energy contributions \cite{Harsch2023,Niemi2011}. Forming the focal point of this work, it is important to notice that locking \cite{Ju2025,Sky2023,Demasi2025} can be interpreted as a structure-preservation issue inherent to Cosserat continua of any dimensionality \cite{Sky2024,Dziubek2025,Boon2025}.


The geometrically exact Cosserat formulation employs the rotation tensor $\Coss:\Vol \to \SO(3)$ to describe the independent finite rotations of material points. The effects of these independent rotations are measured by two tensors \cite{Schek2026}. Strains are measured by the Cosserat strain tensor $\Coss^T \defgrad - \one$, composed of the Biot-type stretch tensor $\Stretch = \Coss^T \defgrad$ and the second-order identity tensor $\one$. Curvatures are measured either by the wryness tensor $\Wry = \axl(\Coss^T \Coss_{,i})\otimes \vb{e}_i$ in the classical Eringen form \cite{Eringen1999}, or by the micro-dislocation tensor $\Curva = \Coss^T \Curl \Coss$ in dislocation form \cite{Ghiba}. 
Inherently, the coupling mechanism of the Cosserat strain tensor can be found in all \textbf{derived} models, and is responsible for well-known locking effects in \AS{the corresponding} beam and shell theories. 
But more generally, the problem is one of structure-preservation and can be stated in terms of how does one construct a discrete $\Coss_h$ that precisely matches \AS{$\polar\defgrad_h = \bm{R}_h$}\AS{, which is given by the polar decomposition $\defgrad_h = \bm{R}_h \bm{U}_h$}.
To clarify, let some material parameter $\muc \to + \infty$, then the coupling energy governed by it $\muc \norm{\Coss_h^T\defgrad_h - \one}^2$ must compensate to stay finite. The latter implies $\norm{\Coss_h^T \defgrad_h- \one}^2 \to 0$ going to zero at at-least the same rate as $\muc$ goes to infinity. 
Now, in standard discretisation approaches the deformation gradient is constructed as $\defgrad_h = \Grad \defmap_h$, where the deformation map $\defmap$ is approximated by classical $\C^0(\Vol)$-continuous Lagrange elements $\defmap_h \in \R^3 \otimes \CG^p(\Vol)$. Clearly, since $p \geq 1$ and $\CG^p(T) = \P^p(T)$ is element-wise a full polynomial space, its gradient already contains the space of constants $\R^{3 \times 3} \subseteq \Grad [\R^3 \otimes \CG^p(T)]$ \cite{SKY2024104155}. However, any discrete $\Coss_h$ must lie in a conforming subspace of $\Hone[\Vol;\SO(3)]$, where the rotational field takes values in $\SO(3)$, a nonlinear Lie-group manifold rather than a vector space.
In other words, the lacking ability to approach the constant identity $\one$ is fully due the mismatch in the discretisation of $\Coss_h^T$ and $ \defgrad_h$, which makes it plainly impossible to reconstruct $\Coss_h^T \Grad \defmap_h = \one$ point-wise. To see why, it is sufficient to apply the polar decomposition while demanding $\Coss_h^T \defgrad_h = \Coss^T \bm{R}_h \bm{U}_h = \one$. Since we established that the discrete deformation gradient $\bm{F}_h$ can generate a constant stretch tensor $\bm{U}_h = \one$, it is the incompatibility between $\bm{R}_h = \polar \defgrad_h$ and $\Coss_h$ \AS{in the weak sense} that leads to what we call locking.
Again, since $\SO(3)$ is not an algebraic space, there is no canonical polynomial interpolation for it. Consequently, all interpolations rely on some mapping from a finite element space. This can be achieved, amongst other methods, with the exponential map via its Euler--Rodrigues closed form formula \cite{Pimenta1993,Ibrahimbegovic1997,bauer2010,bauer2012}, Euler rotation matrices \cite{Argyris1982,Betsch1998}, quaternions \cite{Wasmer2024,Ghosh2008,erdelj2020,erdelj2024,Romero2004}, geodesic elements \cite{Sander2009,Sander2012,Sander2015,Greco2023,Greco2024,Crisfield1999}, or projection-based elements \cite{Grohs2019,Mller2022}. However, none of the above on its own can exactly match $\polar \defgrad_h$ for $\defmap_h \in \R^3 \otimes \CG^p(\Vol)$. The reason is simple: either the regularity, the interpolation method, or both, do not match. Simply put, since $\defmap \in \R^3 \otimes\Hone(\Vol)$ and $\Grad \defmap \in \R^3 \otimes \Hcurl{,\Vol}$ we have
\begin{align}
    \begin{matrix}
        \Coss & \in & \Hone[\Vol;\SO(3)] & \subsetneq & \polar [\R^3 \otimes \Hcurl{,\Vol}] & \ni & \bm{R} & = & \polar \defgrad \, .
        \\
        && \cup && \cup
        \\
        \Coss_h & \in & X_h(\Vol) &\subsetneq & Y_h(\Vol) & \ni & \bm{R}_h & = & \polar \defgrad_h \, ,
    \end{matrix}
\end{align}
for some conforming subspaces $X_h(\Vol)$ and $Y_h(\Vol)$.
By density arguments, the issue may not arise in the infinite-dimensional setting or in the weak sense with sufficiently rich finite element spaces, such as higher-order or macro-spaces \cite{Ainsworth2022,Ainsworth2025}. However, this argument does not generally persist in the discrete subspaces used in practical finite element discretisations. Accordingly, the topic of this work is the re-establishment of an exact match. The solution we propose extends the original ideas of the \textit{Mixed Interpolated Tensorial Components (MITC)} method \cite{Bathe1986,Stenberg1997,Brezzi1989,Falk2008} and Regge-interpolation in shells \cite{Neunteufel2021} to higher dimensions and nonlinear finite element spaces. Namely, we propose the following construction in the discrete action-functional
\begin{align}
    \boxed{
        \begin{aligned}
            \muc\norm{ \Coss_h^T \Grad \defmap_h -\one  }^2 \quad \to \quad \muc\norm{ \polar(\Pi_c^{p-1}\Coss_h)^T \Grad \defmap_h -\one  }^2 \, , && \begin{matrix}
                \defmap_h \in \R^3 \otimes \CG^p(\Vol) \, , \qquad \Coss_h \in \GEO^{q}(\Vol), \\ 
                \Pi_c^{p-1}: \GEO^{q}(\Vol) \to \R^3 \otimes \Ned_{II}^{p-1}(\Vol)
            \end{matrix} \, . 
        \end{aligned}
    }
\end{align}
That is, $\Coss_h$ is given by geodesic finite elements $\GEO^q(\Vol)$, but for coupling-terms it is interpolated into the N\'ed\'elec space $\Ned_{II}^{p-1}(\Vol)$ \cite{nedelec_mixed_1980,ndlec_new_1986,SKY2024104155,sky_higher_2023} via $\Pi_c^{p-1}$, before being projected back into the Lie-group of rotations $\SO(3)$, yielding a lower-regularity projection-based interpolation. 
Since $\defgrad_h = \Grad\defmap_h \in \R^3 \otimes\Ned_{II}^{p-1}(\Vol) \subset \R^3 \otimes \Hcurl{,\Vol}$, its polar lives in $\polar[\R^3 \otimes\Ned_{II}^{p-1}(\Vol)]$, such that an exact match becomes possible. 
By discretising the micro-rotation using geodesic finite elements $\Coss_h \in \GEO^q(\Vol)$, the true rates of change and thus the curvature measures are preserved, and the interpolation remains objective.
In recognition of its geometric foundation (ΓΕΩΜΕΤΡΙΑ) in terms of both geometric continuum mechanics and geometric finite element spaces, we dub the method \textbf{\textit{Geometric Structure-Preserving Interpolation (Γ-SPIN)}}.
We also mention that for Cosserat rods, our approach bears relation to applications of the Petrov--Galerkin method to alleviate shear-locking \cite{Harsch2023,Romero2004}, as by its specific choice of test functions the latter also implies a projection that eliminates spurious couplings. 

This work is structured as follows. First we introduce the Cosserat model, and discuss its limit cases of finite elasticity and couple-stress theory. Then, we present our structure-preserving approach, encompassing a discussion of the discretisation of $\SO(3)$-fields and possible choices of stable discrete pairings. Subsequently, we present numerical examples and discuss their results. Finally, we present our conclusions and outlook. 

\subsection{Notation}

The following notation is employed throughout this work.
Exceptions are made clear in the precise context.
\begin{itemize}
    \item Vectors are defined as bold lower-case letters $\vb{v}, \, \bm{\xi} \in \R^d$.
    \item Second-order tensors are denoted with bold capital letters $\bm{T}\in \R^{d \times d}$.
    \item Fourth-order tensors are designated by the blackboard format $\mathbb{C} \in \R^{d \times d \times d \dots}$.
    \item We denote the Cartesian basis as $\{\vb{e}_1, \, \vb{e}_2, \, \vb{e}_3\}$.
    \item Summation over indices follows the standard rule of repeating indices.
    \item The angle-brackets define scalar products of arbitrary dimensions $\con{\vb{v}}{\vb{u}} = v_i u_i$, $\con{\bm{T}}{\bm{F}} = T_{ij}F_{ij}$.
    \item The matrix product is used to indicate all partial-contractions between a higher-order and a lower-order tensor $\bm{T}\vb{v} = T_{ij} v_j \vb{e}_i$, $\mathbb{C}\bm{T} = C_{ijkl}T_{kl}\vb{e}_i \otimes \vb{e}_j$.
    \item The second-order identity tensor is defined via $\one = \vb{e}_i \otimes \vb{e}_i$, such that $\one \vb{v} = \vb{v}$. 
    \item Volumes, surfaces and curves of the physical domain are identified via $\Vol$, $\Area$ and $\curv$, respectively.   
    \item The space of rotation matrices is given by $\SO(d)$, reading $\SO(d) = \{ \bm{R} \in \R^{d \times d} \; | \; \det \bm{R} = 1 \, , \; \bm{R}^{-1} = \bm{R}^T \}$, and is associated with the polar projection operator $\polar: \R^{d\times d}\rightarrow\SO(d)$.
    \item We define the constant space of skew-symmetric second order tensors as $\so(d) = \{ \bm{T} \in \R^{d \times d} \; | \; \bm{T} = -\bm{T}^{T} \}$.
   \item The space $\so(d)$ is associated with the operators $\skw \bm{T} = (1/2)(\bm{T} - \bm{T}^T) \in \so(d)$, $\Anti \vb{v} = \vb{v} \times \one \in \so(3)$, and its inverse $\axl (\Anti \vb{v}) = \vb{v}$.
   \item The nabla operator is defined as $\nabla = \vb{e}_i \partial_i$.
   \item The left-gradient is given via $\nabla$, such that $\nabla \lambda = \nabla \otimes \lambda$.
    \item The right-gradient is defined for vectors and higher order tensors via $\Grad$, such that $\Grad \vb{v} = \vb{v} \otimes \nabla$.
    \item We define the vectorial divergence as $\di \vb{v} = \con{\nabla}{\vb{v}}$.
    \item The tensor divergence is given by $\Di \bm{T} = \bm{T} \nabla$, implying a single contraction acting row-wise.
    \item The vectorial curl operator reads $\curl \vb{v} = \nabla \times \vb{v}$
    \item For tensors the operator is given by $\Curl \bm{T} = -\bm{T} \times \nabla$, acting row-wise. 
\end{itemize}

\section{The nonlinear micropolar model}

In the following we shortly derive the quasi-static, geometrically nonlinear Cosserat micropolar model using the principle of virtual work. As it stands, there are two prominent formulations of the micropolar model, being the classical Eringen-form \cite{Eringen1999} and the dislocation-form \cite{Ghiba}. Their equivalence is due to the fact that the gradient of a $\SO(3)$-matrix field is fully controlled by its Curl \cite{SO3}. The equivalent material parameters of the two forms and corresponding transformation formulae are given in \cite{Ghiba}. However, neither the equivalence of the strong forms nor the relation between the various couple-stress measures are discussed in these works. 

To derive the Cosserat model, we begin by first introducing the Eulerian formulation.
Let $\Vol_\varphi \subset \R^3$ be an open and bounded domain with a sufficiently smooth boundary $\partial \Vol_\varphi$ representing the current configuration, the body forces are given by $\vb{f}_\varphi : \Vol_\varphi \to \R^3$, the couple-forces by $\vb{m}_\varphi: \Vol_\varphi \to \R^3$, the tractions on the boundary are given by $\vb{t}_\varphi: \partial \Vol_\varphi \to \R^3$, and the couple-tractions by $\bm{\mu}_\varphi: \Vol_\varphi \to \R^3$. The domain along with the forces and tractions is depicted in \cref{fig:dom}.
\begin{figure}
    \centering
    \input{figs/dom}
    \caption{Current configuration of a micropolar body $\Vol_\varphi \subset \R^3$ as given by the deformation map $\defmap : \Vol \to \Vol_\varphi$ and the independent orientation of material points $\Coss:\Vol \to \SO(3)$. Source terms in the domain are the body forces $\vb{f}_\varphi$ and the couple-forces $\vb{m}_\varphi$. Their fluxes on the boundary are given by the traction $\vb{t}_\varphi$ and couple-traction $\bm{\mu}_\varphi$. Notably, the orientation of material points is independent of deformation curves.}
    \label{fig:dom}
\end{figure}
Thus, the global equilibrium of linear momentum of statics reads
\begin{align}
    \int_{\Vol_\varphi} \vb{f}_\varphi \, \dd \Vol_\varphi + \int_{\partial \Vol_\varphi} \vb{t}_\varphi \, \dd \Area = 0 \,.
\end{align}
Now, by Cauchy's postulate there exists a tensor $\bm{\sigma}: \Vol_\varphi \to \R^{3 \times 3}$, such that 
\begin{align}
    \bm{\sigma} \, \vb{n}_\varphi \at_{\partial_{\Vol}} = \vb{t}_\varphi \at_{\partial_{\Vol}} \, ,
\end{align}
where $\vb{n}_\varphi$ is the unit normal vector field of the boundary of the domain $\vb{n}_\varphi: \partial \Vol_\varphi \to \R^3$. As per the postulate, requiring this assumption to hold for every arbitrarily small portion of the domain $\Delta \Vol_\varphi$ with $|\Delta \Vol_\varphi| \to 0$, while applying the divergence theorem implies
\begin{align}
    \int_{\Delta\Vol_\varphi} \vb{f}_\varphi \, \dd \Vol_\varphi + \int_{\partial \Delta \Vol_\varphi} \vb{t}_\varphi \, \dd \Area = \int_{\Delta\Vol_\varphi} \vb{f}_\varphi + \Di_\varphi \bm{\sigma} \, \dd \Vol_\varphi = 0  \, , 
\end{align}
from which we extract
\begin{align}
    -\Di_\varphi \bm{\sigma} = \vb{f}_\varphi \, , 
\end{align}
representing the local form of equilibrium of linear momentum. Note that unlike in the classical Cauchy continuum, we cannot presume that the stress tensor is symmetric. 
Now, the global form of the balance of angular momentum reads
\begin{align}
    \int_{\Vol_\varphi} \bm{\varphi} \times \vb{f}_\varphi + \vb{m}_\varphi \, \dd \Vol_\varphi + \int_{\partial \Vol_\varphi} \bm{\varphi} \times \vb{t}_\varphi + \bm{\mu}_\varphi \, \dd \Area = 0 \, ,
\end{align}
where $\bm{\varphi}: \Vol \to \Vol_\varphi$ is the deformation map from the reference to the current configuration.
In analogy to the balance of linear momentum, postulating the couple stress tensor $\bm{\Sigma} \vb{n}_\varphi$ with $\bm{\Sigma}: \Vol_\varphi \to \R^{3 \times 3}$, and applying the divergence theorem for a portion of the domain $\Delta \Vol_\varphi$ yields
\begin{align}
    \int_{\Delta\Vol_\varphi} \bm{\varphi} \times \vb{f}_\varphi + \vb{m}_\varphi \, \dd \Vol_\varphi + \int_{\partial \Delta\Vol_\varphi} \bm{\varphi} \times \vb{t}_\varphi + \bm{\mu}_\varphi \, \dd \Area = \int_{\Delta\Vol_\varphi} \bm{\varphi} \times \vb{f}_\varphi + \Di_\varphi(\bm{\varphi} \times \bm{\sigma}) + \vb{m}_\varphi + \Di \bm{\Sigma} \, \dd \Vol_\varphi = 0 \, ,
\end{align}
which can be reformulated into
\begin{align}
    \int_{\Delta\Vol_\varphi} \bm{\varphi} \times (\underbrace{\vb{f}_\varphi + \Di_\varphi \bm{\sigma}}_{=0}) + 2 \axl \bm{\sigma} + \vb{m}_\varphi + \Di \bm{\Sigma} \, \dd \Vol_\varphi = 0 \, ,
\end{align}
where we exploited the balance of linear momentum, and used that $\defmap_{,i} = \vb{e}_i$ and $(\vb{e}_i \times \bm{\sigma}) \vb{e}_i = 2 \axl \bm{\sigma}$.
As such, the local balance of angular momentum reads
\begin{align}
    -\Di \bm{\Sigma} - 2 \axl \bm{\sigma} = \vb{m}_\varphi \, ,
\end{align}
from which it is clear that the balance of angular momentum controls the skew-symmetric part of the stress tensor.
In order to derive the Lagrangian equations, we first define the forces $\vb{f} : \Vol \to \R^3$, the couple-forces $\vb{m}: \Vol \to \R^3$, the tractions $\vb{t}: \partial \Vol \to \R^3$, and the couple-tractions $\bm{\mu}:  \partial\Vol \to \R^3$, of the reference configuration. In the second step we exploit that two-forms transform according to the contravariant Piola transformation, such that the stress and couple stress tensors on the reference configuration read
\begin{align}
    \Piola = (\det \defgrad)\bm{\sigma}\defgrad^{-T} \quad \iff \quad \bm{\sigma} = \dfrac{1}{\det \defgrad} \Piola \defgrad^T \, , && \Double = (\det \defgrad)\bm{\Sigma}\defgrad^{-T} \quad \iff \quad \bm{\Sigma} = \dfrac{1}{\det \defgrad} \Double \defgrad^T \, ,
\end{align}
where $\defgrad = \Grad \defmap$ is the deformation gradient. Their respective divergences transform according to
\begin{align}
     \Di_{\varphi} \bm{\sigma} \simeq \Di \Piola \,, && \Di_{\varphi} \bm{\Sigma} \simeq \Di \Double \, ,
\end{align}
by virtue of Piola's identity. Finally, using the transformation of three-forms $\dd \Vol_\varphi = (\det \defgrad) \, \dd \Vol$ of densities
\begin{align}
    \vb{f}_\varphi = \dfrac{1}{\det \defgrad}\vb{f} \, , && \vb{m}_\varphi = \dfrac{1}{\det \defgrad } \vb{m} \, , && \axl \bm{\sigma} = \dfrac{1}{\det \defgrad}\axl(\Piola \defgrad^T) \, ,  
\end{align}
yields the Lagrangian equations 
\begin{subequations}
    \begin{align}
    -\Di\Piola &= \vb{f} \, , \\
    -\Di \Double  -2 \axl(\Piola \defgrad^T) &= \vb{m} \, .
\end{align}
\end{subequations}

By the principle of virtual work, the work conjugate of the balance of linear momentum is a virtual deformation $\delta \defmap$, such that the first equation yields
\begin{align}
    -\int_\Vol \con{\delta \defmap}{\Di\Piola} \,\dd \Vol = \int \con{\delta \defmap}{\vb{f}} \, \dd \Vol \quad \Rightarrow \quad \int_\Vol \con{\Grad \delta \defmap}{\Piola} \, \dd \Vol = \int_\Vol \con{\delta \defmap}{\vb{f}} \, \dd \Vol + \int_{\Area_N^{\defmap}} \con{\delta \defmap}{\underbrace{\bm{P} \vb{n}}_{\vb{t}}} \, \dd \Area \, .  
\end{align}
The work conjugate of the balance of angular momentum is a virtual rotation $\delta \bm{\omega}$, such that $\delta \Coss =  (\Anti \delta \bm{\omega})\Coss = \delta \bm{\omega}\times\Coss$, where $\Anti \delta \bm{\omega}: \Vol \to \so(3)$ is clearly skew-symmetric. Thus, testing with $\delta \bm{\omega}$, the virtual work of the second balance equation reads
\begin{align}
    &-\int_\Vol\con{\delta \bm{\omega}}{\Di \Double + 2 \axl (\Piola \defgrad^T)} \, \dd \Vol = \int_\Vol \con{\delta \bm{\omega}}{\vb{m}} \, \dd \Vol 
    \notag \\
    &\quad \Rightarrow \quad \int_\Vol \con{\Grad \delta \bm{\omega}}{\Double} - \con{\delta \bm{\omega}\times \defgrad}{\Piola} \, \dd \Vol = \int_\Vol \con{\delta \bm{\omega}}{\vb{m}} \, \dd \Vol +  \int_{\Area_N^{\Coss}} \con{\delta \bm{\omega}}{\underbrace{\Double \vb{n}}_{\bm{\mu}}} \, \dd \Area \, . 
\end{align}
Combined, the total virtual work principle reads
\begin{align}
    \int_\Vol\con{\Grad \delta \defmap - \delta \bm{\omega} \times \defgrad}{\Piola} + \con{\Grad \delta \bm{\omega}}{\Double} \, \dd \Vol = \int_\Vol \con{\delta \defmap}{\vb{f}} + \con{\delta \bm{\omega}}{\vb{m}} \, \dd \Vol + \int_{\Area_N^{\defmap}} \con{\delta \defmap}{\vb{t}} \, \dd \Area + \int_{\Area_N^{\Coss}} \con{\delta \bm{\omega}}{\bm{\mu}} \, \dd \Area \, .
\end{align}
 
In the Cosserat model, the kinematical fields are the deformation map $\defmap:\Vol \to \R^3$ and the independent rotation field $\Coss: \Vol \to \SO(3)$, which characterises the orientation of material points. From $\defmap$ we immediately retrieve the already aforementioned deformation gradient
\begin{align}
    \defgrad = \Grad \defmap \, , && \defgrad : \Vol \to \GL^+(3) \, .
\end{align}
Changes in the rotation $\Coss$ are governed by the curvature tensor $\Grad \Coss: \Vol \to \R^{3\times 3 \times 3}$. However, it is also possible to control the curvature using the wryness tensor
\begin{align}
    \Wry = \axl(\Coss^T \Coss_{,i})\otimes \vb{e}_i = -\axl(\Coss_{,i}^T \Coss)\otimes \vb{e}_i \, , && \Wry:\Vol \to \R^{3 \times 3} \, .
\end{align}
The relation to $\Grad \Coss$ is immediately clear from
\begin{align}
    \Grad \Coss = \Coss_{,i} \otimes \vb{e}_i = \overbrace{\Coss \, \Coss^T}^{\one} \Coss_{,i}  \otimes \vb{e}_i = \Coss \Anti \axl(\Coss^T\Coss_{,i}) \otimes \vb{e}_i \quad \Rightarrow \quad \Coss^T \Grad \Coss &= \Anti \axl(\Coss^T\Coss_{,i}) \otimes \vb{e}_i  \notag \\
    &= \Anti \Wry \, .
\end{align}
The deformation tensor $\defgrad$ and the wryness tensor $\Wry$ represent strain measures in the Cosserat model. Thus, their variations represent virtual strain measures. Whereas the variation of the deformation gradient is straight-forward 
\begin{align}
    \delta \defgrad = \delta \Grad \defmap = \Grad \delta \defmap \, ,
\end{align}
the variation of the wryness tensor is more involved, but can be expressed as
\begin{align}
    \delta \Wry &= -\delta\axl(\Coss_{,i}^T\Coss)\otimes \vb{e}_i = -\axl(\delta\Coss_{,i}^T\Coss + \Coss_{,i}^T\delta\Coss)\otimes \vb{e}_i 
    \notag \\
    &= -\axl(\Coss^T \Coss\delta\Coss_{,i}^T\Coss + \Coss^T \Coss\,\Coss_{,i}^T\delta\Coss\,\Coss^T \Coss)\otimes \vb{e}_i = -\axl(\Coss^T \Coss\delta\Coss_{,i}^T\Coss + \Coss^T [\Coss_{,i}\Coss^T][\Coss\delta\Coss^T] \Coss)\otimes \vb{e}_i
    \notag \\
    &= -\axl(\Coss^T \Coss\delta\Coss_{,i}^T\Coss + \Coss^T \Coss_{,i}\delta\Coss^T \Coss)\otimes \vb{e}_i = -\axl(\Coss^T [\Coss\delta\Coss_{,i}^T + \Coss_{,i}\delta\Coss^T]\Coss)\otimes \vb{e}_i 
    \notag \\
    &= -\Coss^T\axl( \Coss\delta\Coss_{,i}^T + \Coss_{,i}\delta\Coss^T)\otimes \vb{e}_i = \Coss^T \Grad \delta \bm{\omega} \, , 
\end{align}
since
\begin{align}
    \delta \bm{\omega} = \axl(\delta \Coss\,\Coss^T) = -\axl(\Coss \delta \Coss^T) \quad \Rightarrow \quad \Grad \delta \bm{\omega} = -\axl([\Coss \delta \Coss^T]_{,i}) \otimes \vb{e}_i = -\axl(\Coss_{,i} \delta \Coss^T + \Coss \delta \Coss^T_{,i}) \otimes \vb{e}_i \, .
\end{align}
To finally retrieve the virtual work principle in terms of the virtual strain measures, we also introduce the Biot-type stretch tensor and its corresponding variation
\begin{align}
    \boxed{
    \Stretch = \Coss^T \defgrad
    } \quad \Rightarrow \quad \delta \Stretch = \delta\Coss^T \defgrad + \Coss^T \delta \defgrad = \Coss^T\Coss(\delta\Coss^T \defgrad + \Coss^T \delta \defgrad) = \Coss^T(\delta \defgrad - \delta \bm{\omega} \times \defgrad)  \, .
\end{align}
Thus, the virtual work reads 
\begin{align}
    \int_\Vol \con{\delta \Stretch}{\Coss^T\Piola} + \con{\delta \Wry}{\Coss^T \Double} \, \dd \Vol = \int_\Vol \con{\delta \defmap}{\vb{f}} + \con{\delta \bm{\omega}}{\vb{m}} \, \dd \Vol + \int_{\Area_N^{\defmap}} \con{\delta \defmap}{\vb{t}} \, \dd \Area + \int_{\Area_N^{\Coss}} \con{\delta \bm{\omega}}{\bm{\mu}} \, \dd \Area \, .
\end{align}

In order to transition to dislocation form we first observe that
\begin{align}
    \Curl(\underbrace{\Coss^T\Coss}_{\one}) &= \Coss^T\Curl \Coss - \Coss^T_{,i}\Coss (\Anti\vb{e}_i) = \Coss^T\Curl \Coss - [\axl(\Coss^T_{,i}\Coss) \times \vb{e}_i \times \vb{e}_j] \otimes \vb{e}_j
    \notag \\
    & = \Coss^T\Curl \Coss - [\con{\axl(\Coss^T_{,i}\Coss)}{\vb{e}_j} \vb{e}_i - \con{\axl(\Coss^T_{,i}\Coss)}{\vb{e}_i}\vb{e}_j]  \otimes \vb{e}_j  
    \notag \\[1ex]
    & = \Coss^T\Curl \Coss - \vb{e}_i \otimes \axl(\Coss^T_{,i}\Coss) + \tr[\vb{e}_i \otimes \axl(\Coss^T_{,i}\Coss)]\one = \Coss^T\Curl \Coss + \Wry^T - \tr(\Wry^T) \one = 0 \, ,
\end{align}
such that
\begin{align}
    \Coss^T\Curl \Coss = \tr(\Wry^T) \one - \Wry^T \, .
\end{align}
Now, by taking the trace on both sides we extract
\begin{align}
    \tr (\Coss^T\Curl \Coss) = 2 \tr (\Wry^T) \quad \Rightarrow \quad \tr (\Wry^T) =  \dfrac{1}{2}\tr (\Coss^T\Curl \Coss) = \dfrac{1}{2}\tr [(\Curl \Coss)^T\Coss] \, ,
\end{align}
such that the wryness tensor can be written in terms of the Curl of the rotation field
\begin{align}
    \Wry = \dfrac{1}{2}\tr [(\Curl \Coss)^T\Coss]\one - (\Curl \Coss)^T\Coss \, .
\end{align}
Defining the dislocation density $\Curva$ as the back-rotated Curl, the final relation reads
\begin{align}
    \boxed{\Curva = \Coss^T \Curl \Coss} \quad \Rightarrow \quad \Wry = \dfrac{1}{2}\tr(\Curva^T)\one - \Curva^T \, ,
\end{align}
and the variation of the wryness tensor can be written as
\begin{align}
    \delta \Wry = \dfrac{1}{2}\tr(\delta\Curva^T)\one - \delta\Curva^T \, ,
\end{align}
which substituted into the virtual work of the couple-stress yields
\begin{align}
    \int_\Vol \con{\delta \Wry}{\Coss^T \Double} \, \dd \Vol &= \int_\Vol \dfrac{1}{2}\con{\tr(\delta\Curva^T) \one}{\Coss^T \Double} - \con{\delta\Curva^T}{\Coss^T \Double} \, \dd \Vol 
    \notag\\
    &= \int_\Vol \con{\delta\Curva^T}{(1/2)\tr(\Coss^T \Double)\one - \Coss^T \Double} \, \dd \Vol = \int_\Vol \con{\delta\Curva}{(1/2)\tr(\Double^T\Coss)\one - \Double^T\Coss} \, \dd \Vol  \, .
\end{align}
The latter motivates the definition of a new tensor, which we dub the dislocation-stress 
\begin{align}
    \Distress = \Coss [ (1/2)\tr(\Double^T\Coss)\one - \Double^T\Coss] \quad \iff \quad \Double = \tr(\Distress^T\Coss) \Coss - \Coss \Distress^T \Coss  \, , 
\end{align}
such that the virtual work can be written as
\begin{align}
    \int_\Vol \con{\delta \Stretch}{\Coss^T\Piola} + \con{\delta \Curva}{\Coss^T \Distress} \, \dd \Vol = \int_\Vol \con{\delta \defmap}{\vb{f}} + \con{\delta \bm{\omega}}{\vb{m}} \, \dd \Vol + \int_{\Area_N^{\defmap}} \con{\delta \defmap}{\vb{t}} \, \dd \Area + \int_{\Area_N^{\Coss}} \con{\delta \bm{\omega}}{\bm{\mu}} \, \dd \Area \, .
\end{align}
The variation of the micro-dislocation tensor reads
\begin{align}
    \delta(\Coss^T\Curl \Coss) &= \delta \Coss^T \Curl \Coss + \Coss^T \Curl \delta \Coss = \Coss^T\Coss \delta \Coss^T \Curl \Coss + \Coss^T \Curl (\delta \Coss\,\Coss^T\Coss) 
    \notag \\
    &= \Coss^T[\Coss \delta \Coss^T \Curl \Coss +  \Curl (\delta \Coss\,\Coss^T\Coss)] = \Coss^T[\Curl(\delta \bm{\omega} \times \Coss) - \delta \bm{\omega} \times  \Curl \Coss] \, .
\end{align}
Using this identity, the virtual work of $\Distress$ transforms according to 
\begin{align}
    \int_\Vol\con{\Coss^T[\Curl(\delta \bm{\omega} \times \Coss) - \delta \bm{\omega} \times  \Curl \Coss]}{\Coss^T\Distress} \, \dd \Vol &= \int_\Vol\con{\Curl(\delta \bm{\omega} \times \Coss) - \delta \bm{\omega} \times  \Curl \Coss}{\Distress} \, \dd \Vol 
    \notag \\
    &= \int_\Vol \con{\Anti \delta \bm{\omega}}{(\Curl \Distress)\Coss^T - \Distress (\Curl\Coss)^T} \, \dd \Vol
    \notag \\
    &= \int_\Vol 2\con{\delta \bm{\omega}}{\axl[(\Curl \Distress)\Coss^T - \Distress (\Curl\Coss)^T]} \, \dd \Vol \, ,
\end{align}
where the boundary term was omitted. By the latter, the strong form of the balance of angular momentum can also be expressed as
\begin{align}
    \axl[(\Curl \Distress)\Coss^T - \Distress (\Curl\Coss)^T] = \dfrac{1}{2}\vb{m} \, ,
\end{align}
using dislocation form. As for the expression of the boundary term using $\Distress$, it can be easily recovered from the relations $\bm{\mu} = \Double \vb{n}|_{\Area_N^{\Coss}}$ and $\Double = \Double(\Distress,\Coss)$. We discuss some further relations between the two forms of the balance of angular momentum in \cref{ap:obs}. 

In the following, we formulate the action-functional of the Cosserat model in dislocation form. The corresponding Eringen-type formulation can be found in \cite{Eringen1999}.
To simplify constitutive laws, we first introduce the Biot-type stress and dislocation-stress tensors
\begin{align}
    \Rpiola = \Coss^T \Piola \,, && \Rdistress = \Coss^T \Distress \, . 
\end{align}
A variational principle can thus be introduced via the hyperelastic micropolar and curvature energy densities
\begin{align}
    \Rpiola = \Grad_{\Stretch} \Psi_\mathrm{mp} (\Stretch) \, , && 
    \Rdistress = \Grad_{\Curva} \Psi_\mathrm{curv} (\Curva) \, ,
\end{align}
such that
\begin{align}
    \delta \Psi_\mathrm{mp} (\Stretch) = \con{\delta \Stretch}{\Grad_{\Stretch} \Psi_\mathrm{mp}} = \con{\delta\Stretch}{\Rpiola} \, , && \delta \Psi_\mathrm{curv}(\Curva) =  \con{\delta \Curva}{\Grad_{\Curva} \Psi_\mathrm{curv}} = \con{\delta \Curva}{\Rdistress} \, .
\end{align}
For a physically linear material, a general non-negative isotropic function of the Biot-type stretch tensor can be defined as
\begin{align}
    \Psi_\mathrm{mp}(\Stretch) = \underbrace{\mu \norm{\sym(\Stretch - \one)}^2 + \dfrac{\lambda}{2} \tr(\Stretch-\one)^2}_{\frac{1}{2}\norm{\sym(\Stretch - \one)}_{\Cm}^2} + \muc\norm{\skw(\Stretch - \one)}^2 = \frac{1}{2}\norm{\sym(\Stretch - \one)}_{\Cm}^2 + \muc\norm{\skw\Stretch}^2 \, ,
\end{align}
such that $\Psi_\mathrm{mp}(\one) = 0\, ,\quad\Grad_{\Stretch}\Psi_\mathrm{mp}|_{\one} = 0$, and the material tensor reads
\begin{align}
    \Cm = 2 \mu \J + \lambda \one \otimes \one  \in \R^{3 \times 3 \times 3 \times 3} \, .   
\end{align}
For the dislocation density we postulate the isotropic energy function to be 
\begin{align}
    \Psi_\mathrm{curv}(\Curva) &= \dfrac{\mu \Lc^2}{2} (\alpha_1 \norm{\dev \sym \Curva}^2 + \alpha_2 \norm{\skw \Curva}^2 + \alpha_3 \norm{\sph \Curva}^2) = \dfrac{\mu \Lc^2}{2} \norm{\Curva}^2_{\Lm} \, ,
\end{align}
where $\Lc \geq 0$ is the characteristic length-scale parameter that governs the energy of higher order effects, and the tensor of dimensionless weights takes the form
\begin{align}
    &\Lm = \alpha_1 \mathbb{D} \mathbb{S} + \alpha_2 \mathbb{A} + \alpha_3 \mathbb{V} \in \R^{3 \times 3\times 3\times 3} \, .
\end{align}
A deduction of the constant identity or any constant-valued rotation is not needed here, since, e.g., for $\Coss = \one$ one already finds $\Curl \Coss = \Curl \one = 0$. 
Material coefficients for $\Cm$ and $\Lm$ can be found in \cite{Lakes2018_stability,Lakes2021_softening,Lakes1987_negative_poisson,Lakes1995_experimental_methods,Lakes1985_microelasticity,YangLakes1981_transient} for the classical Eringen form, and in \cite{Ghiba} for dislocation form.
Now, putting it all together, the internal energy of the Cosserat model reads
\begin{align}
    \boxed{
    \mathcal{I}(\defmap,\Coss) = 
    \dfrac{1}{2}\int_{\Vol} \norm{\sym(\Stretch - \one)}^2_{\Cm} + 2\muc \norm{\skw\Stretch}^2 + \mu \Lc^2 \norm{\Curva}_{\Lm}^2  \, \dd \Vol \, .
    }
\end{align}
By introducing the couple-forces tensor 
\begin{align}
    \bm{M} = \dfrac{1}{2} (\Anti\vb{m})\Coss \quad \iff \quad \vb{m} = 2 \axl(\bm{M}\Coss^T) \, , 
\end{align}
the external work can be written as
\begin{align}
    \boxed{
    \mathcal{W}(\defmap,\Coss) = \int_{\Vol} \con{\defmap}{\vb{f}}  + \con{\Coss}{\bm{M}} \, \dd \Vol \, .
    } 
\end{align}
Thus, the solution is a minimiser of the corresponding action-functional
\begin{align}
    \boxed{
    \mathcal{A}(\defmap,\Coss) = \mathcal{I}(\defmap,\Coss) - \mathcal{W}(\defmap,\Coss) \quad \to \quad \min \quad  \wrt \quad \{\defmap,\Coss\} \,,
    }
\end{align}
for which the proof of existence of minimisers as per 
\AS{
\begin{align}
    \{\defmap,\Coss \} \in [\R^3\otimes\Hone(\Vol)] \times \Hone[\Vol;\SO(3)] \, , && \begin{aligned}
        \Hone(\Vol) &= \{u \in \Le(\Vol) \;|\; \nabla u \in \Le(\Vol) \otimes \R^3\} \, ,  \\  \Le(\Vol) &= \left\{u:\Vol \to \R \;|\; \int_\Vol u^2 \, \dd \Vol < + \infty\right\}
    \end{aligned}  \, ,
\end{align}
}
can be found in \cite{Neff2004,Neff2005,Neff2006b,Neff2014}.

\subsection{The finite elasticity and couple-stress limits}
The Cosserat micropolar model reduces to finite elasticity or couple-stress theory depending on the choice of material coefficients. Let the characteristic length scale parameter approach zero $\Lc \to 0$, while the Cosserat couple modulus is larger or equal to the elastic shear modulus $\muc \geq \mu$ and no couple-forces or couple-tractions are prescribed, then at the limit $\Lc = 0$, the action-functional simplifies to finite elasticity with the Biot stretch tensor \cite{Neff2004,Neff2006b}
\begin{align}
    \Coss \to \polar \defgrad \, , && \Stretch \to \bm{U} \, , && \Energy(\defmap,\Coss) \to \Energy_\mathrm{e}(\defmap) =  \dfrac{1}{2} \int_\Vol \norm{\sym(\bm{U} - \one)}_{\Cm}^2 \, \dd \Vol \, ,
\end{align}
\AS{where the polar operator is given by the polar decomposition
\begin{align}
    \polar \defgrad = \bm{R} \, , && \bm{F} = \bm{R}\bm{U} \quad \Rightarrow \quad \bm{R} = \bm{F} \bm{U}^{-1} \, , && \bm{R}:\Vol \to \SO(3) \, ,&& \bm{U}:\Vol \to \Sym^{++}(3) \, .
\end{align}}
This \AS{regression} happens because the curvature vanishes, leaving $\Coss$ as a local variable, which then collapses onto $\polar \defgrad$ \cite{Neff2004,Neff2006b,Lankeit2016,Neff2006Couple,Fischle2017}.
On the other hand, let $\muc \to + \infty$ with $\Lc > 0$, then $\Coss \to \polar \defgrad$ is inherent for finite coupling  energies
\begin{align}
    \muc \norm{\skw(\Coss^T \defgrad)}^2 = \dfrac{\muc}{4} \norm{\Coss^T \defgrad -  \defgrad^T\Coss}^2 = \dfrac{\muc}{4} \norm{\Coss^T (\polar \defgrad) \bm{U} -  \bm{U}(\polar \defgrad)^T\Coss}^2 \, ,
\end{align}
since $\Coss = \polar \defgrad$ implies 
\begin{align}
    \norm{\Coss^T (\polar \defgrad) \bm{U} -  \bm{U}(\polar \defgrad)^T\Coss}^2 \to \norm{\bm{U} - \bm{U}}^2 = 0 \, .
\end{align}
Consequently, the micropolar model approaches couple-stress theory
\begin{align}
    \Energy(\defmap,\Coss) \to \Energy_\mathrm{cs}(\defmap) = \dfrac{1}{2} \int_\Vol \norm{\sym(\bm{U} - \one)}_{\Cm}^2 + \mu \Lc^2 \norm{\bm{R}^T\Curl\bm{R}}_{\Lm}^2 \, \dd \Vol \, . 
\end{align}
While both cases lead to $\Coss \to \polar \defgrad$, only the case of $\muc \to +\infty$ can lead to locking, since $\Coss$ is not local and the energies produced by deviations between $\Coss$ and $\polar\defgrad$ in $\muc \norm{\skw(\Coss^T \defgrad)}^2$ are massively scaled by $\muc$. Observably, $\Curl (\polar \Grad \defmap)$ implies second-order derivatives on components of the deformation map, such that the case $\muc \to + \infty$ represents a fourth-order limiting problem given by the action-functional
\begin{align}
    \mathcal{A}_\mathrm{cs}(\defmap) = \mathcal{I}_\mathrm{cs}(\defmap) - \mathcal{W}_\mathrm{cs}(\defmap) \quad \to \quad \min \quad  \wrt \quad \defmap \, ,
\end{align}
with the corresponding external work
\begin{align}
    \mathcal{W}_\mathrm{cs}(\defmap) = \int_{\Vol} \con{\defmap}{\vb{f}}  + \con{\bm{R}}{\bm{M}} \, \dd \Vol \, .
\end{align}
Thus, existence of minimisers for the problem requires higher regularities such as $\defmap \in \R^3 \otimes \mathit{H}^2(\Vol)$, to the point that $\Curl \bm{U}$ is well-defined \cite{Lankeit2016}.

\section{Objective and structure-preserving discretisations}

In the following we present a structure-preserving discretisation of the action-functionals that maintains objectivity up to projection effects in the coupling term. We start by discussing discretisations of rotations and objectivity, \AS{covering the Euler-Rodrigues formula, Euler matrices, geodesic interpolation, and projection-based interpolation}. For completeness and to properly introduce the newly proposed ideas, we first recall several well-known notions about interpolations of $\SO(3)$-tensors, before moving on to present the structure-preserving method developed to match the discretisation of the deformation tensor. \AS{The method builds on the relationship between Lagrange and N\'ed\'elec elements in the de Rham complex, which we present briefly as well.} Discretised fields given by Galerkin quasi-projections are indicated by a subscript $(\cdot)_h$, and do not imperatively agree with interpolations of possibly analytical fields, for the which the interpolation operator is denoted with $\Pi^{p}(\cdot)$.

\subsection{Discretisation of rotation tensors and objectivity}

The discretisation of the Cosserat rotation tensor can be rather involved, since it requires a subspace of $\Hone[\Vol;\SO(3)]$, but $\SO(3)$ is not a vector space. Traditionally, rotation tensors were discretised by defining an axial vector \AS{in the $C^0(\Vol)$-continuous Lagrange space} $\bm{\theta}_h \in \R^3 \otimes \CG^p(\Vol) \subset \R^3 \otimes \Hone(\Vol)$ and mapping it to a rotation matrix via the Euler--Rodrigues formula 
\begin{align}
    \Coss_h(\bm{\theta}_h) = \exp \Anti \bm{\theta}_h = \one + \dfrac{\sin\norm{\bm{\theta}_h}}{\norm{\bm{\theta}_h}} \Anti \bm{\theta}_h + \dfrac{1-\cos \norm{\bm{\theta}_h}}{\norm{\bm{\theta}_h}^2} (\Anti \bm{\theta}_h)^2 \, , && \Anti \bm{\theta}_h :\Vol \to \so(3) \, ,
\end{align}
which is a closed-form solution of the exponential map $\exp:\so(3) \to \SO(3)$.
True rotations are updated multiplicatively rather than additively, where a total multiplicative update guarantees path-independence \cite{Jeleni1999}, cf. \cref{ap:mul}.
Thus, in this approach rotations are updated at the nodes before extracting the corresponding axial vector with 
\begin{align}
    \bm{\theta}_h =  \axl\log \Coss_h = \left \{ \begin{aligned}
          &\dfrac{\arccos \left (\dfrac{\tr \Coss - 1}{2} \right )}{\norm{\axl \Coss}} \axl \Coss && \text{for} &  \norm{\axl \Coss} &> 0 \\[1\baselineskip]
          &0  && \text{for} & \tr \Coss &= 3  \\[1\baselineskip]
          &   \pm \dfrac{\pi}{\norm{\Coss \vb{e}_{\xi} + \vb{e}_\xi}} (\Coss \vb{e}_{\xi} + \vb{e}_\xi)  && \text{for} &  \tr \Coss &= -1
    \end{aligned} \right . \, , &&
    \xi = \argmax_{\zeta \in \{1,2,3\}  }  \norm{\Coss \vb{e}_{\zeta} + \vb{e}_\zeta} \, ,
\end{align} 
which is a closed-form formula of the logarithmic map $\log : \SO(3) \to \so(3)$, cf. \cite{Lankeit2014} and \cref{ap:log}. The updated axial vector is then interpolated in $\R^3$.
However, this is non-objective \cite{Crisfield1999,Bauchau2013,Romero2004,Choi2024}, as rigid-body rotations are not maintained in the interpolation of the axial vector. To illustrate, let $\Coss_1$ and $\Coss_2$ be two nodal values with corresponding nodal functions $\phi_1 \in \CG^1(\Vol)$ and $\phi_2\in \CG^1(\Vol)$ \AS{of the Lagrange interpolant $\Pi_g^1$}, and impose the rigid body rotation $\bm{Q}$, then a linear interpolation yields  
\begin{align}
    \exp [\Anti\Pi_g^1(\bm{\theta})] = \exp [ \phi_1 \log (\bm{Q} \Coss_1) + \phi_2 \log( \bm{Q}\Coss_2) ] \neq \bm{Q} \exp [ \phi_1 \log (\Coss_1) + \phi_2 \log(\Coss_2) ]  \,,
\end{align}
since the matrix logarithm can satisfy $\log(\bm{Q} \Coss) = \log \bm{Q} + \log\Coss$ and the exponential map can satisfy $\exp(\log \bm{Q} + \log\Coss) = \exp(\log \bm{Q}) \exp( \log\Coss)$ \textbf{if and only if} the product of the two rotation tensors commutes $\bm{Q}\Coss = \Coss \bm{Q}$, which is generally \textbf{not} the case, cf. \cref{ap:mexp}. 

Another natural appearing approach is given by Euler matrices. Let the components of the axial vector $\bm{\theta}_h = \theta_i \vb{e}_i$ represent rotation angles around the Cartesian axes, where the rotation matrix around each axis reads
\begin{align}
    \Coss_\xi (\theta_\xi) = \vb{e}_\xi \otimes \vb{e}_\xi + (\cos\theta_\xi)(\one - \vb{e}_\xi \otimes \vb{e}_\xi) + (\sin \theta_\xi)\Anti \vb{e}_\xi \,,
\end{align}
with \textbf{no summation} over repeating $\xi$-indices, then the total rotation tensor is defined via
\begin{align}
    \Coss_h (\bm{\theta}_h) = \Coss_z(\theta_z) \Coss_y(\theta_y) \Coss_x(\theta_x) \, . 
\end{align}
However, the same problem of non-objectivity applies \cite{Bauchau2013} and Euler matrices are also known to be susceptible to \AS{the Euler angle singularity \cite{Hemingway2018}, otherwise known as} gimbal locking. Pertinently, since the Cosserat model measures strains via the discrete Biot-type stretch tensor $\Stretch_h = \Coss_h \defgrad_h$ and curvature via the discrete micro-dislocation $\Curva_h = \Coss_h \Curl \Coss_h$, both of the previous choices lead to the introduction of erroneous strains and curvature energies in the presence of rigid body rotations. 

Recognising that the problem lies not in the representation of the rotation tensor as the exponential map of some axial vector $\Coss_h = \exp \Anti \bm{\theta}_h$, but specifically in the interpolation procedure, it was suggested to replace the interpolation with one that respects the true manifold of rotations \cite{Crisfield1999}. This separate definition of nodal values $\Coss_i$ and a geodesic interpolation operator $\Pi_s^p$ is given by geodesic finite elements \cite{Sander2009,Sander2012,Sander2015,Greco2023,Greco2024}, where the interpolation is defined via
\begin{align}
    \Coss_h \simeq \Pi_s^p \Coss = \argmin_{\Coss \in \SO(3)} \dfrac{1}{2} \sum_i \phi_i \norm{  \log( \Coss_i \Coss^T) }^2 \, , && \phi_i \in \CG^p(\Vol) \, , 
\end{align}
which requires the solution of a minimisation problem. For a two-nodes linear interpolation on a curve $\curv$, this can be rewritten as
\begin{align}
    \Pi_s^1\Coss = \Coss_1\exp[ \phi_2 \log(\Coss_1^T \Coss_2) ] \, , && \phi_2 \in \CG^1(\curv) \, , 
\end{align}
which can be easily shown to be objective
\begin{align}
    \Pi_s^1 (\bm{Q}\Coss) = \bm{Q}\Coss_1\exp[ \phi_2 \log([\bm{Q}\Coss_1]^T [\bm{Q}\Coss_2]) ] = \bm{Q}\Coss_1\exp[ \phi_2 \log(\Coss_1^T\bm{Q}^T \bm{Q}\Coss_2) ] = \bm{Q} \Coss_1\exp[ \phi_2 \log(\Coss_1^T \Coss_2) ]  \, .
\end{align}
Geodesic interpolation is inherently the only one that respects the underlying physics of the problem, as it follows the true shortest path between rotations on the hypersphere $\mathrm{S}^4$, and does not introduce any artifacts into the interpolation. While the shortest path between two poles is not unique, this generally does not pose a problem in finite element applications, since either $h$- or $p$-refinements allow to keep the difference between two nodal values below this threshold $\norm{\Delta \bm{\theta}} < \pi$.
More importantly, geodesic interpolation maintains rates of change between rotations. In other words, a linear or quadratic varying rotation stays as such under geodesic interpolation. This is critical, as it implies that the curvature measure $\Curva_h$ maintains its physical validity with this approach. Nevertheless, for the purpose of structure-preservation within this work, this interpolation alone is not sufficient. 

Another approach that maintains the nodal values but interpolates differently is projection-based interpolation \cite{Grohs2019}. Its structure is relatively similar to geodesic interpolation, but the distance-measure is changed to the Euclidean norm
\begin{align}
    \Coss_h \simeq \Pi_r^p \Coss = \argmin_{\Coss \in \SO(3)}\dfrac{1}{2} \sum_i \phi_i \norm{ \Coss_i - \Coss}^2 = \argmin_{\Coss \in \SO(3)}\dfrac{1}{2} \norm{ \Pi_g^p \Coss - \Coss}^2 \, , && \phi_i \in \CG^p(\Vol) \, .  
\end{align}
In essence, the procedure embeds $\Coss_i \in \SO(3)$ into the larger algebraic space of all three-by-three matrices $\R^{3 \times 3} \supset \SO(3)$, interpolates it there, and then projects the result back into $\SO(3)$. 
As it turns out, this projection amounts to an extraction of the polar decomposition, cf. \cref{ap:pol}. As such, it can be rewritten as
\begin{align}
    \Pi_r^p\Coss = \polar (\Pi_g^p \Coss) \, .
\end{align}
By the polar decomposition, this approach is also objective 
\begin{align}
    \Pi_r^1 (\bm{Q}\Coss) = \polar(\underbrace{\bm{Q}[\phi_1 \Coss_1 + \phi_2  \Coss_2]}_{\Pi_g^1(\bm{Q}\Coss)}) = \bm{Q} \polar(\phi_1 \Coss_1 + \phi_2  \Coss_2) = \bm{Q} \polar(\Pi_g^1 \Coss) \, , && \phi_i \in \CG^1(\Vol) \, , 
\end{align}
since for every arbitrarily evaluated $\bm{Z} = \Pi_g^p\Coss : \Vol \to \R^{3 \times 3}$ we find 
\begin{align}
    \underbrace{\bm{Q} \bm{Z}}_{\widetilde{\bm{Z}}}  = \underbrace{\bm{Q} \Coss}_{\widetilde{\bm{R}}} \bm{U} \quad \Rightarrow \quad \polar(\underbrace{\bm{Q} \bm{Z}}_{\widetilde{\bm{Z}}}) = \underbrace{\bm{Q} \Coss}_{\widetilde{\bm{R}}} = \bm{Q} \polar \bm{Z} \, , &&  \bm{U} \in \Sym(3) \, . 
\end{align}
However, distances measured over the Euclidean norm do not generally agree with geodesic distances. Consequently, this interpolation does not maintain rates of changes of $\Coss$, implying that curvature measures, such as $\Curva_h$, are distorted. Nevertheless, this interpolation procedure is critical for the structure-preserving approach of this work.   

\subsection{Discretisation of the deformation map and structure-preservation}
The discretisation of the deformation map $\defmap$ is classically accomplished via the standard $\C^0(\Vol)$-continuous Lagrange finite element space
\begin{align}
    \defmap_h \in \R^3 \otimes \CG^p(\Vol)   \subset \R^3 \otimes \Hone(\Vol)  \, .
\end{align}
By virtue of the de Rham complex \cite{arnold_complexes_2021,Pauly2022}
\begin{align}
    \Hone(\Vol) \supset \CG^p(\Vol) \xrightarrow[]{\nabla} \Ned_{II}^{p-1}(\Vol) \subset \Hcurl{,\Vol} \, ,
\end{align}
the discrete deformation gradient satisfies
\begin{align}
    \defgrad_h = \Grad \defmap_h \in \Grad [\R^3 \otimes \CG^p(\Vol)] \subset \R^3 \otimes \Ned_{II}^{p-1}(\Vol) \subset \R^3 \otimes \Hcurl{,\Vol} \, ,
\end{align}
\AS{where the Ciarlet definition reads
\begin{align}
    \Ned_{II}^1(\Vol) = \{T_k\subset \Vol,\, \P^1(T_k)\otimes\R^3,\,\spa(D_{ij})\} \, , && \Vol = \bigcup_k T_k \,, && D_{ij}(\cdot) = \int_{\curv_{ij}}q \con{\cdot}{\vb{t}} \, \dd \curv \quad \forall \, q \in \P^1(\curv_{ij}) \, ,
\end{align}
with $T_k\subset \R^3$ representing a simplex in the mesh of $\Vol$ and $\curv_{ij}$ being its respective boundary curves of codimension two.} In other words, the infinite-dimensional deformation gradient lives in $\R^3 \otimes \Hcurl{,\Vol}$, and the discrete deformation gradient lives in the corresponding finite dimensional N\'ed\'elec space $\R^3 \otimes \Ned_{II}^{p-1}(\Vol)$. A key characteristic of the N\'ed\'elec space \cite{nedelec_mixed_1980,ndlec_new_1986,sky_higher_2023} is its lower regularity in comparison to the $\C^0(\Vol)$-continuous Lagrange space \cite{SKY2024104155}. Specifically, tensorial fields in $\Ned_{II}^{p-1}(\Vol)$ are only tangentially continuous $\jump{\defgrad_h \Anti \vb{n}}_{\Xi} = 0$, and their normal components $\defgrad_h\vb{n}$ are allowed to jump at material interfaces $\Xi$, as illustrated in \cref{fig:inter}.
\begin{figure}
    \centering
    \input{figs/inter}
    \caption{The deformation tensor is defined such that $\defgrad = \defgrad_1$ in $\Vol_1$ and $\defgrad = \defgrad_2$ in $\Vol_2$ with $\Vol = \Vol_1 \cup \Vol_2$. Its product with the normal vector $\vb{n}$ at the arbitrary interface $\Xi = \Vol_1 \cap \Vol_2$ defines two different vectors, depending on whether $\defgrad_1$ or $\defgrad_2$ is evaluated. Nevertheless, the deformation tensor satisfies $\defgrad \in \R^3 \otimes \Hcurl{,\Vol}$ as long its tangential components $\defgrad \Anti \vb{n}$ are continuous at the interface, which is naturally given for any $\defgrad = \Grad \defmap$ with the deformation map $\defmap \in \R^3 \otimes \Hone(\Vol)$.}
    \label{fig:inter}
\end{figure}
Inherently, this implies that less control is applied over fields in $\Ned_{II}^{p-1}(\Vol)$, whose degrees of freedom correspond to scaled edge-wise and face-wise tangential projections at element interfaces, governing entire vectors. 

We can now state the challenge of structure-preservation. That is, let $\muc \to + \infty$, then the discretisation must satisfy $\Coss_h \to \polar \bm{F}_h$ in the weak sense. Evidently, if $\Coss_h$ is interpolated using geodesic elements $\Pi_s^p \Coss$, then it could never match $\polar\defgrad_h = \polar \Grad\defmap_h$ exactly, since the former is essentially given by a projection-based interpolation $\polar\defgrad_h \simeq \polar(\Grad\Pi_g^p\defmap)$ of lower  regularity. Now, assuming commuting interpolants over the de Rham complex \AS{\cite{arnold_complexes_2021,Pauly2022} as per \cite{Demkowicz2000}}
\begin{align}
    \begin{matrix}
        \dots &\xrightarrow[]{\id} &\Hone(\Vol) &\xrightarrow[]{\nabla} & \Hcurl{,\Vol} &\xrightarrow[]{\curl}&\dots \\[0.5em]
        &&\Pi_g^p \bigg \downarrow & & \Pi_c^{p-1} \bigg \downarrow \\[1em]
        \dots&\xrightarrow[]{\id}&\CG^p(\Vol) &\xrightarrow[]{\nabla} & \Ned_{II}^{p-1}(\Vol) &\xrightarrow[]{\curl}&\dots 
    \end{matrix} \, , 
\end{align}
then $\polar\defgrad_h \simeq \polar(\Grad\Pi_g^p\defmap) = \polar(\Pi_c^{p-1}\Grad\defmap)$, \AS{where $\Pi_c^{p-1}$ is the N\'ed\'elec interpolant given by its functional degrees of freedom from the definition of the space. The latter commutation} provides the key to solving this problem. Namely, we propose to interpolate $\Coss_h$ into the N\'ed\'elec space and then project the result back onto $\SO(3)$ via the polar decomposition. 
The interpolation reduces the regularity of the rotation tensor $\Pi_c^{p-1}\Coss_h \in \R^3 \otimes \Ned_{II}^{p-1}(\Vol)$ to that of the deformation tensor $\defgrad_h \in \R^3 \otimes \Ned_{II}^{p-1}(\Vol)$ but introduces coupling artifacts since the result is no longer a true rotation tensor $\Pi_c^{p-1}\Coss_h: \Vol \to \R^{3 \times 3}$. Its rotation form is subsequently restored via the polar projection $\polar(\Pi_c^{p-1}\Coss_h):\Vol \to \SO(3)$.
In other words, the discrete Biot-type stretch tensor is modified to 
\begin{align}
    \Stretch_h = \Coss_h^T \defgrad_h \simeq  (\polar\Pi_c^{p-1}\Coss)^T \Grad \Pi_g^p \defmap &= (\polar\Pi_c^{p-1}\Coss)^T (\polar\Grad \Pi_g^p \defmap) \bm{U}_h 
    \notag \\
    &= (\polar\Pi_c^{p-1}\Coss)^T (\polar\Pi_c^{p-1}\Grad \defmap) \bm{U}_h \, .  
\end{align}
Intuitively, since $\Coss$ and $\Grad \defmap$ are interpolated in the same space, it should now be possible to construct a Cosserat rotation tensor $\Coss_h$ that matches $\polar \defgrad_h$. Although the $\polar$ operator is not linear, this expectation can be further explained by construction by observing that, in general, for some $\bm{Y}$ and $\bm{Z}$ one finds
\begin{align}
    \polar([\Pi_c^{p-1}\bm{Z}]^{-1} \Pi_c^{p-1}\bm{Y}) &= \polar([\Pi_c^{p-1}\bm{Z}]^{-1}) \polar (\Pi_c^{p-1}\bm{Y}) = \polar(\Pi_c^{p-1}\bm{Z})^{-1} \polar (\Pi_c^{p-1}\bm{Y}) 
    \notag \\
    &= \polar(\Pi_c^{p-1}\bm{Z})^{T} \polar (\Pi_c^{p-1}\bm{Y}) \, ,
\end{align}
since inversion and transposition are synonymous on $\SO(3)$ and the polar operator satisfies
\begin{align}
    \bm{R} &= \polar(\underbrace{\bm{R}\bm{U}}_{\bm{Z}})  \, , \qquad  \bm{R}^T = \polar(\underbrace{\bm{U}\bm{R}^T}_{\bm{Z}^T}) = \polar(\bm{R}^T \bm{V}\bm{R}\bm{R}^T) = \polar(\bm{R}^T \bm{V}) \, , \qquad \Sym(3) \ni \bm{V} =  \bm{R}\bm{U}\bm{R}^T \, , 
    \notag \\
    \bm{R}^T &= \polar( \underbrace{\bm{U}^{-1} \bm{R}^{T}}_{\bm{Z}^{-1}} ) = \polar([\bm{R}^T \bm{V}\bm{R}]^{-1} \bm{R}^T) = \polar(\bm{R}^T \bm{V}^{-1}\bm{R} \bm{R}^T) = \polar(\bm{R}^T \bm{V}^{-1}) 
     \, .\qquad
\end{align}
Evidently, let $\bm{Z}_h \simeq \Pi_c^{p-1}\bm{Z}$ and $\bm{Y}_h\simeq \Pi_c^{p-1}\bm{Y}$, then the two are spanned by the exact same polynomial space with the same regularity, such that it suffices to set $\bm{Z}_h = \bm{Y}_h$ to get 
\begin{align}
    \polar(\bm{Z}_h^{-1} \bm{Z}_h) = (\polar\bm{Z}_h)^T (\polar \bm{Z}_h) = \one \, .
\end{align}
In short, what we propose is a mixed-interpolation strategy for $\Coss$. Namely, a projection-based interpolation via N\'ed\'elec elements $\polar \Pi_c^{p-1}\Coss_h$ for interactions with $\defgrad_h$
\begin{align}
    \boxed{
    \begin{aligned}
&\norm{\sym(\overbrace{\Coss_h^T\defgrad_h}^{\Stretch_h} - \one)}_\Cm^2 +  2\muc \norm{ \skw (\overbrace{\Coss_h^T\defgrad_h}^{\Stretch_h} ) }^2 
\\[1\baselineskip]
&\quad\to \quad \norm{\sym([ \polar\Pi_c^{p-1}\Coss_h]^T \defgrad_h - \one)}_\Cm^2 + 2\muc \norm{ \skw ([ \polar\Pi_c^{p-1}\Coss_h]^T\defgrad_h ) }^2  \, .
    \end{aligned}
    }
\end{align}
and geodesic elements for the discretisation of $\Coss$ and the computation of curvature measures
\begin{align}
    \boxed{
    \Coss_h \in \GEO^q(\Vol) \subset \Hone[\Vol;\SO(3)] \, ,
    } 
\end{align}
where $\GEO^q(\Vol)$ is defined by the basis functions of the Lagrange space $\CG^q(\Vol)$ and the geodesic interpolation operator $\Pi_s^q$. 
Notwithstanding, it remains to define the underlying $q$-polynomial order of the geodesic elements, such that the computation remains robust. Further, at this stage it is also important to emphasize that the choice of interpolation operator is not arbitrary, but rather motivated by the de Rham complex. Namely, one could instead consider replacing $\Pi_c$ with an interpolant of even lower regularity. However, this may result in a loss of control over certain coupled modes and consequently lead to hour-glassing, cf. \cite{WRIGGERS1996201,Pfefferkorn,Bieber} and the construction within \cite{Neunteufel2021}. 

The here presented approach draws its inspiration from the the \textit{Mixed Interpolated Tensorial Components} method, i.e., the \textit{MITC} elements \cite{Bathe1986}, and can be viewed as an extension of this approach to more general three-dimensional continua and nonlinear finite element spaces. Namely, the MITC elements target the two-dimensional manifolds of plates and shells, alleviating shear-locking of an additive coupling over two linear spaces. 
In the formal analysis of the MITC elements, their stability is given by a commuting interpolation operator and stable Stokes-pairs \cite{Stenberg1997,Brezzi1989}. In the original intuitive view \cite{Bathe1986}, the stability of MITC elements stems from their construction to match shear-modes on the triangulation and eliminate any non-matching modes by separation. The number of possible modes is fundamentally tied to the degrees of freedom and can be estimated element-wise by the dimension of the space. 
In this work, we draw on aspects of both perspectives and combine them to propose two possible choices of discretisation spaces, which we validate through numerical benchmarks. A formal analysis is left for future work. The first choice is by analogy to MITC6 \cite{LEE2004945}
\begin{align}
    \boxed{
    \begin{aligned}
        \defmap_h \in \R^3 \otimes \CG^2(\Vol) \, , && \Coss_h \in \GEO^2(\Vol) \,, && \Pi_c^1 : \GEO^2(\Vol) \to \R^3 \otimes\Ned_{II}^1 (\Vol) \, ,
    \end{aligned}
    }
\end{align}
where the deformation is discretised by quadratic Lagrange elements, the Cosserat rotation by quadratic geodesic elements, and the interpolation operator is given by the degrees of freedom of the linear N\'ed\'elec element of the second type. With this choice of polynomial orders, the interpolation operator satisfies the commuting property $\Grad(\Pi_g^2\defgrad) = \Pi_c^1 (\Grad\defgrad)$. The element-wise dimensions of each of the spaces read 
\begin{align}
    \dim [\R^3 \otimes \CG^2(T)] = 30 \, , && \dim \R^3 \otimes\Ned_{II}^1 (T) = 36 \, , && \dim \GEO^2(T) = 30 \, .
\end{align}
The dimension of the space of gradients on a single element can be approximated as
\begin{align}
    \dim \Grad  [\R^3 \otimes \CG^2(T)] \leq 27 \, , 
\end{align}
by recalling that $\nabla \CG^2(T) \subset \Ned_{II}^1(T) = \Ned_{I}^0(T) \bigoplus_{j \in \mathcal{J}}\nabla\mathcal{E}_j^2(T)$ where $\mathcal{E}_j^2(T)$ are the edge-wise spaces of the quadratic basis functions of the Lagrange element with the multi-index $j \in \mathcal{J} = \{(0,1)(0,2),(0,3),\dots,(2,3)\}$, as per \cite{sky_higher_2023,SKY2024104155,Sky2025}. Namely, there follows the dimensionality of solenoidal fields $\curl \Ned_{II}^1(T) = \curl \Ned_{I}^0(T) = \curl (\R^3 \oplus \R^3 \times \vb{x}) = \R^3$ with $\dim \R^3 = 3$, implying that the linear N\'ed\'elec element of the second type contains approximately nine gradients, such that three copies of it yield exactly twenty-seven gradient modes. From the latter dimension counts it is clear that $\Pi_c^1:\GEO^2(T)$ is not surjective onto $\R^3 \otimes \Ned_{II}^1(T)$, which may lead to instability for a lone interpolation. However, the subsequent projection reduces the dimension of the space substantially. In other words, the relevant space for stability considerations is $\polar [\R^3 \otimes \Ned_{II}^1(T)]$ which is never constructed explicitly. Still, its dimension can be estimated as 
\begin{align}
    \dim \polar [\R^3 \otimes \Ned_{II}^1(T)] = \dim \polar [\R^{3 \times 3} \otimes \P^1(T)] \approx 12 \, , && \dim \P^1(T) = 4 \, , 
\end{align}
considering that $\dim \SO(3) = 3$. Accordingly, we can estimate the dimension of the space of rotational gradients with $\dim \polar \Grad  [\R^3 \otimes \CG^2(T)] \approx 9$. Now, observe that the N\'ed\'elec interpolant is associated with two degrees of freedom per edge, and the quadratic Lagrange space contains one degree of freedom per edge. Thus, three copies of the N\'ed\'elec space boast six degrees of freedom per edge and the quadratic geodesic space features three degrees of freedom per edge. Consequently, it is reasonable to expect that at least three edge degrees of freedom remain after the interpolation and subsequent projection, implying that $\dim \polar[\Pi_c \GEO^2(T)] \geq 18$, which suffices to span the entire target space. To numerically assess and validate these notions, we consider the additional choice of discretisation spaces    
\begin{align}
    \boxed{
    \begin{aligned}
        \defmap_h \in \R^3 \otimes \CG^2(\Vol) \, , && \Coss_h \in \GEO^{2+}(\Vol) \,, && \Pi_c^1 : \GEO^{2}(\Vol) \to \R^3 \otimes\Ned_{II}^1 (\Vol) \, ,
    \end{aligned}
    }
\end{align}
where $\GEO^{2+}(\Vol)$ is defined over via the enriched Lagrange space
\begin{align}
    \CG^{2+}(T) = \CG^{2}(T) \oplus \mathcal{B}^4(T) \, , && \mathcal{B}^4(T) = \spa\{\lambda_0\lambda_1\lambda_2\lambda_3\} \, ,
\end{align}
which contains one additional quartic bubble function per element, e.g., given by the multiplication of the barycentric coordinates $\lambda_i$. The latter is roughly akin to MITC7 \cite{LEE2004945}, and implies that the geodesic space is now slightly larger, allowing it to accommodate additional interaction modes should they arise.

\section{Numerical examples}

In the following we explore several benchmarks to observe and quantify our structure-preserving method. For simplicity of implementation, we employ here Euler matrices to avoid the multiplicative update and the extraction of geodesics, which requires the intricate element-wise purely \textbf{algorithmic} solution of minimisation problems \cite{Sander2015}. To maintain the validity of the tests while using Euler matrices, we restrict our benchmarks to moderate changes over uniaxial rotations. Thus, the gimbal lock does not arise and the additive update becomes a reasonable approximation of the multiplicative approach, cf. \cref{ap:mul}. Moreover, all methods compared within this work employ these same Euler rotation matrices, ensuring that they do not account for any of the observed differences.
Nevertheless, it should be noted that this simplification is not advisable for more general scenarios.

For the implementation of the polar extraction we make use of the recently introduced eigenvalue decomposition formulae from \cite{Habera1,Habera2}, which provide closed-form expressions that are applicable to automated solvers of partial differential equations. 

In the following sections, physical units are omitted from notation for clarity. The units corresponding to each quantity are as follows throughout all benchmarks.
\begin{center}
    \begin{tabular}{c|c|c}
    Symbols & Quantity & Physical unit \\
    \hline
    $\vb{x}=x\vb{e}_1+y\vb{e}_2+z\vb{e}_3$ & Length/Width/Height & $\mathrm{m}$ \\
    $E,\,\mu,\,\lambda$ & Classical material parameters & $\mathrm{Pa} = \mathrm{N}\,\mathrm{m}^{-2}$\\
    $\muc$ & Cosserat couple modulus & $\mathrm{Pa} = \mathrm{N}\,\mathrm{m}^{-2}$\\
    $\Lc$ & Characteristic length scale & $\mathrm{m}$\\
    $F$ & Force & $\mathrm{N}$\\
    $M_\mathrm{t}$ & Torque & $\mathrm{N}\,\mathrm{m}$\\
    $\vb{t}$ & Traction vector & $\mathrm{Pa} = \mathrm{N}\,\mathrm{m}^{-2}$\\
    \end{tabular}
\end{center}
All the benchmarks are computed using NGSolve\footnote{www.ngsolve.org} \cite{Sch2014,Sch1997}, and are available online\footnote{www.github.com/lschek/structure\_preserving\_cosserat}.

\subsection{Rigid body rotation of an axis-symmetric cube}

\begin{figure}
    \centering
    \begin{subfigure}{0.19\textwidth}
        \centering
        \includegraphics[width=\linewidth]{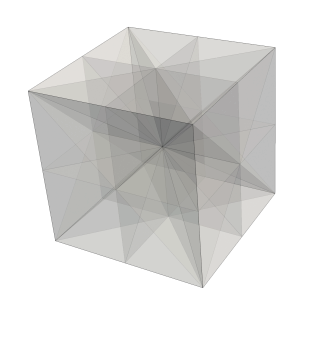}
        \caption{}
    \end{subfigure}\hfill
    \begin{subfigure}{0.19\textwidth}
        \centering
        \includegraphics[width=\linewidth]{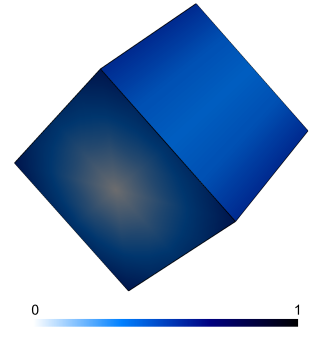}
        \caption{}
      \end{subfigure}\hfill
      \begin{subfigure}{0.19\textwidth}
        \centering
        \includegraphics[width=\linewidth]{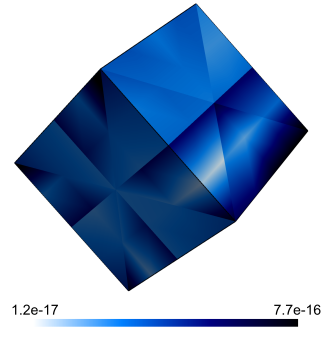}
        \caption{}
    \end{subfigure}\hfill
    \begin{subfigure}{0.19\textwidth}
        \centering
        \includegraphics[width=\linewidth]{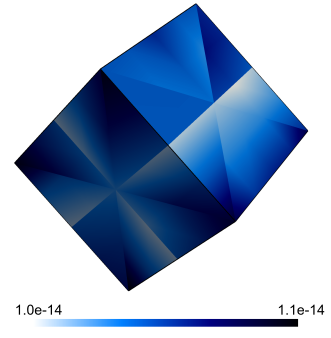}
        \caption{}
      \end{subfigure}\hfill
      \begin{subfigure}{0.19\textwidth}
        \centering
        \includegraphics[width=\linewidth]{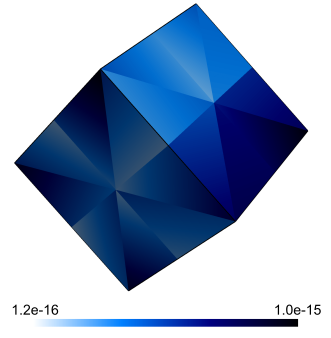}
        \caption{}
      \end{subfigure}\hfill
    \caption{\AS{Finite element mesh in the reference configuration (a) and displacement magnitude after rigid body rotation (b). Norm of the Cosserat strain, $\norm{\Coss_h^T \defgrad_h - \one}$, computed using the standard approach (c), using the interpolation $\Coss_h \mapsto \Pi_c^{p-1} \Coss_h$ (d), and using interpolation followed by polar projection, $\Coss_h \mapsto \polar(\Pi_c^{p-1}\Coss_h)$ (e).}
    }    \label{fig:unit_cube}
\end{figure}

In the first example we perform a brief consistency check with respect to superimposed rigid body motions. The domain is defined as the axis-symmetric unit cube $\overline{\Vol} = [-0.5,0.5]^3$. Its entire boundary is set to Dirichlet for the deformation field $\Area_D^{\defmap} = \partial \Vol$ and natural Neumann for the Cosserat rotation tensor $\Area_N^{\Coss} = \partial \Vol$. On the Dirichlet boundary of the deformation we prescribe $\defmap|_{\partial \Vol} = \bm{Q}\vb{x}$ with $\bm{Q} = \vb{e}_2 \otimes \vb{e}_2 + [\cos(\pi/4)](\one - \vb{e}_2 \otimes \vb{e}_2) + [\sin(\pi/4)]\Anti \vb{e}_2$, which represents a rigid body rotation of $45^\circ$ around the $y$-axis. 
Throughout the entire domain, the classical elastic Lam\'e parameters are chosen as $\lambda = 692.398$ and $\mu = 461.538$. The weights associated with the micro-dislocation are specified as $\alpha_1 \Lc^2 = 0.288$, $\alpha_2 \Lc^2 = 1.711$, and $\alpha_3 \Lc^2 = 0.42$, and the Cosserat couple modulus is taken as $\muc = \mu$.
The results shown in \cref{fig:unit_cube} depict the norm of the Cosserat strain tensor in its respective form for three approaches, which differ in the discrete Biot-type stretch tensor. Specifically, we compare the quadratic discretisation across the standard approach based on Euler matrices, a formulation using only the N\'ed\'elec interpolation, and the proposed method with the subsequent polar projection. In this simple example, all three approaches manage to reconstruct a rigid body rotation with numerically vanishing strains. Although extremely small deteriorations in the representation of zero appear, their origin may lie in the additional floating-point operations inherent to the interpolation and the projection. Even so, it is noteworthy that the polar projection reduces the error of the sole interpolation. 
\subsection{Convergence in a bending problem}

Next we consider simple bending in a thin and wide \AS{plate}. Its domain is defined as $\overline{\Vol} = [0,20]\times[0,10]\times[0,1]$, as depicted in \cref{fig:beam_dom}.
\begin{figure}
    \centering
    \input{figs/beam_dom}
    \caption{Wide cantilever \AS{plate} domain with the corresponding Dirichlet and loaded Neumann surfaces.}
    \label{fig:beam_dom}
\end{figure}
The left end of the \AS{plate} a is defined as the Dirichlet boundary for the deformation field $\Area_D^{\defmap} = \{\vb{x} \in \partial\Vol \; | \; x = 0 \}$. The remaining boundary $\Area_N^{\defmap} = \partial \Vol \setminus \Area_D^{\defmap}$ is Neumann, with tractions imposed at $x = 20$. For the Cosserat rotation tensor, the left end of the \AS{plate} is Dirichlet $\Area_{D}^{\Coss} = \Area_{D}^{\defmap}$ and the remaining boundary is natural Neumann $\Area_N^{\Coss} = \partial \Vol \setminus \Area_{D}^{\Coss}$.  
Finally, to test the behaviour and stability of the formulations we vary the Cosserat couple modulus $\muc$ with respect to the shear modulus $\muc/\mu\in\{1,10,10^2,10^3,10^4\}$. All other material parameters are defined as in the previous example.

To study convergence, we first consider a very small traction at the right end of the \AS{plate}. This amounts to a minor perturbation of the system, for which there clearly exists a unique solution. The traction is given by $\vb{t} = 0.2(\vb{e}_3-\vb{e}_1)$ implying a total force of $F = 2\sqrt{2}$ acting in a $45^\circ$ angle to the \AS{plate} axis. We start with a rough mesh of $192$ elements, and perform in-plane refinement. The error is measured relative to a high-fidelity reference solution obtained from a cubic-ordered finite element discretisation of both the displacement and the rotations on the finest mesh. The measurement itself is taken as the maximal displacement magnitude at the top corner of the loaded end. 
\begin{figure}
    \centering
  \begin{subfigure}{0.32\textwidth}
    \centering    \input{figs/bending_err_convergence_do_nothing}
    \caption{}
  \end{subfigure}\hfill
  \begin{subfigure}{0.32\textwidth}
    \centering    \input{figs/bending_err_convergence_interp_Ned2_no_polar}
    \caption{}
  \end{subfigure}\hfill
  \begin{subfigure}{0.32\textwidth}
    \centering    \input{figs/bending_err_convergence_interp_Ned2}
    \caption{}
  \end{subfigure}
    \caption{Relative error with mesh refinement (a) without interpolation, (b) with interpolation onto $\Ned_{II}^1(\Vol)$ and (c) with interpolation onto $\Ned_{II}^1(\Vol)$ and the subsequent polar extraction. Solid lines represent the formulation with $\Coss_h\in\GEO^2(\Vol)$ and dashed lines the enrichment with $\Coss_h\in\GEO^{2+}(\Vol)$.}
    \label{fig:3D_bending_convergence}
\end{figure}

The results are listed in \cref{fig:3D_bending_convergence}, where we again compare the classical approach without any interpolation or projection, an approach with only an interpolation into the N\'ed\'elec space, and our proposed method which combines the interpolation with a subsequent projection onto the Lie-group of rotations. From the first graph it becomes apparent that the classical method deteriorates for ever increasing $\muc/\mu$ ratios. Namely, we observe a drop the in convergence rates and consequently, in the absolute error, with a factor of about $10^2$ between the ratios $\mu/\muc =1$ and $\mu/\muc = 10^4$. Applying only the interpolation onto $\Ned_{II}^1(\Vol)$ in the formulation yields inconsistent results. On the one hand, the convergence rate appears to improve, while on the other hand they are unstable, producing ambiguous jumps. In contrast, the proposed methodology of interpolating into N\'ed\'elec and subsequently projecting onto $\SO(3)$ yields stable results with optimal convergence. Overall, these results already allow for several conclusions. Firstly, the classical formulation does indeed lock, as indicated by the loss of convergence rates. Secondly, the interpolation by itself behaves unpredictably. Thirdly, the proposed Γ-SPIN method is stable and optimal. More critically, it converges \textbf{in all scenarios}, thereby demonstrating that it does not represent an over-relaxation of the original problem. Finally, the fact that no fluctuations are observed for the $\GEO^2(\Vol)$ formulation and that no improved convergence rates are obtained for the enriched formulation $\GEO^{2+}(\Vol)$, marked by the dashed lines, is a credible indicator that the targeted space $\polar[\R^3 \otimes \Ned_{II}^1(\Vol)]$ is fully spanned by $\polar[\Pi_c^1 \GEO^2(\Vol)]$. 

\begin{figure}
    \centering
    \begin{subfigure}{0.32\textwidth}
        \centering
        \includegraphics[width=0.9\linewidth]{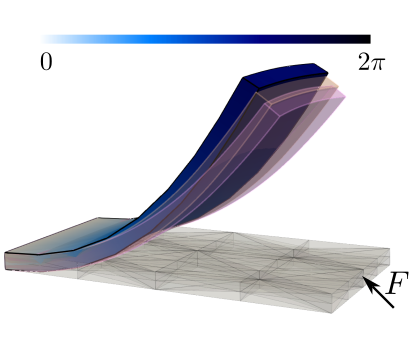}
        \caption{}
    \end{subfigure}\hfill
    \begin{subfigure}{0.32\textwidth}
        \centering \input{figs/huge_bending_convergence_muc1e0mu}
        \caption{}
      \end{subfigure}\hfill
      \begin{subfigure}{0.32\textwidth}
        \centering \input{figs/huge_bending_convergence_muc1e4mu}
        \caption{}
    \end{subfigure}  
    \caption{Deformation of the three formulations on the coarsest mesh with $\muc/\mu = 10^4$ (a). The largest deformation in blue is obtained with the proposed Γ-SPIN method for which the intensity of the rotation angle is measured. It is followed by the orange-coloured deformation with the sole interpolation and of the classical approach in violet.
    Convergence results for $\muc/\mu = 1$ (b) and $\muc/\mu = 10^4$ (c). Solid lines correspond to the formulation with $\Coss_h \in \GEO^2(\Vol)$, whereas dashed lines indicate the enriched variant with $\Coss_h \in \GEO^{2+}(\Vol)$.}    \label{fig:3D_bending_convergence_large}
\end{figure}

We repeat the study for the two extreme ratios $\mu/\muc = 1$ and $\muc/\mu = 10^4$, but increase the applied traction to yield a total force of $F = 20\sqrt{2}$, in order to induce large deformations. We use the same initial mesh as in the previous study, but this time refine it in-plane and out-of-plane. The results for this benchmark are visualised in \cref{fig:3D_bending_convergence_large}.
We immediately observe that all three formulations converge for $\muc/\mu = 1$ without loss of optimality, whereas the classical method decreases in its convergence rate for $\muc/\mu = 10^4$. In the coarsest mesh with $\muc/\mu = 10^4$, the difference in deformation is clearly visible and results in a measurable deviation. The benchmark reconfirms the emergence of locking in the classical method, revealing even poorer convergence than in the small-deformation test, with rates falling below $\mathcal{O}(h^2)$. More crucially, it again demonstrates the stable convergence of our proposed method in both extreme parameter regimes. 

\subsection{Torsion of a beam}

In the second benchmark we consider the next major curvature-dependent effect, namely torsion. We define a simple cantilever beam domain with $\overline{\Vol} = [0,8]\times[0,1]^2$, as depicted in \cref{fig:torsion_beam}.
\begin{figure}
    \centering
    \input{figs/torsion_beam}
    \caption{Long and slender cantilever beam domain with the corresponding Dirichlet boundary and loaded Neumann boundary surfaces.}
    \label{fig:torsion_beam}
\end{figure}
The division of its boundary into Dirichlet and Neumann surfaces is as in the previous example. The same holds true for the material parameters. We vary the Cosserat couple modulus $\muc$ with respect to the shear modulus between the two cases of $\muc/\mu\in\{1,10^3\}$ to observe the influence of the interpolation and projection. For the convergence analysis we start from a coarse discretisation of 48 elements, followed by uniform mesh refinements. The error is defined as the difference in the displacement magnitude at the top corner of the loaded end relative to a high-fidelity reference solution computed on the finest mesh using a cubic finite element discretisation of both the displacement and the microrotation.

\begin{figure}
    \centering
    \begin{subfigure}{0.32\textwidth}
    \centering
    \includegraphics[width=\linewidth]{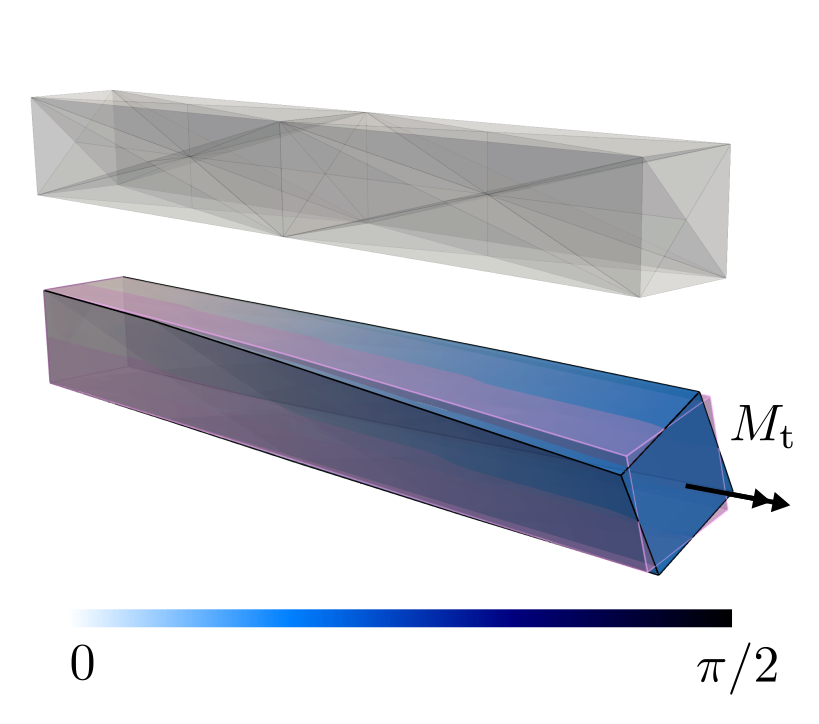}
    \caption{}
  \end{subfigure}\hfill
    \begin{subfigure}{0.32\textwidth}
        \centering
        \input{figs/torsion_convergence_muc1e0mu}
        \caption{}
      \end{subfigure}
      \begin{subfigure}{0.32\textwidth}
        \centering
        \input{figs/torsion_convergence_muc1e4mu}
        \caption{}
    \end{subfigure}
    \caption{The roughest mesh along with the corresponding deformations exhibited by all three formulations for $M_\mathrm{t}=20$ applied via couple-tractions (a) with a measurement of the rotation angle. The largest deformation exhibited by our proposed method is coloured blue, the one corresponding with the interpolation alone is in orange, and the classical approach is coloured violet. Convergence for the ratios $\muc/\mu =1$ (b), and $\muc/\mu = 10^3$ (c). Solid lines correspond to the formulation with $\Coss_h \in \GEO^2(\Vol)$, whereas dashed lines indicate the enriched variant with $\Coss_h \in \GEO^{2+}(\Vol)$.}
    \label{fig:torsion_convergence}
\end{figure}

In the first scenario, in order to this time test the coupling stemming from $\Coss_h$ over to $\defgrad_h$ we set the couple traction to $\bm{\mu} =20\vb{e}_1$ to yield a total torsion moment of $M_t = 20$ over the surface. The prescribed torsion moment induces a small deformation along the beam, for which the results are listed in \cref{fig:torsion_convergence}.  
In both extreme material parameter regimes, we observe that all formulations converge at approximately the same rate. Once more, the latter is indicative of our formulation satisfying consistency with the underlying model. 
For $\muc/\mu = 10^3$, we initially observe a reduced convergence rate of the standard approach, which is later recovered, although the overall error remains more pronounced. At the finest refinement, the relative error between the standard approach and our method differs by a factor of about $10$. On the roughest mesh, the difference in deformation is visible even for this small load case. Finally, no fluctuations for the $\GEO^2(\Vol)$-formulation nor increased convergence rates for the enriched formulation $\GEO^{2+}(\Vol)$ are observed, again indicating that the targeted space $\polar[\R^3 \otimes \Ned_{II}^1(\Vol)]$ is already spanned by $\polar[\Pi_c^1 \GEO^2(\Vol)]$.  

\begin{figure}
    \centering
  \begin{subfigure}{0.32\textwidth}
    \centering
    \includegraphics[width=\linewidth]{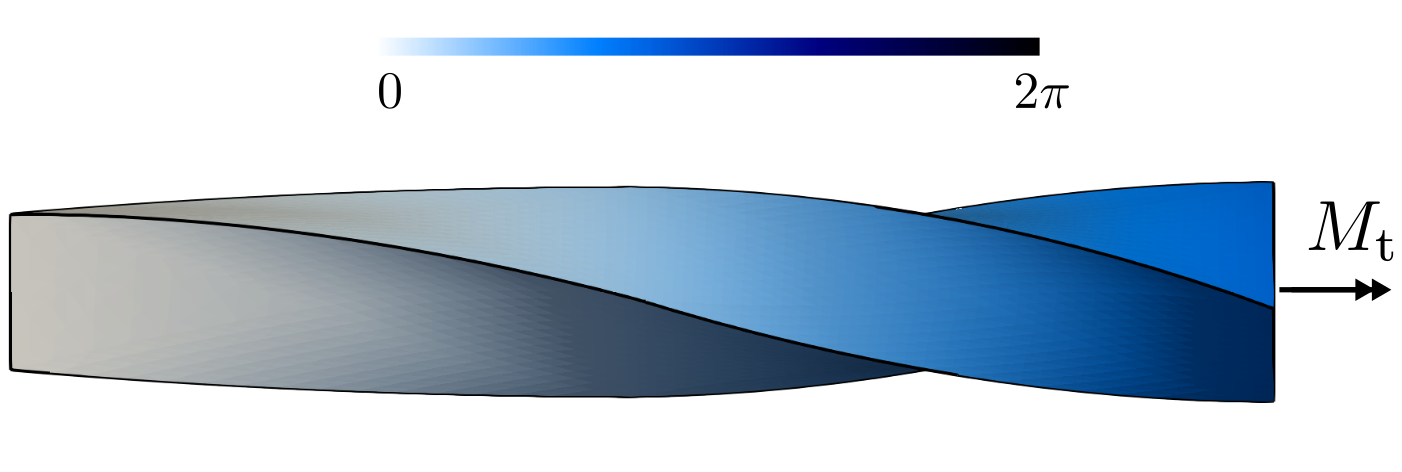}
    \caption{}
  \end{subfigure}\hfill
  \begin{subfigure}{0.32\textwidth}
    \centering
    \includegraphics[width=\linewidth]{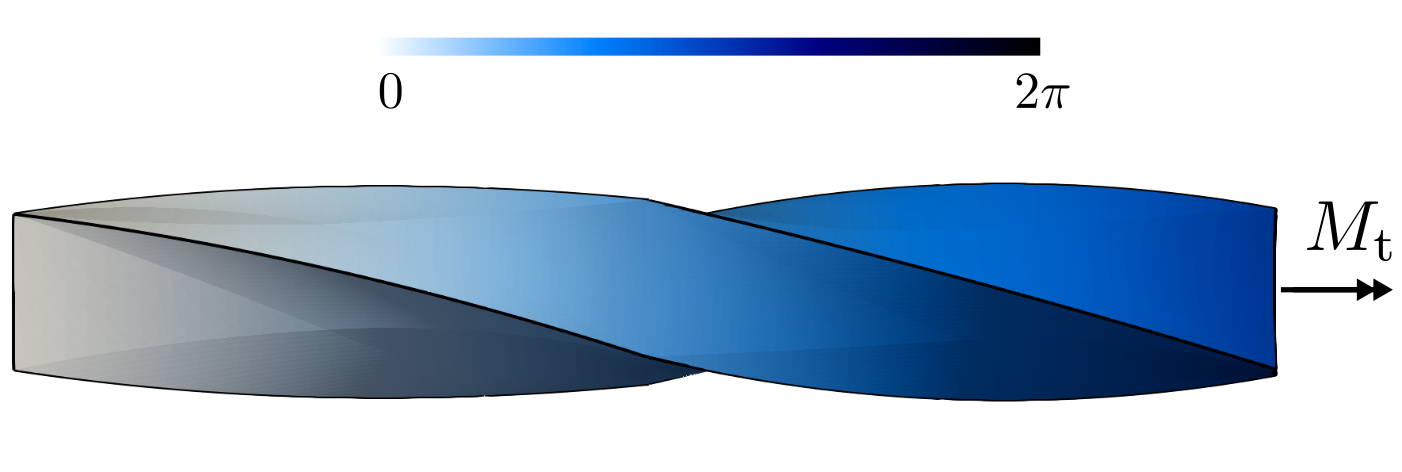}
    \caption{}
  \end{subfigure}\hfill
  \begin{subfigure}{0.32\textwidth}
    \centering
    \includegraphics[width=\linewidth]{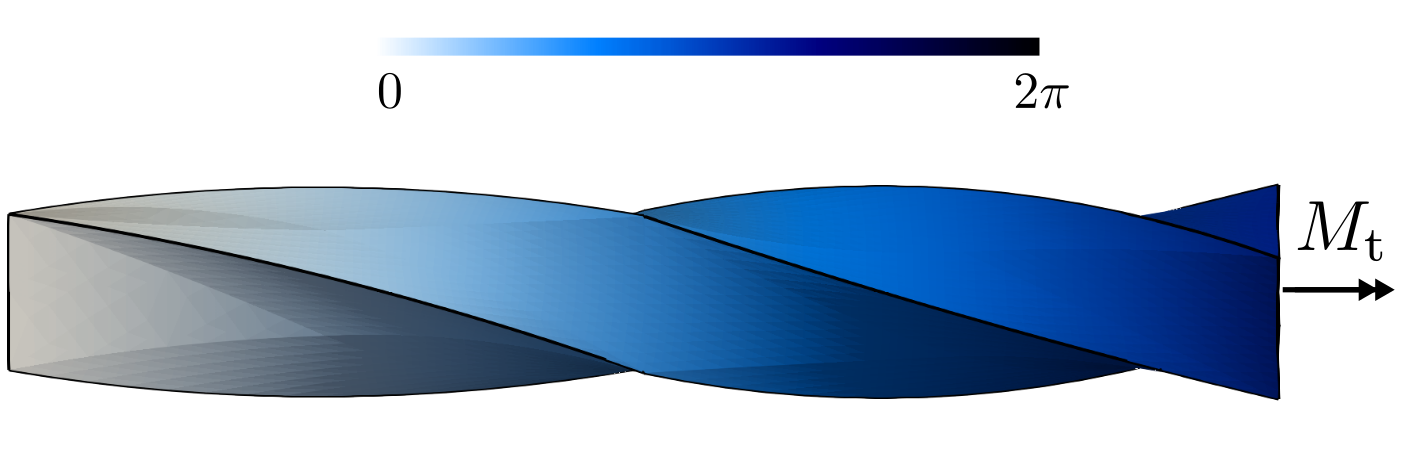}
    \caption{}
  \end{subfigure}
    \caption{Resulting deformation for $M_\mathrm{t}=400$ applied via linearly distributed tractions on the roughest mesh in the standard approach (a), with the interpolation (b), and with the interpolation and subsequent projection of the (c).}
    \label{fig:3D_torsion_comparison}
\end{figure}

Using the $\GEO^{2}(\Vol)$-type formulation, we repeat the test for an increased torsion moment totalling $M_t = 400$ on the right end surface of the reference configuration, in order to provoke large twists. In contrast to the previous scenario, the total torsion moment is now applied by setting the traction to $\vb{t} =2400[(y-0.5)\vb{e}_3-(z-0.5)\vb{e_2}]$. To provoke a type of follower-load behaviour, the corresponding traction vector is redefined as $\vb{t} \mapsto \defgrad \vb{t}$. The resulting deformations on the roughest mesh are depicted in \cref{fig:3D_torsion_comparison}.
There is clear disparity between the results of the three methods, which is exhibited by the number of twists the beam undergoes. For the standard approach we observe half a twist. The interpolation softens the formulation to the point where a full twist is triggered. Finally, the combined interpolation and projection yield one and half twists. It is important to note that these pronounced differences could not be observed in the small deformation regime, and that they clearly imply that the interpolation alone cannot fully eliminate locking. 

\subsection{\AS{Biaxial plate bending}}
\begin{figure}
    \centering
    \begin{subfigure}{0.24\textwidth}
        \centering
        \includegraphics[width=\linewidth]{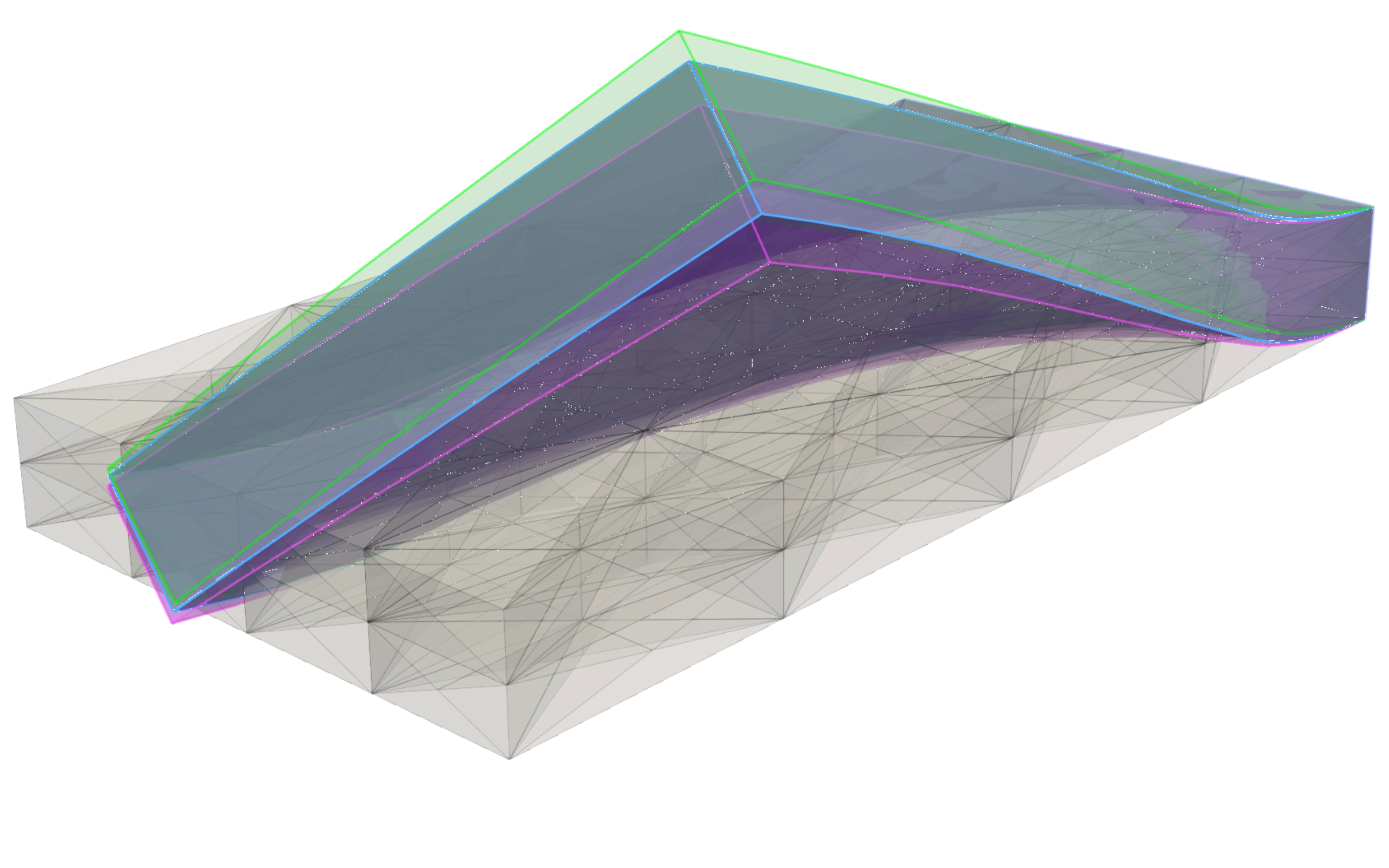}
        \caption{}
    \end{subfigure}\hfill
    \begin{subfigure}{0.24\textwidth}
        \centering
    
        \begin{tikzpicture}
        \node[anchor=south west, inner sep=0] (image) at (0,0)
            {
            \includegraphics[width=\linewidth]{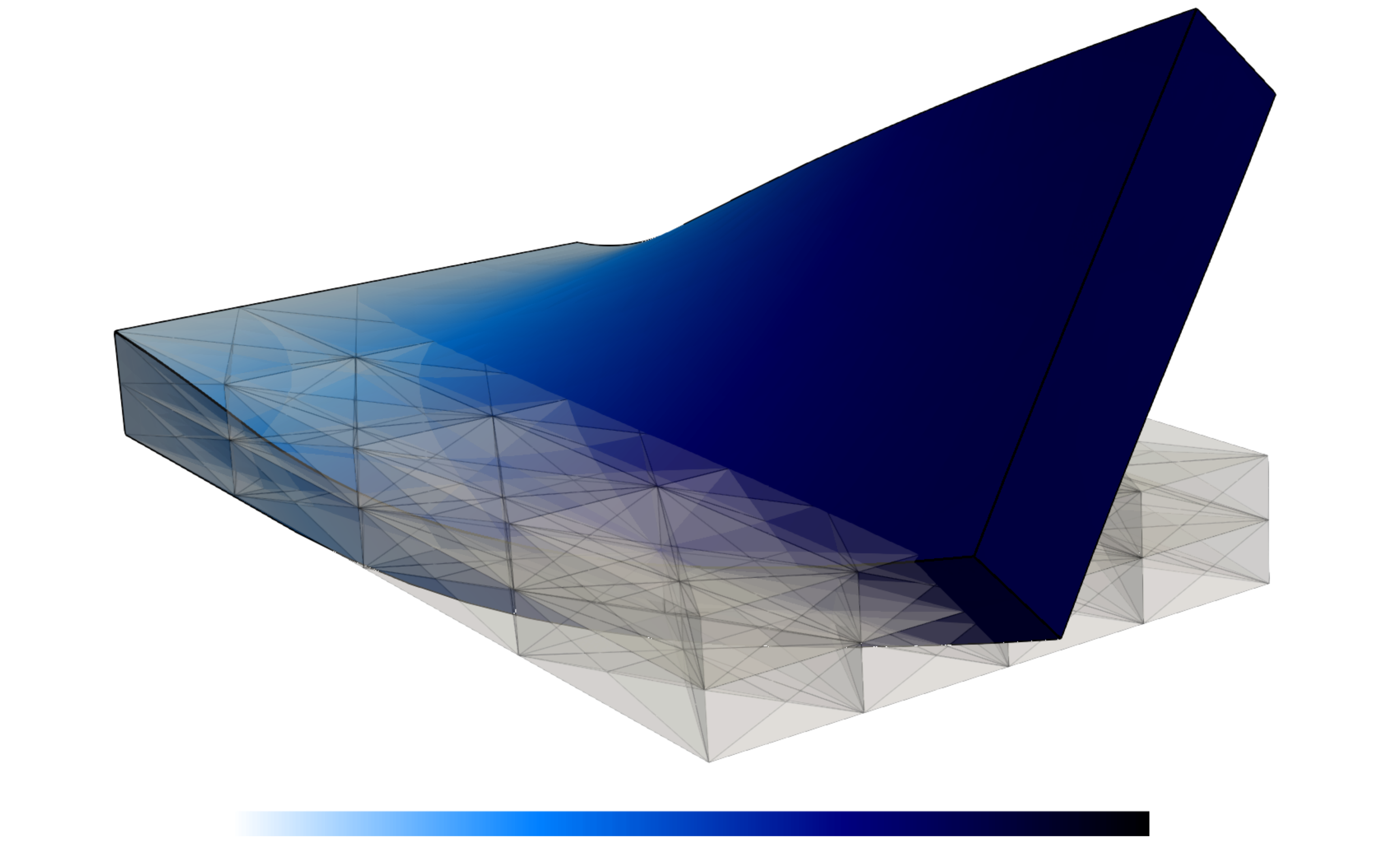}
            };
        \begin{scope}[x={(image.south east)}, y={(image.north west)}]
        \node[black, font=\small] at (0.09, 0.1) {-0.1};
        \node[black, font=\small] at (0.89, 0.1) {1};
        \end{scope}
        \end{tikzpicture}
        \caption{}
      \end{subfigure}\hfill
      \begin{subfigure}{0.24\textwidth}
        \centering
        \begin{tikzpicture}
        \node[anchor=south west, inner sep=0] (image) at (0,0)
            {
            \includegraphics[width=\linewidth]{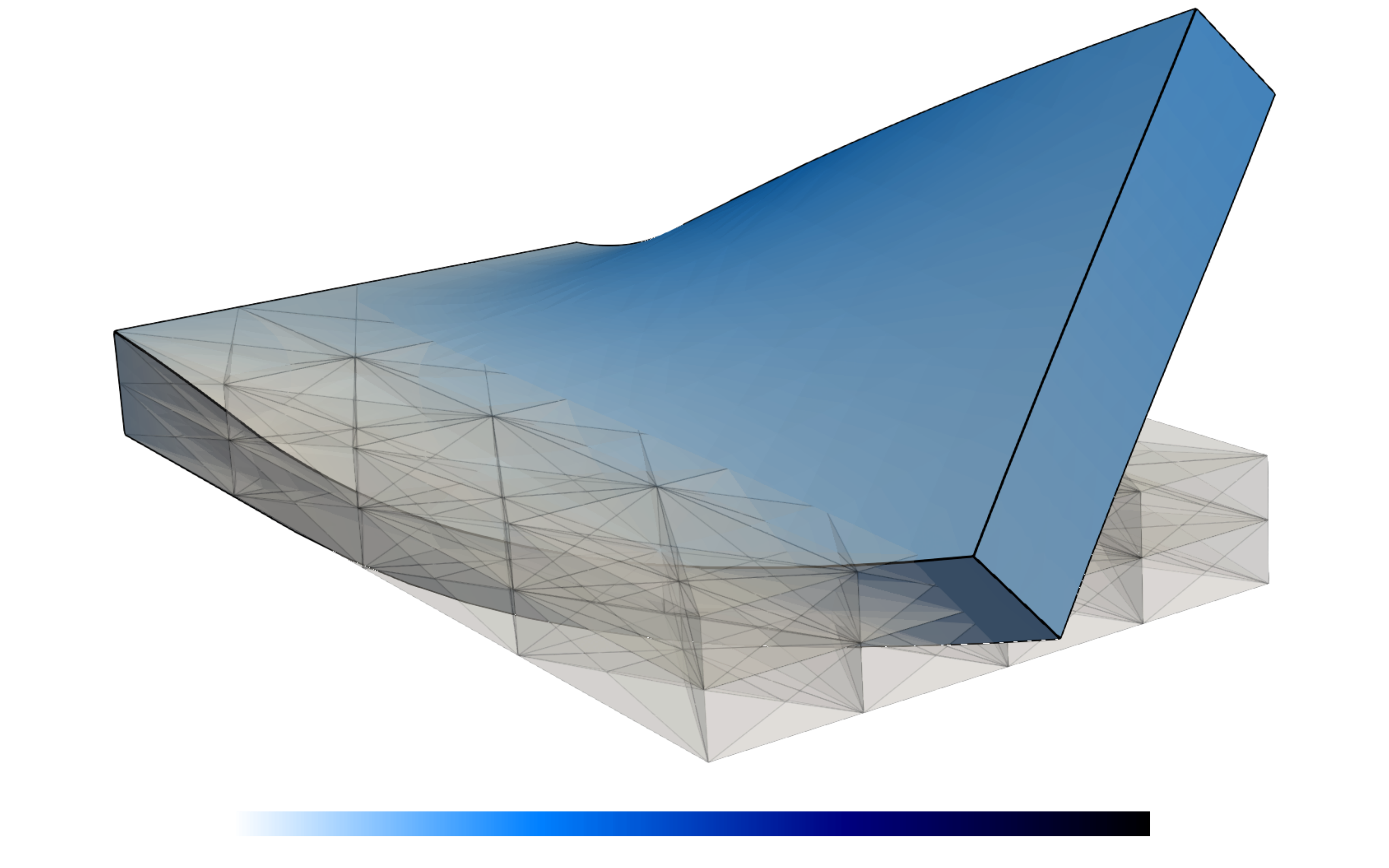}
            };
        \begin{scope}[x={(image.south east)}, y={(image.north west)}]
        \node[black, font=\small] at (0.09, 0.1) {0.1};
        \node[black, font=\small] at (0.89, 0.1) {-1};
        \end{scope}
        \end{tikzpicture}
        \caption{}
    \end{subfigure}\hfill
    \begin{subfigure}{0.24\textwidth}
        \centering
        \begin{tikzpicture}
        \node[anchor=south west, inner sep=0] (image) at (0,0)
            {
            \includegraphics[width=\linewidth]{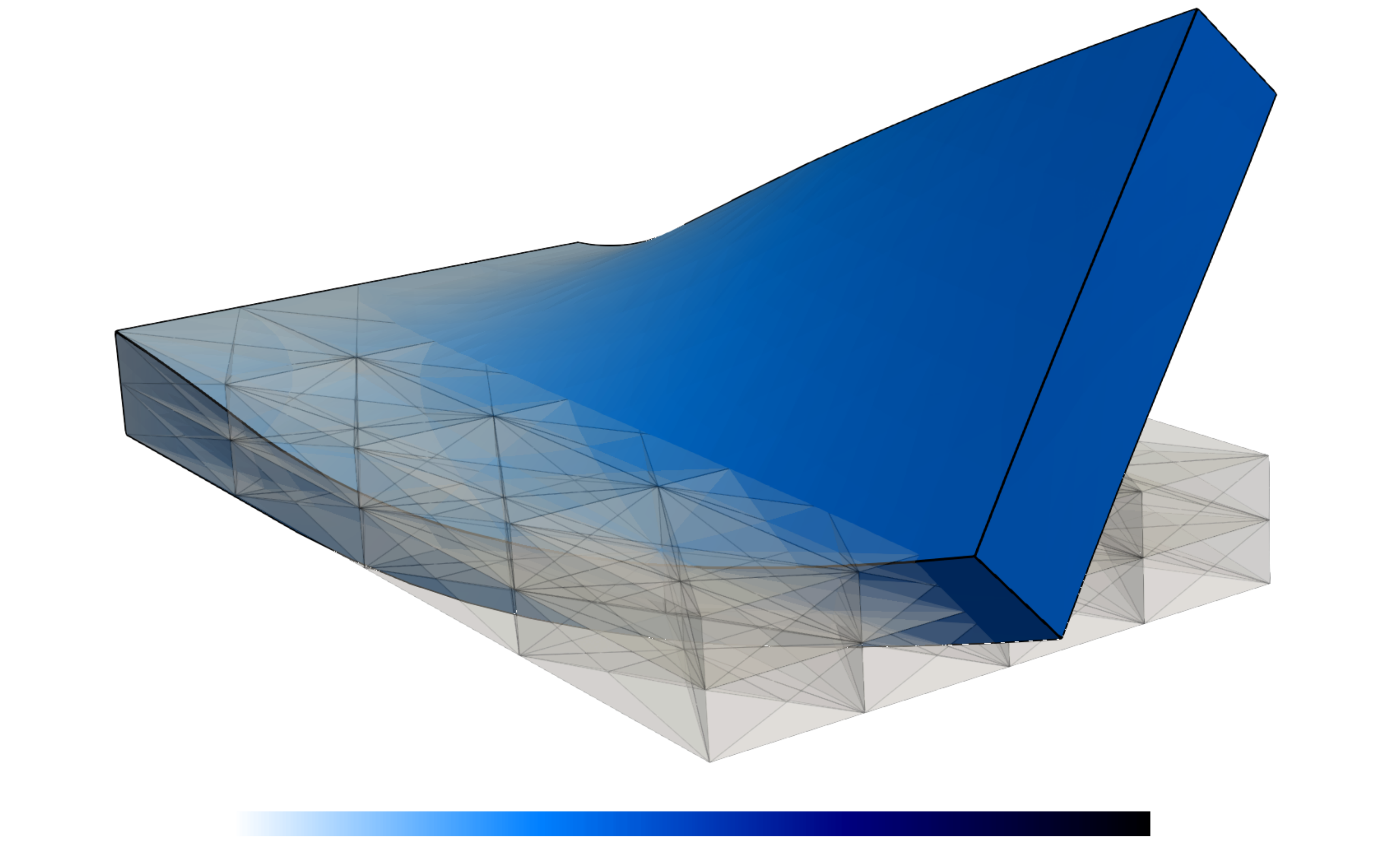}
            };
        \begin{scope}[x={(image.south east)}, y={(image.north west)}]
        \node[black, font=\small] at (0.09, 0.1) {-0.1};
        \node[black, font=\small] at (0.89, 0.1) {1};
        \end{scope}
        \end{tikzpicture}
        \caption{}
      \end{subfigure}\hfill
    \caption{\AS{Comparison of deformed configurations for the biaxial deformation with $\muc/\mu=10^3$, showing the high-fidelity reference solution (green), Γ-SPIN method (blue), and standard approach (violet) in (a). Panels (b)–(d) use a different perspective to show the rotation distribution of the method about (b) $\bm{e}_1$, (c) $\bm{e}_2$, and (d) $\bm{e}_3$.
    }}
    \label{fig:combined_torsion}
\end{figure}

\AS{In the fourth benchmark, we again consider a wide plate with domain $\overline{\Vol} = [0,2]\times[0,1]\times[0,0.2]$, as shown in \cref{fig:thin_plate}.}
\begin{figure}
    \centering
    \input{figs/thin_plate}
    \caption{\AS{Thin and wide plate domain with corresponding homogeneous Dirichlet and Neumann boundary surfaces. The domain is loaded by a distributed volume couple $\vb{m}$.}}
    \label{fig:thin_plate}
\end{figure}
\AS{The left end of the plate is prescribed as a homogeneous Dirichlet surface, while all other boundaries are treated as natural Neumann surfaces. The Cosserat couple modulus is set to $\muc/\mu=10^3$, with all other material parameters retained from the previous examples. In this benchmark, the loading is chosen to induce rotations about multiple axes, allowing the proposed formulation to be assessed under general three-dimensional rotations. The load is given by the constant distributed couple moment $\vb{m}=100\bm{e}_1+50\bm{e}_3$. The resulting deformation patterns are shown in \cref{fig:combined_torsion}.}
\AS{
We compare three cases consisting of a high-fidelity reference formulation with $\Coss_h \in \GEO^3(\Vol)$ and $\defmap_h \in \mathbb{R}^3 \otimes \CG^3(\Vol)$, the newly proposed formulation, and the standard lower-order formulation with $\Coss_h \in \GEO^2(\Vol)$ and $\defmap_h \in \mathbb{R}^3 \otimes \CG^2(\Vol)$. The high-fidelity solution exhibits the largest deformation and serves as the reference. Relative differences are evaluated using the displacement magnitude at the corner of the right end where the maximum displacement occurs. The proposed structure-preserving formulation yields a relative difference of $5.4\%$ compared with $12.8\%$ for the standard formulation. In contrast to the previous examples, this deformation involves rotations of comparable magnitude about all coordinate axes.
}

\subsection{A \AS{precurved} wide spring}
In all previous benchmarks, the interpolation alone performed relatively well, which may lead to the false impression that it is sufficient in general. In the final benchmark, however, we present a clear example in which it fails decisively. To understand this behaviour, it is important to recall that the earlier tests involved almost purely simple curvature about a single principal axis. In other words, no complex membrane-bending-torsion coupling modes were activated. Now, in order to induce such a mode, we a define the more intricate domain, given by 
$\overline{\Vol} = \{\bm{\xi} \in [0,1]^3 \; | \; \vb{x} = 200 \xi \vb{e}_1 + 10[2\eta-1][1.5+\cos(2\pi[2\xi-1])]\vb{e}_2 + 4[\zeta - 2.5\sin(2\pi(2\xi-1)]\vb{e}_3\}$, which is depicted in \cref{fig:spring}. 
\begin{figure}
    \centering
    \input{figs/spring}
    \caption{Domain of a curved wide spring with Dirichlet and Neumann boundaries. Importantly, the Dirichlet boundary of the Cosserat rotation tensor is comprised of both ends of the domain in this benchmark.}
    \label{fig:spring}
\end{figure}
The domain can be understood as a form of a \AS{precurved} spring with a varying width. 
The left end of the domain, given by $\{\vb{x} \in \partial \Vol \; | \; x = 0\}$, is defined as a Dirichlet boundary for both the displacement field and the rotation tensor. For the rotation tensor, also the right end of the domain, given by $\{\vb{x} \in \partial \Vol \; | \; x = 200\}$, is Dirichlet. On that same surface we impose the traction for the deformation field $\vb{t} = -2.5 \vb{e}_1$. In total, it yields a force of $F = 250$. The elastic Lam\'e parameters of the domain are set to $\lambda=692.398$ and $\mu=461.538$. The weights acting on the micro-dislocation tensor are given by $\alpha_1\Lc^2=0.288$, $\alpha_2\Lc^2=1.711$, and $\alpha_3\Lc^2=0.42$. Finally, the Cosserat couple modulus is defined as $\muc = 500 \mu$, representing a moderate coupling ratio. For the comparison we use the $\GEO^{2}(\Vol)$-type formulation.

\begin{figure}
    \centering
  \begin{subfigure}{\textwidth}
    \centering
    \includegraphics[width=\linewidth]{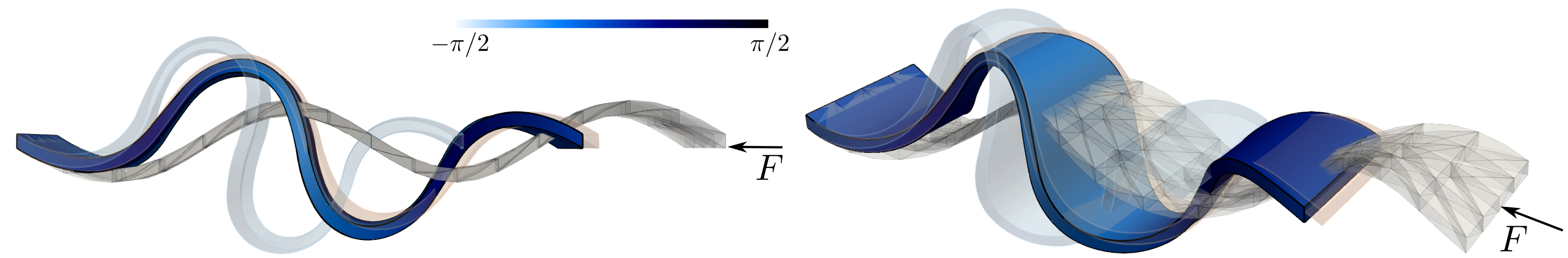}
    \caption{}
  \end{subfigure}\hfill
  \begin{subfigure}{\textwidth}
    \centering
    \includegraphics[width=\linewidth]{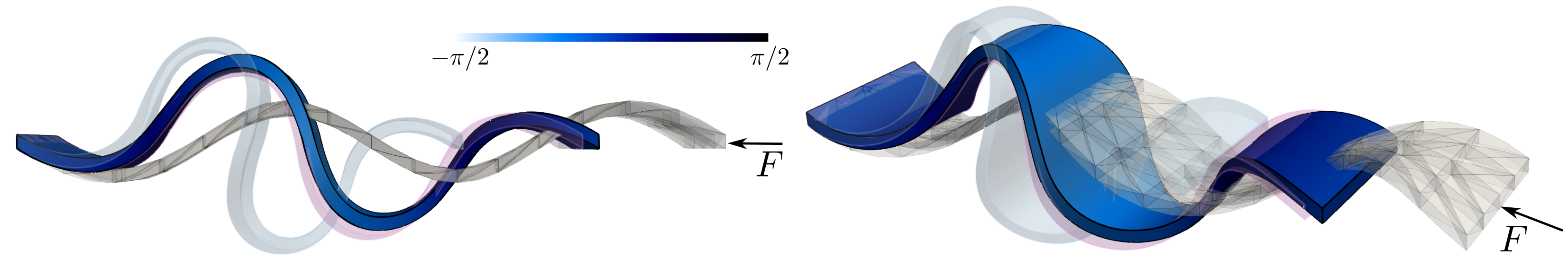}
    \caption{}
  \end{subfigure}\hfill
  \begin{subfigure}{\textwidth}
    \centering
    \includegraphics[width=\linewidth]{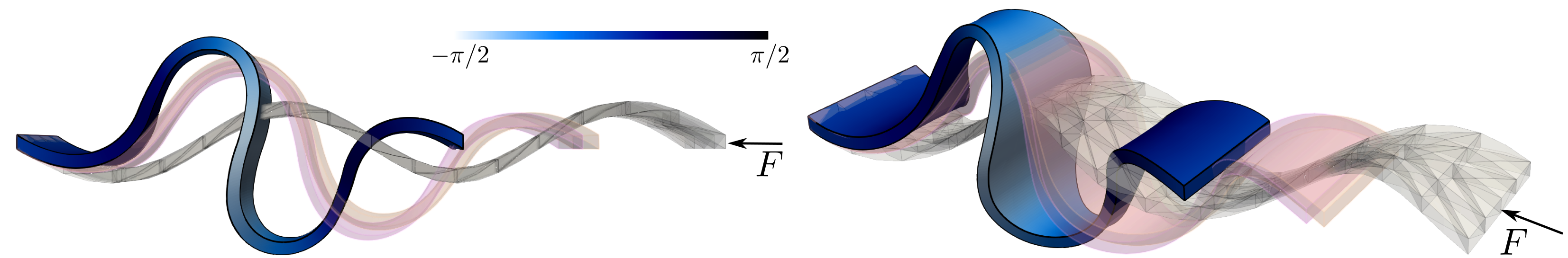}
    \caption{}
  \end{subfigure}
    \caption{
    \AS{
    Deformation of the curved spring obtained with the standard approach (a), with interpolation only (b), and with Γ-SPIN using interpolation followed by subsequent projection (c).
    }
    }
    \label{fig:spring_deformation}
\end{figure}

\begin{figure}
    \centering
    \begin{tikzpicture}
    \node[anchor=south west, inner sep=0] (image) at (0,0)
            {
            \includegraphics[width=\linewidth]{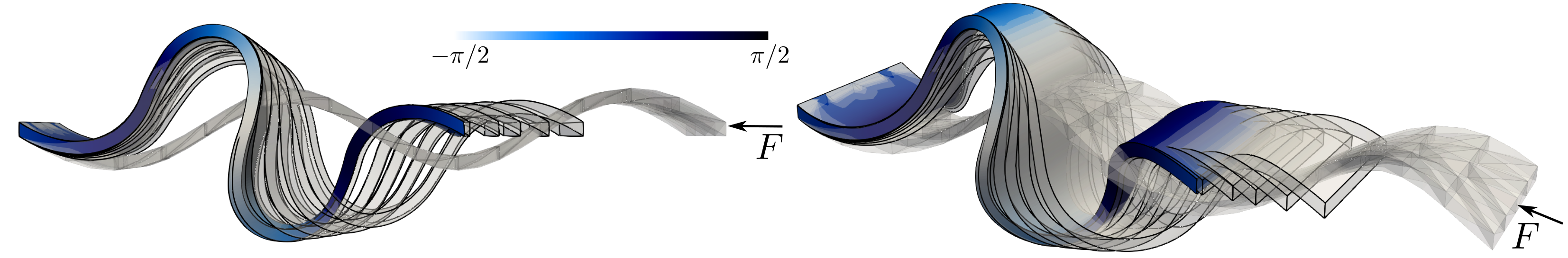}
            };
    \begin{scope}[x={(image.south east)}, y={(image.north west)}]
        \node[black, inner sep=1pt] (n1) at (0.88, 0.75) {$p=2$ with Γ-SPIN};
        \draw[] (n1.west) -- ++(-0.02,-0.14);

        \node[black, inner sep=1pt] (n2) at (0.4, 0.74) {$p=2$};
        \draw[] (n2.west) -- ++(-0.013, -0.2);

        \node[black, inner sep=1pt] (n3) at (0.4, 0.4) {$p=4$};
        \draw[] (n3.west) -- ++(-0.025, 0.1);

        \node[black, inner sep=1pt] (n4) at (0.39, 0.28) {$p=6$};
        \draw[] (n4.west) -- ++(-0.034, 0.22);

        \node[black, inner sep=1pt] (n5) at (0.38, 0.16) {$p=8$};
        \draw[] (n5.west) -- ++(-0.04, 0.34);
        
        \node[black, inner sep=1pt] (n6) at (0.32, 0.08) {$p=12$};
        \draw[] (n6.north) -- ++(-0.0215, 0.36);
        
    \end{scope}
    \end{tikzpicture}
    \caption{\AS{Deformation of the precurved spring obtained with Γ-SPIN in blue compared with solutions employing $\Coss_h \in \GEO^p(\Vol)$ with $p\in{2,4,6,8,12}$, shown in transparent grey.}}
\label{fig:spring_deformation_high_fidelity}
\end{figure}

The resulting deformations of all three formulations are depicted in \cref{fig:spring_deformation}, for which the corresponding mesh \AS{consists of} $384$ elements.  
As illustrated in the first and second figures, both the standard approach and the method that applies only interpolation produce significantly smaller deformations than the one obtained with Γ-SPIN, shown in the third figure. This pronounced discrepancy clearly reveals that locking can produce substantially different results, \AS{here amounting to about $45\%$ in the relative deformation}, and must not be neglected.
In contrast to the previous examples, this is the first instance in which the interpolation-only method performs worse than the standard approach.
The root of this behaviour lies in the nature of the interpolation procedure itself. When applied in isolation, the rotation field is embedded and interpolated in the full matrix space $\mathbb{R}^{3\times 3}$ rather than being constrained to the rotation group. As a result, its physical interpretation changes. Namely, instead of representing a pure rotation, it effectively becomes a general second-order tensor resembling a deformation measure.
This modification has important consequences. The interpolated field is no longer restricted to interacting only with the rotational modes of the discrete deformation tensor $\defgrad_h$. Instead, it can also couple with shear and dilatational modes. Put simply, by interpolating the rotation field in the ambient matrix space without enforcing the rotational constraint, the method introduces additional deformation mechanisms that are not physically admissible. Thus, applying only the interpolation procedure can yield \textbf{unphysical solutions}. \AS{Finally, the validity of the softer response obtained with the Γ-SPIN method is assessed by comparison with high-fidelity solutions in which the polynomial degree of the rotation field is increased, $\Coss_h \in \GEO^p(\Vol)$ with $p\in{2,4,6,8,12}$, while retaining $\defmap_h \in \mathbb{R}^3 \otimes \CG^2(\Vol)$ for the displacement field. Enriching the rotation space progressively reduces the weak discretisation mismatch between $\Coss_h$ and $\defgrad_h$, cf. \cite{Ainsworth2022}. The resulting deformations are shown in \cref{fig:spring_deformation_high_fidelity}.
}
\AS{
The response softens with increasing polynomial order and decreasing weak mismatch. High-fidelity results from the standard approach are compared with the Γ-SPIN solution obtained by interpolation followed by projection with $\Coss_h \in \GEO^2(\Vol)$. For $p=2$, the relative difference in displacement magnitude at the top corner of the loaded end is $45.49\%$ compared with the structure-preserving formulation. This discrepancy decreases with increasing polynomial order for $\Coss_h$, reaching $2.11\%$ at the highest order considered, $p=12$. This behaviour clearly indicates convergence.
}

\section{Conclusions and outlook}
In this work, we introduced the Geometric Structure-Preserving Interpolation (Γ-SPIN) method for structure-preservation in the finite-strain Cosserat micropolar model. The method advocates the use of geodesic finite elements for the interpolation $\Coss_h$. This choice ensures objectivity under superimposed rigid body motions and preserves the correct measurement of curvature rates. 
In the limit as the Cosserat couple modulus tends to infinity, $\muc \to +\infty$, the model approaches its descendant couple-stress theory. This limit induces locking in the numerical solution, as $\Coss_h \to \polar \defgrad_h$ cannot be sufficiently satisfied in the discrete sense.
To mitigate this effect, the method introduces a structure-preserving relaxation strategy in the coupling between $\Coss_h$ and $\defgrad_h$.
The relaxation mechanism proceeds in two steps. First, $\Coss_h$ is interpolated into the intrinsic regularity of $\defgrad_h$ using the N\'ed\'elec interpolant. Second, the resulting field is projected back onto the Lie-group $\SO(3)$ via the polar extraction. This procedure yields a projection-based interpolation with reduced regularity. The discrete objectivity of the this approach was briefly verified in the first benchmark. To further assess its correctness and effectiveness, we additionally benchmarked two possible discretisations within the method on beam bending and torsion problems, as well as on a multi-mode response induced by a curved geometry. The numerical results demonstrate that the method maintains optimal convergence rates in all considered scenarios and successfully circumvents locking. Moreover, since neither improved convergence rates for the enriched formulation $\GEO^{2+}(\Vol)$ nor oscillations for the $\GEO^{2}(\Vol)$-formulation are observed, we conjecture that the target space $\polar[\R^3 \otimes \Ned_{II}^1(\Vol)]$ is already fully spanned by $\polar[\Pi_c^1 \GEO^2(\Vol)]$. 
In cases involving large applied forces, the deformation predicted by the standard approach differs substantially from that obtained with our proposed method, highlighting the pronounced impact of locking effects at finite strains. Finally, we showed that a naive attempt to match regularities without respecting the true kinematics of the model leads to unphysical solutions.

While this work provides evidence of stability, consistency and optimality through numerical benchmarks, supplemented by analytical indicators based on dimension counting, a rigorous stability analysis is still required and will be addressed in future work. Further, to avoid the implementation complexity associated with geodesic elements, we employed Euler rotation matrices instead. Although this substitution introduces approximation errors, it was applied consistently across all tested methods to ensure that it does not account for the observed differences. Nevertheless, in future work the approach will be revalidated using true geodesic elements, which require a dedicated low-level algorithmic implementation.

\section*{Acknowledgements}
Adam Sky acknowledges support by the European Commission within the \textit{Marie Skłodowska-Curie Actions} postdoctoral fellowship programme, project AGE2M, grant 101152236. Lucca Schek acknowledges support by the DFG under Project-No. 525235558.
Patrizio Neff acknowledges support in the framework of the DFG-Priority Programme 2256 “Variational Methods for Predicting Complex Phenomena in Engineering Structures and Materials”, Neff 902/10-1, Project-No. 440935806. 

\bibliographystyle{spmpsci}   

\footnotesize
\bibliography{ref}   

\normalsize
\appendix

\section{Derivations and formulae}
In the following we specify certain helpful formulae and elaborate on relevant relations and procedures concerning rotation matrices. 

\subsection{The matrix exponential and logarithm} \label{ap:mexp}

In the general case, matrix products do not commute $\bm{A}\bm{B} \neq \bm{B}\bm{A}$. As a consequence, the exponential function of matrices does not coincide with the scalar one, and is instead given by the Lie formula 
\begin{align}
    \exp(\bm{A} + \bm{B}) = \lim_{n \to \infty} [ \exp(n^{-1}\bm{A}) \exp(n^{-1}\bm{B}) ]^n \, .
\end{align}
Let the commutator be defined in terms of the adjoint operator
\begin{align}
        \ad_{\bm{A}} \bm{B} = \bm{A}\bm{B} - \bm{B}\bm{A} \, , && \ad_{\bm{A}}^2 \bm{B} = \bm{A}(\bm{A}\bm{B} - \bm{B}\bm{A}) - (\bm{A}\bm{B} - \bm{B}\bm{A}) \bm{A} \, ,  
\end{align}
the product of two matrix exponentials can be expressed over the Baker--Campbell--Hausdorff formula
\begin{align}
    (\exp\bm{A})(\exp\bm{B}) = \exp\bigg[\bm{A} + \bm{B} + \dfrac{1}{2}\ad_{\bm{A}}(\bm{B}) + \dfrac{1}{12}\ad^2_{\bm{A}}(\bm{B}) + \dfrac{1}{12}\ad^2_{\bm{B}}(\bm{A}) + \dots \bigg] \, , 
\end{align}
which clearly reveals how the non-commutation shifts the result from the classical exponential.  
By defining $\bm{A} = \log \bm{X}$ and $\bm{B} = \log \bm{Y}$, the same can be shown for the matrix logarithm
    \begin{align}
        \log(\bm{X}\bm{Y}) = \log\bm{X} + \log\bm{Y} + \dfrac{1}{2}\ad_{\log\bm{X}}(\log\bm{Y}) + \dfrac{1}{12}\ad^2_{\log\bm{X}}(\log\bm{Y}) + \dfrac{1}{12}\ad^2_{\log\bm{Y}}(\log\bm{X}) + \dots  \, .  
    \end{align}

\subsection{Total multiplicative updating} \label{ap:mul}

We start by recalling multiplicative updating before differing between incremental and total updates. Observably, rotations define a multiplicative group such that any rotation can be written as the product of previous rotations 
\begin{align}
    \bm{R} = \prod_k \bm{R}_k \, , && \bm{R},\,\bm{R}_k \in \SO(3) \, .
\end{align}
Consequently, given a rotation tensor $\bm{R}_k$ and an incremental rotation $\Delta \bm{R}_{k+1}$, the new total rotation reads  
\begin{align}
    \bm{R}_{k+1} = \Delta\bm{R}_{k+1} \bm{R}_k \, .
\end{align}
If the rotations are parametrised by axial vectors, this implies that the new rotation must be computed as
\begin{align}
    \bm{R}_{k+1}(\bm{\theta}_{k+1}) = \exp(\Anti \Delta \bm{\theta}_{k+1}) \exp(\Anti  \bm{\theta}_{k}) \, , && \bm{\theta}_{k+1} = \log \bm{R}_{k+1} = \log[\exp(\Anti \Delta \bm{\theta}_{k+1}) \exp(\Anti  \bm{\theta}_{k})] \, .     
\end{align}
The matrix exponential and logarithm do not general obey the same rules as their scalar counterparts 
\begin{align}
    \exp(\Anti \Delta \bm{\theta}_{k+1}) \exp(\Anti  \bm{\theta}_{k}) \neq \exp(\Anti \Delta \bm{\theta}_{k+1} + \Anti  \bm{\theta}_{k}) = \exp[\Anti  (\underbrace{\Delta \bm{\theta}_{k+1}+ \bm{\theta}_{k}}_{\widetilde{\bm{\theta}}_{k+1}})] \, ,
\end{align}
such that an additive update does \textbf{not} yield the true new axial vector $\bm{\theta}_{k+1} \neq \widetilde{\bm{\theta}}_{k+1}$. Nevertheless, by the Baker--Campbell--Hausdorff formula
\begin{align}
    \bm{\theta}_{k+1} = \underbrace{\bm{\theta}_{k} + \Delta\bm{\theta}_{k+1}}_{\widetilde{\bm{\theta}}_{k+1}} + \dfrac{1}{2}\bm{\theta}_{k} \times \Delta \bm{\theta}_{k+1} + \dots \, , 
\end{align}
the additive update represents a first order approximation of the true axial vector.

A related problem is that of path-independence, which is the independence of the final solution step from intermediate computations. This is important, as the final deformation result of an elastic system should in general not depend on intermediate configurations. Path-independence is guaranteed by only interpolating the total multiplicative updates. This is readily observed with an example. Let $\bm{R}^1_{k} = \Delta\bm{R}^1_{k} \Delta\bm{R}^1_{k-1}\dots \bm{R}_0^1$ and $\bm{R}^2_{k} = \Delta\bm{R}^2_{k} \Delta\bm{R}^2_{k-1}\dots \bm{R}_0^2$ be two neighbouring nodal values on nodes $1$ and $2$ at iteration step $k$, the total linear geodesic interpolation on a curve is given by
\begin{align}
    \Pi_s^1\bm{R}_k = \bm{R}_k^1\exp[ \lambda_2 \log([\bm{R}_k^1]^T \bm{R}_k^2) ] \, , && \lambda_2 \in \CG^1(s) \, .
\end{align}
Evidently, this interpolation is blind to the composition of $\bm{R}_k^i = \Delta \bm{R}_k^i \Delta\bm{R}_{k-1}^i\dots\bm{R}_{0}^i$ from the previous steps and thus independent of the path the computation takes. In comparison, in an incremental updating procedure the increments are essentially interpolated separately 
\begin{align}
    (\Pi_s^1 \Delta \bm{R}_k)(\Pi_s^1 \Delta \bm{R}_{k-1})\dots (\Pi_s^1 \bm{R}_0) \,  \neq \Pi_s^1\bm{R}_k \, ,
\end{align}
which is not equal to the total update. Namely, set $\bm{R}_1^1 = \one_{1}\one_{0} = \one$ and $\bm{R}_1^2 = \Delta\bm{R}_{1} \bm{Q}_{0}$ we generally find
\begin{align}
     \Pi_s^1\bm{R}_1 = \exp[ \lambda_2 \log( \Delta\bm{R}_{1} \bm{Q}_0) ] \neq \exp[ \lambda_2 \log (\Delta\bm{R}_{1}) ] \exp[ \lambda_2 \log \bm{Q}_0 ]  \, .  
\end{align}
Let $\lambda_2:\curv \to [0,1]$, then the two methods agree at the nodes $\lambda_2 = 0$ and $\lambda_2 = 1$. However, they do not agree for $\lambda_2:\curv \to (0,1)$ since 
\begin{align}
    \exp[ \lambda_2 \log (\Delta\bm{R}_{1}) ] \exp[ \lambda_2 \log \bm{Q}_0 ] &\neq \exp[\lambda_2 \log (\Delta\bm{R}_{1})  + \lambda_2 \log \bm{Q}_0] = \exp[\lambda_2 \log (\Delta\bm{R}_{1})  + \lambda_2 \log \bm{Q}_0] \, , 
    \notag \\
     \exp(\lambda_2 [\log (\Delta\bm{R}_{1})  +  \log \bm{Q}_0]) &\neq \exp[\lambda_2 \log (\Delta\bm{R}_{1}\bm{Q}_0) ] \, .
\end{align}
This is crucial, since it implies that any strains and curvatures evaluated between nodes will depend on the history of deformations. By composition, this implies that after initial steps also the nodal values will differ. Thus, these two effects combined lead to path-dependence of the final result.

\subsection{Logarithm of a rotation tensor} \label{ap:log}
The Euler--Rodrigues formula gives a closed-form formulate for the representation of every rotation tensor. By expanding $(\Anti \bm{\theta})^2$, it can be rewritten as
\begin{align}
    \bm{R}(\bm{\theta}) = (\cos \norm{\bm{\theta}}) \one + \dfrac{\sin\norm{\bm{\theta}}}{\norm{\bm{\theta}}} \Anti \bm{\theta} + \dfrac{1-\cos \norm{\bm{\theta}}}{\norm{\bm{\theta}}^2} \bm{\theta} \otimes \bm{\theta} \, .
\end{align}
Let $0< \norm{\bm{\theta}} < \pi$, then its skew-symmetric component does not vanish and we can extract
\begin{align}
    \axl \bm{R} = \dfrac{\sin\norm{\bm{\theta}}}{\norm{\bm{\theta}}} \bm{\theta} \quad \Rightarrow \quad \bm{\theta} = \dfrac{ \norm{\bm{\theta}}  }{\norm{\axl \bm{R}}} \axl \bm{R} \, ,
\end{align}
since $\sin\norm{\bm{\theta}} > 0$ is positive in that range, which implies that the two vector are parallel $\axl \bm{R} \parallel \bm{\theta}$. Naturally, the latter implies $\sin \norm{\bm{\theta}} = \norm{\axl \bm{R}}$ and thus with the inverse sine function $\norm{\bm{\theta}} = \arcsin\norm{\axl \bm{R}}$. However, since $\arcsin:[-1,1] \to [-\pi/2, \pi/2] \ni \norm{\bm{\theta}}$ it can violate the original assumption on the range and the resulting parallel orientation of the vectors. Consequently, the intensity of the rotation is be extracted in a different manner via
\begin{align}
    \tr \bm{R} = 1 + 2 \cos  \norm{\bm{\theta}} \quad \Rightarrow \quad \norm{\bm{\theta}} = \arccos \left ( \dfrac{\tr \bm{R} - 1}{2} \right ) \, ,
\end{align}
where $\arccos:[-1,1] \to [0,\pi] \ni \norm{\bm{\theta}}$ satisfies the chosen range. As such, the axial vector reads 
\begin{align}
    \bm{\theta} = \dfrac{ \arccos \left ( \dfrac{\tr \bm{R} - 1}{2} \right )  }{\norm{\axl \bm{R}}} \axl \bm{R} \, . 
\end{align}
While the range of $\norm{\bm{\theta}} \in [0,\pi]$ may seem limited, it is in fact sufficient, as any larger rotation $\norm{\bm{\theta}} > \pi$ can be reinterpreted as a smaller one around the flipped axis. 

The case $\bm{\norm{\theta}} = 0$ implies no rotation $\bm{\theta} = 0$. Although it leads to a singularity in the formula, it can be recovered by redefining the angle as an independent parameter $\alpha = \bm{\norm{\theta}}$. Observably, the zero rotation case is characterised by $\tr \bm{R} = 3$.

The last case is a half-circle rotation $\norm{\bm{\theta}} = \pi$, which is characterised by $\tr \bm{R} = -1$.
In this scenario, the representation of the rotation tensor reduces to 
\begin{align}
    \bm{R}(\bm{\theta}) = \dfrac{2}{\pi^2} \bm{\theta} \otimes \bm{\theta} -\one   \, ,
\end{align}
which can be rearranged into
\begin{align}
    \dfrac{1}{\pi^2} \bm{\theta} \otimes \bm{\theta} = \dfrac{1}{2}(\bm{R} + \one) \, . 
\end{align}
The latter is projection matrix onto the axis $\bm{\theta}$, since $\pi^{-1}\bm{\theta}$ is a unit vector. Thus, the vector can be recovered by testing against the Cartesian basis and taking any non-vanishing projection
\begin{align}
    \dfrac{1}{\pi} \bm{\theta} = \pm \dfrac{2^{-1}}{\norm{2^{-1}(\bm{R} + \one) \vb{e}_\xi}} (\bm{R} + \one) \vb{e}_\xi = \pm\dfrac{1}{\norm{\bm{R} \vb{e}_\xi + \vb{e}_\xi}} \bm{R} \vb{e}_\xi + \vb{e}_\xi \quad \Rightarrow \quad \bm{\theta} =\pm\dfrac{\pi}{\norm{\bm{R} \vb{e}_\xi + \vb{e}_\xi}} \bm{R} \vb{e}_\xi + \vb{e}_\xi \, .   
\end{align}
The sign cannot be recovered from the projection. But, assuming temporal continuity in an iterative solution procedure, one can fix the sign via $\con{\bm{\theta}_{k+1}}{\bm{\theta}_k} > 0$ for moderate change $\norm{\Delta \bm{\theta}}$. This implies the closed-form solution
\begin{align}
    \bm{\theta}_{k+1} =\dfrac{\pi}{\norm{\bm{R} \bm{\theta}_{k} + \bm{\theta}_{k}}} \bm{R} \bm{\theta}_{k} + \bm{\theta}_{k} \, .
\end{align}
\begin{observation}
    Let $\norm{\bm{\theta}_{k+1}} \in [0,\pi]$, then  
    \begin{align}
        \sym\bm{R}_{k+1} - \dfrac{\tr \bm{R}_{k+1} - 1}{2} \one = \dfrac{1-\cos\norm{\bm{\theta}_{k+1}}}{\norm{\bm{\theta}_{k+1}}^2} \bm{\theta}_{k+1} \otimes \bm{\theta}_{k+1} \, ,
    \end{align}
    is a weighted projection matrix. Under temporal continuity we assume $\con{\bm{\theta}_{k+1}}{\bm{\theta}_{k}} > 0$, such that the current axial vector can be recovered via
    \begin{align}
        \boxed{
        \bm{\theta}_{k+1} = \dfrac{\arccos \left ( \dfrac{\tr \bm{R}_{k+1} - 1}{2} \right )}{ \left \| \left ( \sym\bm{R}_{k+1} - \dfrac{\tr \bm{R}_{k+1} - 1}{2} \one   \right ) \bm{\theta}_k \right \| } \left ( \sym\bm{R}_{k+1} - \dfrac{\tr \bm{R}_{k+1} - 1}{2} \one   \right ) \bm{\theta}_k \, .
        }
    \end{align}
    This appears to be a \textbf{novel} formula, which is exact and does not require case-distinction.  
\end{observation}

\subsection{The polar decomposition is a minimiser} \label{ap:pol}
Let $\defgrad \in \GL^+(3)$, then the minimisation problem
\begin{align}
    \argmin_{\bm{Q} \in \SO(3)}\dfrac{1}{2} \norm{\bm{F} - \bm{Q}}^2 = \argmin_{\bm{Q} \in \SO(3)}\dfrac{1}{2} \norm{\bm{Q}^T\bm{F} - \one}^2 \, , 
\end{align}
defines the polar extraction \cite{Neff2014Grioli}. Namely, from
\begin{align}
    \dfrac{1}{2} \norm{\defgrad - \bm{Q}}^2 = \dfrac{1}{2} \norm{\bm{R}\bm{U} - \bm{Q}}^2 = \dfrac{1}{2} \norm{\bm{R}(\bm{U} - \bm{R}^T\bm{Q})}^2 = \dfrac{1}{2} \norm{\bm{U} - \bm{R}^T\bm{Q}}^2 = \dfrac{1}{2}\norm{\bm{U}}^2 + \dfrac{1}{2}\norm{\bm{R}^T \bm{Q}} - \con{\bm{U}}{\bm{R}^T \bm{Q}} \, ,
\end{align}
we get that the minimiser maximises $\con{\bm{U}}{\bm{R}^T \bm{Q}}$ with $\bm{U} \in \Sym^{++}(3)$. We employ the spectral decomposition and redefine
\begin{align}
    \bm{U} = \bm{V} \bm{\Lambda} \bm{V}^T \, , && \bm{W} = \bm{V}^T \bm{R}^T \bm{Q} \bm{V} \, , && \bm{V} \in \SO(3) \, , && \bm{\Lambda} \in \Diag(3) \, .
\end{align}
Consequently, the product now reads
\begin{align}
    \con{\bm{U}}{\bm{R}^T \bm{Q}} = \con{\bm{V} \bm{\Lambda} \bm{V}^T}{\bm{V} \bm{W} \bm{V}^T} = \con{\bm{\Lambda}}{\bm{W}} \leq \sum_i \sigma_i(\bm{\Lambda}) \sigma_i(\bm{W}) = \sum_i \Lambda_i = \tr \bm{\Lambda} = \con{\bm{\Lambda}}{\one}  \, ,
\end{align}
with Neumann's trace inequality, since the singular values of $\bm{W}\in\SO(3)$ are all $\sigma_i(\bm{W}) =  \sqrt{\Lambda_i(\bm{W}^T\bm{W})} =  \sqrt{\Lambda_i(\one)} = 1$ and those of $\bm{\Lambda}$ are precisely its eigenvalues $\sigma_i(\bm{\Lambda}) = \sqrt{\Lambda_i(\bm{\Lambda}^T\bm{\Lambda})} = \Lambda_i(\bm{\Lambda}) > 0$. Equality, and thus the maximum value, are given for $\bm{W} = \one$, which implies $\bm{Q} = \bm{R} = \polar \defgrad$. 

The more general problem is given by
\begin{align}
    \argmin_{\bm{Q} \in \SO(3)} (\mu \norm{\sym(\bm{Q}^T\bm{F} - \one)}^2 + \muc \norm{\skw(\bm{Q}^T\bm{F} - \one)}^2)  \, , 
\end{align}
which still yields the polar projection $\bm{Q} = \bm{R} = \polar \bm{F}$ for $\muc \geq \mu > 0$ \cite{Fischle2017}.

\section{Observations on the balance of angular momentum} \label{ap:obs}
We algebraically lift the balance of angular momentum to a matrix-valued equation via
\begin{align}
    -\dfrac{1}{2}\Anti[\Di \Double + 2\axl(\Piola \defgrad^T)] = -\dfrac{1}{2}\Anti(\Di \Double) - \skw(\Piola \defgrad^T) &= \dfrac{1}{2}\Anti\vb{m} \, , 
\end{align}
where the half-factor stems from $\con{\delta\bm{\omega}}{\delta\bm{\omega}} = (1/2)\con{\Anti\delta\bm{\omega}}{\Anti\delta\bm{\omega}}$. The latter demonstrates that a complete formulation of the balance of angular momentum should make the implied skew-symmetry explicit.
To finally make the connection to dislocation form, obverse that there holds \cite{Sky2024s}
\begin{align}
    &\skw \Curl (\sym  \Double) = \dfrac{1}{2}\Anti[\Di (\sym \Double) - \nabla (\tr \Double)] 
    \notag \\
    &\Rightarrow \quad \Anti [\Di (\sym \Double)] = 2 \skw \Curl(\sym \Double) + \Anti [\nabla (\tr \Double)] \, ,
\end{align}
for symmetric tensors, and \cite{Lewintan2021}
\begin{align}
    \Di (\skw \Double) = -\curl (\axl \Double)  \quad \Rightarrow \quad\Anti [\Di (\skw \Double)] = -\Anti[\curl (\axl \Double)] 
\end{align}
by Room's formula for skew-symmetric tensors. Now, by Nye's formula 
\begin{align}
    \Curl (\skw \Double) = \tr (\Grad \axl \Double) \one - (\Grad \axl \Double)^T \, ,
\end{align}
and another of Room's formulae
\begin{align}
    \Anti [\curl (\axl \Double)] = 2 \skw \Grad (\axl \Double) \, ,
\end{align}
there holds
\begin{align}
    \skw \Curl (\skw \Double) = \skw (\Grad \axl \Double) = \dfrac{1}{2} \Anti [\curl (\axl \Double)] \, ,
\end{align}
such that we can recombine the decomposition to find
\begin{align}
    \Anti [\Di (\underbrace{\sym \Double + \skw \Double}_{=\Double})] =  2 \skw \Curl(\underbrace{\sym \Double - \skw \Double}_{=\Double^T}) + \underbrace{\Anti [\nabla (\tr \Double)]}_{=3\Curl(\sph \Double)}
    \, ,
\end{align}
representing a general relationship between the divergence and the Curl of a second-order tensor
\begin{align}
    \Di \Double = 2\axl \Curl(\Double^T) + \nabla (\tr \Double) \quad \iff \quad \skw \Curl (\Double) = \dfrac{1}{2}[\Anti (\Di \Double^T) - \nabla (\tr \Double)] \, .
\end{align}
Consequently, it is possible to express the balance of angular momentum as 
\begin{align}
    - \skw\Curl(\Double^T) - \dfrac{1}{2}\Anti[\nabla (\tr \Double)]  - \skw (\Piola \defgrad^T) = \dfrac{1}{2}\Anti \vb{m} \, . 
\end{align}
Clearly, if the couple-stress is deviatoric $\Double:\Vol \to \sl(3)$, then one finds
\begin{align}
    - \skw\Curl(\Double^T) - \skw (\Piola \defgrad^T) = \dfrac{1}{2}\Anti \vb{m} \, . 
\end{align}
If the couple-stress is also symmetric $\Double:\Vol \to \sl(3) \cap \Sym(3)$, then the transposition can be dropped and there holds the complete equality $\Anti (\Di \Double)= 2 \skw \Curl \Double$.

\end{document}

%% file: figs/dom.tex
\definecolor{qqqqff}{rgb}{0.,0.,1.}
\definecolor{qqwwzz}{rgb}{0.,0.4,0.6}
\definecolor{xfqqff}{rgb}{0.4980392156862745,0.,1.}
\begin{tikzpicture}[scale=0.7,line cap=round,line join=round,>=triangle 45,x=1.0cm,y=1.0cm]
\draw [-to,line width=0.7pt] (-4.,-1.) -- (-2.,-1.);
\draw [-to,line width=0.7pt] (-4.,-1.) -- (-4.,1.);
\draw [-to,line width=0.7pt,color=xfqqff] (-4.,-1.) -- (3.,4.);
\draw [-to,line width=0.7pt] (3.,4.) -- (3.5,5.);
\draw [-to,line width=0.7pt] (3.,4.) -- (4.,3.5);
\draw[line width=0.7pt,dashed,color=qqwwzz, smooth,samples=100,domain=0.0:0.5546146682700169] plot[parametric] function{-2.6934406765449803*t**(3.0)+7.260917604035623*t-0.5675156898573737,-1.8886702459478892*t**(3.0)-5.4145245008536635*t+7.325178164784006};
\draw[line width=0.7pt,dashed,color=qqwwzz, smooth,samples=100,domain=0.5546146682700169:1.0] plot[parametric] function{3.353998438890217*t**(3.0)-10.061995316670652*t**(2.0)+12.841447798725381*t-1.5991969907569312,2.3518606188918723*t**(3.0)-7.055581856675618*t**(2.0)-1.5013953099615658*t+6.601751882082559};
\draw[line width=0.7pt,dashed,color=qqwwzz, smooth,samples=100,domain=0.0:0.2687701594778118] plot[parametric] function{0.29542127152800696*t**(3.0)+16.582104430633116*t-1.4625105322529643,-7.350428129646393*t**(3.0)+3.17705267716275*t+3.2888134975811463};
\draw[line width=0.7pt,dashed,color=qqwwzz, smooth,samples=100,domain=0.2687701594778118:1.0] plot[parametric] function{-0.10858476755409584*t**(3.0)+0.3257543026622875*t**(2.0)+16.49455139475599*t-1.454666651114476,2.7017165207924925*t**(3.0)-8.105149562377477*t**(2.0)+5.355475017634461*t+3.0936485242949425};
\draw [-to,line width=0.7pt,color=qqqqff] (11.505074918500053,8.432113153914122) -- (12.16,7.12);
\draw [-to,shift={(13.540008018208527,1.4762271419804813)},line width=0.7pt,color=qqqqff]  plot[domain=-2.995837369692063:0.4941669386155021,variable=\t]({1.*0.7643151038537889*cos(\t r)+0.*0.7643151038537889*sin(\t r)},{0.*0.7643151038537889*cos(\t r)+1.*0.7643151038537889*sin(\t r)});
\draw [-to, shift={(10.764884455388097,4.830054134930684)},line width=0.7pt,color=qqqqff]  plot[domain=-3.061393963625155:1.2036224929766797,variable=\t]({1.*0.7277036652303261*cos(\t r)+0.*0.7277036652303261*sin(\t r)},{0.*0.7277036652303261*cos(\t r)+1.*0.7277036652303261*sin(\t r)});
\draw [-to,line width=0.7pt,color=qqqqff] (14.,6.) -- (15.,5.);
\draw[line width=0.7pt,color=qqwwzz,fill=qqwwzz,fill opacity=0.05000000074505806] (0.9859719520200088,0.8345810181044726)(3.722870688910685E-7,8.000000438638764) -- (7.445741577102429E-7,8.000000877277412) -- (1.4891483951329018E-6,8.000001754554381) -- (2.9782971091154048E-6,8.000003509107016) -- (5.956595493628717E-6,8.000007018207059) -- (1.1913196088845086E-5,8.00001403638625) -- (2.3826412584008935E-5,8.000028072661056) -- (4.76529067930382E-5,8.00005614487636) -- (9.530614008411985E-5,8.000112287969744) -- (1.9061358614411083E-4,8.00022456880767) -- (3.8123239606128464E-4,8.000449109088697) -- (7.624856861714485E-4,8.000898104075644) -- (0.0015250549401914328,8.001795751782803) -- (0.0030504440850011617,8.003589678399635) -- (0.006102224454268727,8.007172058598984) -- (0.012209789771948818,8.014314944105399) -- (0.024440908808309115,8.028513353508972) -- (0.0489668611603923,8.056561829558083) -- (0.09827170838490609,8.111274240056218) -- (0.19787785613849992,8.215231530108431) -- (0.29878343348518505,8.31203332262817) -- (0.40095343064936484,8.401841070086958) -- (0.5043528378554428,8.484816224956317) -- (0.6089466453278227,8.561120239707769) -- (0.7146998432909076,8.630914566812837) -- (0.8215774219691016,8.694360658743044) -- (0.9295443715868077,8.751619967969912) -- (1.0385656823684295,8.802853946964966) -- (1.1486063445383705,8.848224048199725) -- (1.2596313483210348,8.887891724145714) -- (1.371605683940825,8.922018427274455) -- (1.4844943416221454,8.950765610057472) -- (1.5982623115893988,8.974294724966285) -- (1.712874584066989,8.992767224472418) -- (1.8282961492793197,9.006344561047394) -- (1.9444919974507942,9.015188187162735) -- (2.061427118805816,9.019459555289965) -- (2.1790665035687886,9.019320117900605) -- (2.2973751419641157,9.014931327466178) -- (2.4163180242162,9.006454636458207) -- (2.5358601405494463,8.994051497348213) -- (2.655966481188257,8.977883362607724) -- (2.7766020363570365,8.958111684708255) -- (2.897731796280188,8.934897916121333) -- (2.9584710624697785,8.922050701717428) -- (3.019320751182114,8.908403509318482) -- (3.080276486195244,8.893976520483431) -- (3.1413338912872195,8.878789916771222) -- (3.2024885902360904,8.862863879740788) -- (3.2637362068199076,8.846218590951075) -- (3.3250723648167204,8.82887423196102) -- (3.386492688004581,8.810850984329566) -- (3.4479928001615385,8.79216902961565) -- (3.509568325065644,8.772848549378217) -- (3.571214886494948,8.752909725176202) -- (3.6329281082275005,8.732372738568548) -- (3.6947036140413516,8.711257771114198) -- (3.7565370277145522,8.689585004372088) -- (3.818423973025153,8.667374619901159) -- (3.8803600737512047,8.644646799260352) -- (3.942340953670757,8.621421724008608) -- (4.00436223656186,8.597719575704868) -- (4.066419546202566,8.57356053590807) -- (4.128508506370923,8.548964786177157) -- (4.190624740844982,8.523952508071066) -- (4.252763873402795,8.498543883148741) -- (4.3149215278224125,8.472759092969119) -- (4.377093327881883,8.446618319091142) -- (4.439274897359258,8.420141743073751) -- (4.501461860032587,8.393349546475886) -- (4.563649839679922,8.366261910856485) -- (4.625834460079313,8.33889901777449) -- (4.688011345008809,8.311281048788846) -- (4.750176118246462,8.283428185458485) -- (4.812324403570322,8.25536060934235) -- (4.87445182475844,8.227098501999386) -- (4.9365540055888655,8.198662044988527) -- (4.998626569839649,8.17007141986872) -- (5.060665141288841,8.141346808198897) -- (5.122665343714493,8.112508391538004) -- (5.184622800894655,8.083576351444982) -- (5.2465331366073755,8.054570869478766) -- (5.308391974630707,8.025512127198304) -- (5.3701949387427,7.996420306162531) -- (5.4319376527214045,7.967315587930387) -- (5.49361574034487,7.938218154060816) -- (5.555224825391147,7.909148186112753) -- (5.616760531638288,7.880125865645144) -- (5.678218482864342,7.851171374216925) -- (5.739594302847358,7.822304893387042) -- (5.8008836153653895,7.793546604714431) -- (5.862082044196485,7.764916689758026) -- (5.923185213118694,7.736435330076783) -- (5.984188745910069,7.708122707229627) -- (6.045088266348661,7.6799990027755065) -- (6.105879398212519,7.652084398273363) -- (6.2271189913282345,7.596963215360754) -- (6.347872515481622,7.542920610963326) -- (6.468104960897081,7.490118037552599) -- (6.587781317799018,7.4387169476001) -- (6.7068667156611905,7.388878505003496) -- (6.825347943146864,7.34071898814369) -- (6.94325672738761,7.2942615469703265) -- (7.060630383607172,7.249517750925587) -- (7.177506227029298,7.206499169451654) -- (7.2939215728777285,7.165217371990714) -- (7.4099137363762075,7.125683927984949) -- (7.5255200327484815,7.087910406876545) -- (7.640777777218293,7.051908378107682) -- (7.755724285009386,7.017689411120548) -- (7.8703968713455055,6.985265075357326) -- (7.9848328514503955,6.954646940260195) -- (8.099069540547797,6.925846575271347) -- (8.213144253861456,6.89887554983296) -- (8.32709430661512,6.873745433387217) -- (8.440957014032529,6.850467795376306) -- (8.554769691337427,6.829054205242408) -- (8.668569653753561,6.809516232427709) -- (8.782394216504672,6.791865446374391) -- (8.896280694814509,6.7761134165246375) -- (9.010266403906812,6.762271712320634) -- (9.124388659005321,6.750351903204561) -- (9.23868477533379,6.740365558618608) -- (9.353192068115955,6.7323242480049545) -- (9.467947852575561,6.726239540805784) -- (9.582989443936356,6.722123006463283) -- (9.698354157422084,6.719986214419634) -- (9.814079308256485,6.719840734117019) -- (9.930202211663305,6.7216981349976255) -- (10.04676018286629,6.725569986503634) -- (10.163790537089177,6.731467858077231) -- (10.281330589555722,6.739403319160599) -- (10.39941765548966,6.74938793919592) -- (10.518089050114735,6.761433287625381) -- (10.637382088654698,6.7755509338911635) -- (10.757334086333287,6.791752447435453) -- (10.877982358374245,6.810049397700434) -- (10.999364220001322,6.830453354128286) -- (11.06034190791707,6.841449075103974) -- (11.121516986438257,6.8529758861611985) -- (11.182894119967838,6.86503523348023) -- (11.244477972908795,6.877628563241349) -- (11.306273209664083,6.890757321624824) -- (11.368284494636681,6.904422954810926) -- (11.430516492229549,6.918626908979932) -- (11.49297386684566,6.933370630312114) -- (11.55566128288798,6.948655564987741) -- (11.618583404759473,6.9644831591870915) -- (11.681744896863119,6.980854859090436) -- (11.745150423601867,6.9977721108780475) -- (11.808804649378704,7.0152363607302) -- (11.872712238596588,7.033249054827165) -- (11.936877855658489,7.051811639349216) -- (12.001306164967374,7.070925560476625) -- (12.066001830926213,7.090592264389668) -- (12.130969517937976,7.110813197268614) -- (12.196213108098434,7.131588763588439) -- (12.26172237771749,7.152900584828899) -- (12.32747486894894,7.174713991701168) -- (12.393447718344703,7.196993774824676) -- (12.45961806245667,7.219704724818797) -- (12.525963037836789,7.24281163230296) -- (12.59245978103695,7.266279287896566) -- (12.659085428609089,7.2900724822190455) -- (12.72581711710508,7.314156005889785) -- (12.792631983076873,7.338494649528201) -- (12.859507163076358,7.363053203753751) -- (12.926419793655455,7.387796459185765) -- (12.993347011366097,7.4126892064437016) -- (13.060265952760147,7.437696236146962) -- (13.127153754389582,7.462782338914991) -- (13.193987552806291,7.487912305367161) -- (13.260744484562167,7.513050926122901) -- (13.327401686209129,7.538162991801585) -- (13.393936294299095,7.563213293022699) -- (13.460325445384,7.588166620405616) -- (13.526546276015765,7.612987764569766) -- (13.59257592274625,7.637641516134465) -- (13.658391522127374,7.6620926657192) -- (13.723970210711101,7.686306003943429) -- (13.789289125049336,7.710246321426553) -- (13.854325401693913,7.733878408787888) -- (13.919056177196879,7.757167056646949) -- (13.98345858811004,7.780077055623053) -- (14.047509770985343,7.802573196335686) -- (14.111186862374709,7.8246202694042495) -- (14.174466998830056,7.846183065448116) -- (14.237327316903276,7.867226375086744) -- (14.299744953146302,7.887714988939507) -- (14.361697044111011,7.907613697625834) -- (14.423160726349366,7.9268872917651265) -- (14.484113136413242,7.945500561976814) -- (14.604392686225282,7.98060529309501) -- (14.722352785962428,8.012646215935689) -- (14.83781052803998,8.041341655454175) -- (14.950583004873295,8.066409936605709) -- (15.060487308877626,8.087569384345585) -- (15.167340532468302,8.104538323629157) -- (15.270959768060592,8.117035079411636) -- (15.371162108069868,8.12477797664846) -- (15.46776464491137,8.127485340294669) -- (15.56058447100051,8.124875495305844) -- (15.64943867875246,8.116666766637053) -- (15.734144360582647,8.102577479243735) -- (15.814518608906312,8.082325958081128) -- (15.890378516138753,8.055630528104473) -- (15.961541174695355,8.022209514269178) -- (16.027826977932136,7.981788163315315) -- (16.08916513491451,7.934319899035671) -- (16.197247147193167,7.8191714656000215) -- (16.26490089416211,7.71765680915658) -- (16.32238707026613,7.6030942331222775) -- (16.34738600623001,7.541166902538947) -- (16.369925473861173,7.4762739047944535) -- (16.39003294795411,7.408514010800985) -- (16.40773590330329,7.337985991470731) -- (16.42306181470326,7.264788617715766) -- (16.436038156948513,7.189020660448392) -- (16.446692404833584,7.110780890580855) -- (16.455052033152896,7.030168079025117) -- (16.46114451670101,6.947280996693479) -- (16.464997330272425,6.86221841449813) -- (16.466637948661656,6.775079103351317) -- (16.46609384666317,6.685961834165113) -- (16.463392499071517,6.594965377851821) -- (16.458561380681132,6.502188505323403) -- (16.451627966286622,6.407729987492331) -- (16.442619730682367,6.311688595270454) -- (16.431564148663,6.2141630995703565) -- (16.41848869502286,6.115252271303831) -- (16.40342084455665,6.015054881383293) -- (16.386388072058708,5.913669700720874) -- (16.367417852323655,5.811195500228763) -- (16.346537660145884,5.707731050819149) -- (16.323774970320017,5.603375123404106) -- (16.299157257640445,5.498226488895995) -- (16.272711996901734,5.392383918206889) -- (16.24446666289839,5.285946182248864) -- (16.21444873042492,5.17901205193445) -- (16.182685674275774,5.071680298175494) -- (16.149204969245574,4.9640496918842985) -- (16.114034090128655,4.8562190039729956) -- (16.07720051171964,4.748287005353944) -- (16.03873170881303,4.640352466939106) -- (15.998655156203284,4.532514159640897) -- (15.956998328684904,4.424870854371164) -- (15.913788701052425,4.3175213220423245) -- (15.869053748100328,4.210564333566623) -- (15.822820944623146,4.104098659856021) -- (15.775117765415331,3.998223071822707) -- (15.725971685271475,3.8930363403792114) -- (15.675410178985999,3.788637236437438) -- (15.623460721353439,3.6851245309096328) -- (15.570150787168274,3.5825969947078136) -- (15.515507851225038,3.4811533987444534) -- (15.459559388318183,3.3808925139314567) -- (15.40233287324233,3.2819131111812965) -- (15.34385578079187,3.1843139614058202) -- (15.28415558576134,3.088193835517444) -- (15.223259762945247,2.9936515044284135) -- (15.161195787138098,2.9007857390507468) -- (15.0979911331344,2.8096953102965188) -- (15.033673275728631,2.7204789890782024) -- (14.9682696897153,2.6332355463078727) -- (14.901807849888911,2.548063752897491) -- (14.83431523104403,2.465062379759587) -- (14.765819307975107,2.384330197806179) -- (14.696347555476649,2.305965977949512) -- (14.625927448343162,2.230068491101605) -- (14.554586461369098,2.1567365081748164) -- (14.482352069349048,2.0860688000812786) -- (14.409251747077548,2.018164137733237) -- (14.335312969348934,1.953121292042539) -- (14.260563210957827,1.8910390339217713) -- (14.185029946698734,1.8320161342830659) -- (14.108740651366134,1.776151364038384) -- (14.031722799754476,1.7235434941000847) -- (13.95400386665844,1.6742912953803568) -- (13.875611326872303,1.6284935387911048) -- (13.796572655190715,1.586248995244972) -- (13.716915326408099,1.5476564356538063) -- (13.636666815319103,1.512814630929796) -- (13.555854596718063,1.4818223519853007) -- (13.474505519662245,1.4547573029509522) -- (13.392641877686742,1.4315438165864407) -- (13.31028394924482,1.4120383836249175) -- (13.227452003843524,1.396097193599985) -- (13.144166310989789,1.3835764360450185) -- (13.060447140190433,1.3743323004940748) -- (12.976314760952476,1.3682209764800746) -- (12.891789442782851,1.3650986535373022) -- (12.806891455188492,1.364821521198678) -- (12.721641067676245,1.3672457689978046) -- (12.636058549753187,1.3722275864691937) -- (12.550164170926195,1.3796231631455385) -- (12.463978200702115,1.3892886885602138) -- (12.377520908588025,1.401080352247277) -- (12.290812564090828,1.414854343740103) -- (12.203873436717316,1.430466852572522) -- (12.116723795974622,1.4477740682779086) -- (12.029383911369564,1.4666321803898654) -- (11.941874052409133,1.4868973784413129) -- (11.854214488600178,1.5084258519669902) -- (11.766425489449801,1.5310737904998177) -- (11.678527324464739,1.554697383573398) -- (11.590540263152036,1.5791528207211059) -- (11.502484575018599,1.604296291477226) -- (11.414380529571417,1.6299839853742242) -- (11.32624839631734,1.6560720919468395) -- (11.238108444763327,1.682416800727765) -- (11.149980944416342,1.708874301251285) -- (11.061886164783346,1.7353007830500928) -- (10.973844375371158,1.7615524356584729) -- (10.885875845686883,1.7874854486098002) -- (10.798000845237311,1.8129560114374499) -- (10.710239643529434,1.837820313674797) -- (10.622612510070184,1.8619345448563536) -- (10.535139714366522,1.8851548945141303) -- (10.447841525925327,1.9073375521833213) -- (10.360738214253558,1.928338707396847) -- (10.273850048858122,1.9480145496883097) -- (10.187197299246009,1.9662212685910845) -- (10.10080023492418,1.9828150536383191) -- (10.014679125399482,1.9976520943645255) -- (9.928851158013742,2.0106054866390934) -- (9.843312762439382,2.021662186594739) -- (9.758051800740645,2.030856145531203) -- (9.673056108053544,2.038221462454459) -- (9.588313519514287,2.04379223636991) -- (9.503811870259142,2.047602566283416) -- (9.41953899542412,2.049686551200267) -- (9.335482730145486,2.0500782901260948) -- (9.251630909559452,2.0488118820665875) -- (9.167971368802,2.0459214260274905) -- (9.084491943009454,2.0414410210139806) -- (9.001180467317965,2.0354047660320873) -- (8.91802477686366,2.027846760087101) -- (8.835012706782749,2.0188011021846535) -- (8.752132092211326,2.0083018913304898) -- (8.669370768285603,1.9963832265299573) -- (8.58671657014176,1.983079206788858) -- (8.504157332915952,1.9684239311128806) -- (8.421680891744302,1.9524514985073722) -- (8.33927508176302,1.9351960079779076) -- (8.25692773810826,1.9166915585304025) -- (8.17462669591626,1.8969722491700907) -- (8.09235979032303,1.8760721789027173) -- (8.010114856464895,1.8540254467338002) -- (7.927879729477922,1.8308661516691416) -- (7.845642244498293,1.8066283927142024) -- (7.763390236662218,1.7813462688743869) -- (7.681111541105793,1.7550538791557813) -- (7.598793992965284,1.7277853225632214) -- (7.51642542737676,1.6995746981030209) -- (7.433993679476401,1.6704561047804702) -- (7.351486584400476,1.6404636416010874) -- (7.268891977284994,1.6096314075706744) -- (7.1861976932662515,1.5779935016945217) -- (7.103391567480372,1.5455840229784883) -- (7.020461435063453,1.5124370704282057) -- (6.9373951311517885,1.4785867430490214) -- (6.854180490881447,1.4440671398466236) -- (6.770805349388638,1.4089123598268145) -- (6.687257541809487,1.3731565019948277) -- (6.603524903280174,1.3368336653564654) -- (6.519595268936911,1.2999779489173022) -- (6.435456473915849,1.262623451682913) -- (6.351096353353171,1.2248042726588722) -- (6.266502742384944,1.1865545108509252) -- (6.181663476147435,1.1479082652643058) -- (6.09656638977674,1.1088996349049296) -- (6.011199318409069,1.0695627187783714) -- (5.92555009718069,1.0299316158899785) -- (5.839606561227555,0.9900404252456667) -- (5.753356545685932,0.949923245850897) -- (5.666787885691974,0.9096141767110169) -- (5.579888416381891,0.8691473168320556) -- (5.492645972891864,0.8285567652192469) -- (5.405048390358047,0.7878766208782793) -- (5.3170835039164785,0.747140982814841) -- (5.228739148703482,0.7063839500346205) -- (5.140003159855098,0.6656396215428515) -- (5.050863372507621,0.6249420963455634) -- (4.961315032767516,0.5843332708387834) -- (4.871396207428688,0.543900094849846) -- (4.781160308471499,0.5037456635174067) -- (4.690660767513009,0.4639730926386392) -- (4.599951016168461,0.42468549801117206) -- (4.509084486056054,0.3859859954335434) -- (4.418114608790347,0.3479777007020175) -- (4.327094815988971,0.3107637296153598) -- (4.2360785392676235,0.2744471979709715) -- (4.145119210243138,0.23913122156625377) -- (4.05427026053178,0.2049189161995173) -- (3.9635851217499294,0.17191339766748115) -- (3.8731172255136244,0.14021778176822863) -- (3.7829200034398127,0.10993518430007043) -- (3.6930468871445328,0.08116872105972561) -- (3.6035513082447324,0.054021507845959604) -- (3.514486698355654,0.028596660455718848) -- (3.4259064890949276,0.00499729468640453) -- (3.3378641120781367,-0.016673473663445293) -- (3.2504129989221155,-0.03631252879665681) -- (3.1636065812429024,-0.05381675491560145) -- (3.077498290657104,-0.06908303622265066) -- (2.9921415587810998,-0.08200825691972113) -- (2.907589817231269,-0.09248930120986643) -- (2.823896497624105,-0.10042305329432111) -- (2.741115031575873,-0.10570639737682086) -- (2.6592988507027258,-0.108236217659055) -- (2.5785013866214967,-0.10790939834271285) -- (2.498776070948338,-0.10462282363084796) -- (2.420176335299857,-0.09827337772514966) -- (2.34275561129175,-0.0887579448286715) -- (2.266567330541534,-0.0759734091431028) -- (2.1916649246646784,-0.05981665487149712) -- (2.1181018252781314,-0.040184566215089035) -- (2.045931463997931,-0.016974027377159473) -- (1.9752062428601675,0.00991621615639815) -- (1.905924685869195,0.04049001783096173) -- (1.8380060555123237,0.07460794157759665) -- (1.7713632696431887,0.11211908145469351) -- (1.7059092461140608,0.1528725315210977) -- (1.6415569027783476,0.19671738583747356) -- (1.57821915748832,0.24350273845993797) -- (1.5158089280980676,0.2930776834482458) -- (1.4542391324582695,0.3452913148612424) -- (1.3934226884248346,0.3999927267586827) -- (1.273701526582954,0.5162552682368187) -- (1.2146226444806416,0.5775145859370241) -- (1.1559487853955943,0.6406580603552356) -- (1.0975928671791735,0.7055347855502987) -- (1.0394678076859236,0.7719938555819681) -- (0.9814865247676607,0.8398843645072702) -- (0.9235619362771104,0.9090554063877789) -- (0.865606960068817,0.9793560752787016) -- (0.8075345139941419,1.0506354652425216) -- (0.7492575159053558,1.1227426703353558) -- (0.6906888836574581,1.1955267846160496) -- (0.6317415351022646,1.2688369021443577) -- (0.5723283880915915,1.3425221169777615) -- (0.5123623604818022,1.4164315231782894) -- (0.4517563701219842,1.4904142148002393) -- (0.3904233348666821,1.564319285905185) -- (0.32827617256907615,1.6379958305519722) -- (0.26522780108234656,1.7112929427994459) -- (0.2011911382587641,1.7840597167032684) -- (0.13607910195105433,1.8561452463254682) -- (0.06980461001239746,1.9273986257239812) -- (0.0022805802946095355,1.9976689489567434) -- (-0.06654272044579557,2.0668370679746886) -- (-0.1365498544428192,2.1349242801776427) -- (-0.20757973957915965,2.2019906945012053) -- (-0.2794712921941027,2.2680964211982655) -- (-0.3520634286205677,2.333301570520348) -- (-0.42519506519874994,2.397666252718068) -- (-0.4987051182652067,2.4612505780456786) -- (-0.5724325041555858,2.524114656753113) -- (-0.6462161392064445,2.586318599092806) -- (-0.7198949397588876,2.647922515315827) -- (-0.7933078221440155,2.7089865156755195) -- (-0.8662937027029329,2.769570710422272) -- (-0.9386914977721972,2.829735209808973) -- (-1.010340123687456,2.889540124087148) -- (-1.081078496787086,2.949045563507866) -- (-1.150745533407644,3.008311638323562) -- (-1.2191801498838686,3.0673984587857603) -- (-1.2862212625550455,3.12636613514735) -- (-1.3517077877577321,3.1852747776589467) -- (-1.4154786418312142,3.2441844965725295) -- (-1.477372741108411,3.3031554021405327) -- (-1.594886340625635,3.4215212142458995) -- (-1.7029599170118672,3.540853095989519) -- (-1.800304800957747,3.661631929386317) -- (-1.8445512741100174,3.722714228583868) -- (-1.8856323231630086,3.7843385964505387) -- (-1.9233868644550967,3.846565143237626) -- (-1.9576538143219295,3.909453979196428) -- (-1.988272089099155,3.973065214580515) -- (-2.0150883129622343,4.037456489640249) -- (-2.0380740038103795,4.102645422036403) -- (-2.0573019607168135,4.1686171735569815) -- (-2.072847902963531,4.235355970200999) -- (-2.084787549827979,4.3028460379661055) -- (-2.0931966205939716,4.371071602851657) -- (-2.0981508345412294,4.4400168908555315) -- (-2.099725910948564,4.509666127976516) -- (-2.097997569099789,4.580003540212942) -- (-2.0930415282723516,4.651013353563144) -- (-2.0849335077477917,4.7226797940259075) -- (-2.0737492268076494,4.794987087599679) -- (-2.059564404730736,4.867919460282678) -- (-2.0424547607985915,4.941461138073578) -- (-2.0224960142918462,5.015596346970824) -- (-1.9997638844920402,5.090309312972977) -- (-1.9743340906766207,5.165584262078369) -- (-1.9462823521294013,5.241405420285446) -- (-1.9156843881305576,5.317757013592654) -- (-1.882615917957537,5.394623267998895) -- (-1.8471526608950626,5.471988409501932) -- (-1.8093703362205815,5.549836664101008) -- (-1.769344663215179,5.6281522577939995) -- (-1.7271513611608498,5.706919416579694) -- (-1.6828661493377695,5.786122366456425) -- (-1.6365647470252043,5.865745333422751) -- (-1.588322873505149,5.945772543477119) -- (-1.5382162480568695,6.026188222617975) -- (-1.4863205899614513,6.106976596843765) -- (-1.432711618500889,6.188121892153276) -- (-1.3774650529535393,6.2696083345445) -- (-1.3206566125995778,6.351420150016338) -- (-1.2623620167223635,6.433541564566781) -- (-1.2026569846007078,6.515956804195071) -- (-1.141617235515696,6.59865009489863) -- (-1.0793184887456846,6.68160566267693) -- (-1.015836463575397,6.764807733527846) -- (-0.9512468792818254,6.848240533450166) -- (-0.8856254551469647,6.931888288442224) -- (-0.8190479104509905,7.01573522450235) -- (-0.7515899644754427,7.099765567629333) -- (-0.6833273364991328,7.1839635438215055) -- (-0.6143357458036007,7.268313379077199) -- (-0.5446909116685674,7.352799299395315) -- (-0.4744685533764823,7.437405530773617) -- (-0.4037443902070663,7.522116299211234) -- (-0.3325941414391309,7.606915830706612) -- (-0.2610935263551255,7.691788351257628) -- (-0.18931826423659004,7.776718086863525) -- (-0.11734407436188121,7.86168926352218) -- (-0.04524667600981047,7.946686107232608) -- (-0.009175435932775144,7.989189224856432) -- (-1.5688003577452037E-4,7.999815159644868) -- (-1.596469428477576E-5,7.999981189960522) -- (-7.157485470088432E-6,7.999991566855215) -- (-2.7538817448657937E-6,7.999996755302732) -- (-5.520792001334485E-7,7.99999934952632) -- (-1.628450263524428E-9,7.99999999808233)(0.9859737910219337,0.8345788454316789);
\draw [color=qqwwzz](1.8,7.8) node[anchor=north west] {$V_\varphi \subset \mathbb{R}^3$};
\draw (3.5,4.9) node[anchor=north west] {$\Coss$};
\draw [color=xfqqff](-1.2,1.1) node[anchor=north west] {$\boldsymbol{\varphi}$};
\draw (-1.8942371003009975,-0.7211750651048283) node[anchor=north west] {$x$};
\draw (-4.3,1.8) node[anchor=north west] {$y$};
\draw [color=qqqqff](11.8,8.4) node[anchor=north west] {$\mathbf{t}_\varphi$};
\draw [color=qqqqff](14.4,6.2) node[anchor=north west] {$\mathbf{f}_\varphi$};
\draw [color=qqqqff](10.2,5.2) node[anchor=north west] {$\vb{m}_\varphi$};
\draw [color=qqqqff](14.5,1.9) node[anchor=north west] {$\boldsymbol{\mu}_\varphi$};
\draw [color=qqwwzz](6,8.4) node[anchor=north west] {$\partial V_\varphi$};
\begin{scriptsize}
\draw [fill=black] (-4.,-1.) circle (2.5pt);
\draw [fill=black] (3.,4.) circle (2.5pt);
\end{scriptsize}
\end{tikzpicture}

%% file: figs/inter.tex
\definecolor{xfqqff}{rgb}{0.4980392156862745,0.,1.}
\definecolor{wqwqwq}{rgb}{0.3764705882352941,0.3764705882352941,0.3764705882352941}
\definecolor{qqwwzz}{rgb}{0.,0.4,0.6}
\definecolor{qqqqff}{rgb}{0.,0.,1.}
\begin{tikzpicture}[scale=0.45,line cap=round,line join=round,>=triangle 45,x=1.0cm,y=1.0cm]
\draw [line width=1.5pt,color=qqqqff] (5.065327132806232,1.2100310945710084)-- (8.678439555960983,7.708031731562509);
\draw[line width=0.7pt,color=qqwwzz,fill=qqwwzz,fill opacity=0.05, dashed] (5.0653271328062575,1.210031094571006) -- (4.944529995837541,1.1360418919808353) -- (4.823754075012105,1.0620963717240905) -- (4.703020586473231,0.9882382161341969) -- (4.5823507463642015,0.9145111075445806) -- (4.461765770828295,0.8409587282886671) -- (4.341286876008794,0.7676247606998823) -- (4.220935278048978,0.6945528871116516) -- (4.100732193092131,0.6217867898574011) -- (3.9806988372815315,0.5493701512705562) -- (3.8608564267604617,0.4773466536845429) -- (3.7412261776722016,0.40575997943278663) -- (3.6218293061600333,0.33465381084871315) -- (3.5026870283672373,0.26407183026574843) -- (3.383820560437095,0.19405772001731791) -- (3.265251118512887,0.12465516243684749) -- (3.1469999187378948,0.05590783985776282) -- (3.029088177255399,-0.012140565386510316) -- (2.911537110208681,-0.07944637096254636) -- (2.7943679337410208,-0.14596589453691955) -- (2.677601863995701,-0.2116554537762041) -- (2.5612601171160017,-0.27647136634697445) -- (2.445363909245204,-0.3403699499158046) -- (2.3299344565265896,-0.403307522149269) -- (2.2149929751034385,-0.4652404007139421) -- (2.1005606811190325,-0.5261249032763979) -- (1.986658790716652,-0.5859173475032109) -- (1.8733085200395787,-0.644574051060955) -- (1.7605310852310936,-0.7020513316162051) -- (1.648347702434477,-0.7583055068355347) -- (1.5367795877930108,-0.8132928943855191) -- (1.4258479574499758,-0.8669698119327318) -- (1.3155740275486525,-0.919292577143747) -- (1.2059790142323226,-0.9702175076851391) -- (1.0970841336442674,-1.019700921223483) -- (0.9889106019277669,-1.0676991354253524) -- (0.8814796352261025,-1.114168467957322) -- (0.7748124496825559,-1.1590652364859655) -- (0.6689302614404076,-1.2023457586778576) -- (0.5638542866429388,-1.2439663521995725) -- (0.45960574143342914,-1.2838833347176841) -- (0.3562058419551626,-1.3220530238987673) -- (0.2536758043514169,-1.3584317374093962) -- (0.1520368447654752,-1.392975792916145) -- (0.05131017934061877,-1.4256415080855884) -- (-0.04848297577987282,-1.4563852005842999) -- (-0.14732140445271735,-1.4851631880788547) -- (-0.2451838905346353,-1.511931788235826) -- (-0.34204921788234355,-1.5366473187217888) -- (-0.43789617035256256,-1.5592660972033174) -- (-0.532703531802011,-1.5797444413469859) -- (-0.6264500860874067,-1.5980386688193686) -- (-0.719114617065471,-1.61410509728704) -- (-0.8106759085929207,-1.6279000444165734) -- (-0.9011127445264755,-1.6393798278745448) -- (-0.9904039087228549,-1.648500765327527) -- (-1.0785281850387767,-1.6552191744420952) -- (-1.1654643573309613,-1.6594913728848226) -- (-1.2511912094561266,-1.661273678322285) -- (-1.335687525270992,-1.6605224084210555) -- (-1.4189320886322756,-1.657193880847709) -- (-1.5009036833966976,-1.6512444132688193) -- (-1.5815810934209766,-1.6426303233509603) -- (-1.6609431025618315,-1.6313079287607075) -- (-1.738968494675981,-1.617233547164635) -- (-1.8156360536201435,-1.600363496229316) -- (-1.8909245632510396,-1.5806540936213263) -- (-1.9648128074253872,-1.5580616570072388) -- (-2.0372795699999067,-1.5325425040536285) -- (-2.108303634831314,-1.5040529524270698) -- (-2.1778637857763297,-1.4725493197941364) -- (-2.2459388066916732,-1.437987923821403) -- (-2.312507481434065,-1.4003250821754438) -- (-2.4410409278268608,-1.315520332530145) -- (-2.5632944607861567,-1.2177856987428584) -- (-2.6792285379237892,-1.1069451340746639) -- (-2.7892063405602094,-0.9833590295639958) -- (-2.842099406213226,-0.9169313544403757) -- (-2.8936687239682524,-0.8474912398433112) -- (-2.9439694007344066,-0.7750966676018543) -- (-2.993056543420778,-0.6998056195450566) -- (-3.0409852589364803,-0.6216760775019807) -- (-3.087810654190619,-0.540766023301682) -- (-3.1335878360922997,-0.4571334387732122) -- (-3.178371911550631,-0.3708363057456232) -- (-3.222217987474707,-0.28193260604798454) -- (-3.265181170773648,-0.19048032150934446) -- (-3.3073165683565513,-0.09653743395875125) -- (-3.348679287132523,-1.6192522527447295E-4) -- (-3.3893244340106747,0.09858822286204116) -- (-3.429307115900098,0.1996550284741332) -- (-3.4686824397099194,0.3029805097819427) -- (-3.50750551234923,0.4085066849564285) -- (-3.545831440727138,0.5161755721685211) -- (-3.583715331752739,0.6259291895891756) -- (-3.6212122923351693,0.7377095553893156) -- (-3.6583774293835027,0.8514586877399068) -- (-3.695265849806855,0.9671186048118869) -- (-3.7319326605143424,1.0846313247762076) -- (-3.768432968415052,1.2039388658037922) -- (-3.8048218804181033,1.3249832460655995) -- (-3.841154503432602,1.4477064837325742) -- (-3.877485944367642,1.5720505969756537) -- (-3.913871310132347,1.69795760396579) -- (-3.9503657076358074,1.8253695228739133) -- (-3.987024243787129,1.9542283718709967) -- (-4.02387002138309,2.0844504585205357) -- (-4.042262180320762,2.149959667353045) -- (-4.060526806204479,2.2156312816956927) -- (-4.078582321307607,2.2813869842986207) -- (-4.096347147903572,2.347148457912084) -- (-4.1137397082652285,2.4128373852860534) -- (-4.130678424666229,2.478375449170784) -- (-4.162868014678992,2.6086857174731506) -- (-4.192263296127891,2.7374527248204856) -- (-4.2182116471991264,2.864049933213977) -- (-4.240060446078843,2.987850804654755) -- (-4.257157070953184,3.1082288011443495) -- (-4.268848900008521,3.2245573846839477) -- (-4.274483311430771,3.336210017274567) -- (-4.273407683406248,3.4425601609177363) -- (-4.264969394121209,3.5429812776146434) -- (-4.2485158217617425,3.6368468293664193) -- (-4.223394344514048,3.7235302781744792) -- (-4.188952340564384,3.802405086039954) -- (-4.1445371880988375,3.872844714964259) -- (-4.023176950364672,3.9859122839939403) -- (-3.9450721831049123,4.027431836178721) -- (-3.8562478861226737,4.059862976411296) -- (-3.7587309465970975,4.085242595141153) -- (-3.6545620129027157,4.105621261429405) -- (-3.5457817334167885,4.123049544336936) -- (-3.4344307565149848,4.1395780129264494) -- (-3.322549730572746,4.157257236257919) -- (-3.212179303967332,4.178137783394277) -- (-3.1053601250728207,4.20427022339527) -- (-3.0041328422669267,4.237705125324283) -- (-2.9105128762399204,4.280466216241393) -- (-2.825413048872008,4.333404069299945) -- (-2.748239080523433,4.395765716978758) -- (-2.6782880782956795,4.466682628356807) -- (-2.6148571492902306,4.545286272511589) -- (-2.557243400609366,4.6307081185220795) -- (-2.504743939354114,4.722079635466457) -- (-2.456655872626527,4.818532292422901) -- (-2.4122763075284297,4.919197558469591) -- (-2.370902351161078,5.0232069026851605) -- (-2.331831110626183,5.129691794147675) -- (-2.294359693025683,5.2377837019357685) -- (-2.2577852054617438,5.34661409512762) -- (-2.221404755034712,5.455314442801523) -- (-2.1845154488473213,5.563016214035542) -- (-2.1464143940007148,5.668850877908085) -- (-2.1063986975971716,5.771949903497784) -- (-2.0637654667378342,5.871444759882593) -- (-2.017811808524641,5.966466916140803) -- (-1.9678699167544664,6.056191759280296) -- (-1.9137331537385194,6.14037191978457) -- (-1.8555198367960202,6.219166773509414) -- (-1.7933550762993775,6.292744199158847) -- (-1.727363982621398,6.361272075436716) -- (-1.6576716661348314,6.424918281046871) -- (-1.5844032372121433,6.483850694693274) -- (-1.5076838062262254,6.538237195079887) -- (-1.4276384835496003,6.588245660910502) -- (-1.3443923795550177,6.634043970889081) -- (-1.2580706046151704,6.675800003719473) -- (-1.1687982691026377,6.713681638105697) -- (-1.0767004833902547,6.747856752751659) -- (-0.9819023578505721,6.778493226361206) -- (-0.8845290028561976,6.805758937638188) -- (-0.7847055287800231,6.829821765286624) -- (-0.6825570459945425,6.850849588010192) -- (-0.578208664872534,6.869010284513138) -- (-0.47178549578666207,6.884471733499254) -- (-0.36341264910959126,6.897401813672332) -- (-0.25321523521387235,6.90796840373639) -- (-0.14131836447234036,6.916339382395222) -- (-0.0278471472576598,6.922682628352732) -- (0.0870733060575617,6.927166020313109) -- (0.20331788510060278,6.929957436979976) -- (0.32076147949879896,6.931224757057294) -- (0.43927897887948575,6.931135859249025) -- (0.5587452728699986,6.929858622258962) -- (0.6790352510975595,6.927560924791123) -- (0.8000238031896743,6.924410645549358) -- (0.9215858187735648,6.920575663237628) -- (1.0435961874765098,6.916223856559782) -- (1.1659297989260153,6.911523104219896) -- (1.2884615427490758,6.906641284921648) -- (1.4110663085733108,6.901746277368943) -- (1.5336189860259424,6.897005960265915) -- (1.655994464734306,6.892588212316241) -- (1.77806763432568,6.888660912224054) -- (1.8997133844275709,6.885391938692919) -- (2.0208066046670297,6.882949170427082) -- (2.1412221846714488,6.881500486130335) -- (2.2608350140682774,6.881213764506583) -- (2.379519982484794,6.882256884259675) -- (2.4971519795483346,6.88479772409363) -- (2.613605894886234,6.889004162712126) -- (2.728756618125601,6.89504407881941) -- (2.84250254659338,6.903051152439389) -- (2.954870647108635,6.912972022453346) -- (3.065931619182095,6.9246897060212405) -- (3.1757561986618157,6.938087167439164) -- (3.2844151213962505,6.953047371002867) -- (3.391979123233341,6.969453281008725) -- (3.4985189400214836,6.987187861752545) -- (3.6041053076089327,7.006134077530646) -- (3.7088089618436584,7.026174892639062) -- (3.812700638574114,7.047193271373658) -- (3.915851073648412,7.069072178030979) -- (4.018331002914664,7.0916945769066615) -- (4.12021116222121,7.114943432297025) -- (4.221562287416305,7.13870170849799) -- (4.322455114348031,7.162852369805819) -- (4.422960378864701,7.187278380516489) -- (4.523148816814398,7.211862704926148) -- (4.623091164045491,7.236488307330717) -- (4.722858156406119,7.2610381520263445) -- (4.822520529744452,7.2853952033092355) -- (4.922149019908716,7.309442425475481) -- (5.021814362747136,7.333062782820889) -- (5.121587294107883,7.356139239641777) -- (5.221538549839181,7.378554760234238) -- (5.321738865789484,7.400192308894248) -- (5.422258977806564,7.420934849918012) -- (5.523165749574229,7.440672302707213) -- (5.624478338676678,7.459380275419562) -- (5.726183796027652,7.477092045580946) -- (5.828268574446152,7.493841965003611) -- (5.93071912675121,7.509664385499519) -- (6.033521905761813,7.524593658880889) -- (6.136663364297046,7.538664136959682) -- (6.240129955175842,7.5519101715481725) -- (6.343908131217287,7.564366114458295) -- (6.447984345240339,7.57606631750221) -- (6.552345050064073,7.5870451324920225) -- (6.65697669850746,7.597336911239978) -- (6.7618657433895315,7.606976005557868) -- (6.866998637529274,7.615996767258082) -- (6.972361833745804,7.6244335481525525) -- (7.077941784858009,7.632320700053384) -- (7.183724943685004,7.63969257477271) -- (7.289697763045734,7.646583524122747) -- (7.395846695759271,7.653027899915372) -- (7.502158194644604,7.659060053962861) -- (7.608618712520748,7.664714338077232) -- (7.715214702206719,7.670025104070618) -- (7.8219326165215755,7.675026703755094) -- (7.928758908284262,7.679753488942794) -- (8.035680030313866,7.684239811445764) -- (8.14268243542936,7.688520023076109) -- (8.249752576449747,7.692628475645989) -- (8.356876906194085,7.69659952096751) -- (8.464041877481375,7.70046751085269) -- (8.571233943130633,7.704266797113661) -- (8.678439555960875,7.7080317315624995)(5.0653271328062575,1.210031094571006)(5.065327132806211,1.2100310945709778) -- (5.120913477777843,1.3100003351401084) -- (5.176499822749455,1.4099695757092086) -- (5.232086167721066,1.5099388162783085) -- (5.287672512692678,1.6099080568474085) -- (5.343258857664289,1.7098772974165084) -- (5.398845202635901,1.8098465379856086) -- (5.454431547607513,1.9098157785547085) -- (5.510017892579124,2.0097850191238087) -- (5.565604237550736,2.1097542596929086) -- (5.621190582522347,2.2097235002620086) -- (5.676776927493959,2.3096927408311085) -- (5.73236327246557,2.409661981400209) -- (5.787949617437182,2.5096312219693084) -- (5.843535962408794,2.609600462538409) -- (5.899122307380405,2.7095697031075088) -- (5.954708652352017,2.8095389436766087) -- (6.010294997323628,2.9095081842457087) -- (6.06588134229524,3.0094774248148086) -- (6.121467687266851,3.109446665383909) -- (6.177054032238463,3.2094159059530085) -- (6.232640377210075,3.309385146522109) -- (6.288226722181686,3.4093543870912084) -- (6.343813067153298,3.509323627660309) -- (6.399399412124909,3.609292868229409) -- (6.454985757096521,3.7092621087985087) -- (6.510572102068132,3.809231349367609) -- (6.566158447039744,3.9092005899367086) -- (6.621744792011356,4.009169830505809) -- (6.677331136982967,4.1091390710749085) -- (6.732917481954578,4.209108311644009) -- (6.78850382692619,4.309077552213109) -- (6.844090171897802,4.409046792782209) -- (6.8996765168694125,4.509016033351308) -- (6.955262861841025,4.608985273920409) -- (7.010849206812636,4.708954514489508) -- (7.066435551784247,4.808923755058608) -- (7.122021896755859,4.908892995627707) -- (7.177608241727469,5.008862236196807) -- (7.233194586699081,5.108831476765906) -- (7.288780931670693,5.208800717335007) -- (7.344367276642304,5.308769957904106) -- (7.399953621613916,5.408739198473206) -- (7.455539966585526,5.508708439042305) -- (7.511126311557138,5.608677679611405) -- (7.5667126565287495,5.708646920180505) -- (7.622299001500361,5.8086161607496045) -- (7.677885346471973,5.908585401318704) -- (7.733471691443583,6.0085546418878035) -- (7.789058036415195,6.108523882456903) -- (7.844644381386806,6.2084931230260025) -- (7.9002307263584175,6.308462363595103) -- (7.95581707133003,6.408431604164202) -- (8.01140341630164,6.508400844733302) -- (8.066989761273252,6.608370085302401) -- (8.122576106244864,6.708339325871501) -- (8.178162451216474,6.8083085664406005) -- (8.233748796188085,6.908277807009701) -- (8.289335141159697,7.0082470475788) -- (8.344921486131309,7.1082162881479) -- (8.40050783110292,7.208185528716999) -- (8.456094176074531,7.308154769286099) -- (8.511680521046141,7.408124009855198) -- (8.567266866017754,7.508093250424299) -- (8.622853210989366,7.608062490993398) -- (8.678439555960978,7.708031731562498)(5.065327132806232,1.2100310945710084);
\draw[line width=0.7pt,color=wqwqwq,fill=wqwqwq,fill opacity=0.05, dashed] (5.065430500595812,1.2102169963321148) -- (5.120913477777843,1.3100003351401084) -- (5.176499822749455,1.4099695757092086) -- (5.232086167721066,1.5099388162783085) -- (5.287672512692678,1.6099080568474085) -- (5.343258857664289,1.7098772974165084) -- (5.398845202635901,1.8098465379856086) -- (5.454431547607513,1.9098157785547085) -- (5.510017892579124,2.0097850191238087) -- (5.565604237550736,2.1097542596929086) -- (5.621190582522347,2.2097235002620086) -- (5.676776927493959,2.3096927408311085) -- (5.73236327246557,2.409661981400209) -- (5.787949617437182,2.5096312219693084) -- (5.843535962408794,2.609600462538409) -- (5.899122307380405,2.7095697031075088) -- (5.954708652352017,2.8095389436766087) -- (6.010294997323628,2.9095081842457087) -- (6.06588134229524,3.0094774248148086) -- (6.121467687266851,3.109446665383909) -- (6.177054032238463,3.2094159059530085) -- (6.232640377210075,3.309385146522109) -- (6.288226722181686,3.4093543870912084) -- (6.343813067153298,3.509323627660309) -- (6.399399412124909,3.609292868229409) -- (6.454985757096521,3.7092621087985087) -- (6.510572102068132,3.809231349367609) -- (6.566158447039744,3.9092005899367086) -- (6.621744792011356,4.009169830505809) -- (6.677331136982967,4.1091390710749085) -- (6.732917481954578,4.209108311644009) -- (6.78850382692619,4.309077552213109) -- (6.844090171897802,4.409046792782209) -- (6.8996765168694125,4.509016033351308) -- (6.955262861841025,4.608985273920409) -- (7.010849206812636,4.708954514489508) -- (7.066435551784247,4.808923755058608) -- (7.122021896755859,4.908892995627707) -- (7.177608241727469,5.008862236196807) -- (7.233194586699081,5.108831476765906) -- (7.288780931670693,5.208800717335007) -- (7.344367276642304,5.308769957904106) -- (7.399953621613916,5.408739198473206) -- (7.455539966585526,5.508708439042305) -- (7.511126311557138,5.608677679611405) -- (7.5667126565287495,5.708646920180505) -- (7.622299001500361,5.8086161607496045) -- (7.677885346471973,5.908585401318704) -- (7.733471691443583,6.0085546418878035) -- (7.789058036415195,6.108523882456903) -- (7.844644381386806,6.2084931230260025) -- (7.9002307263584175,6.308462363595103) -- (7.95581707133003,6.408431604164202) -- (8.01140341630164,6.508400844733302) -- (8.066989761273252,6.608370085302401) -- (8.122576106244864,6.708339325871501) -- (8.178162451216474,6.8083085664406005) -- (8.233748796188085,6.908277807009701) -- (8.289335141159697,7.0082470475788) -- (8.344921486131309,7.1082162881479) -- (8.40050783110292,7.208185528716999) -- (8.456094176074531,7.308154769286099) -- (8.511680521046141,7.408124009855198) -- (8.567266866017754,7.508093250424299) -- (8.622853210989366,7.608062490993398) -- (8.678439555960978,7.708031731562498)(5.065430500595812,1.2102169963321148)(8.678416246898982,7.707989811388988) -- (8.793574439702656,7.710990938592333) -- (8.908695150554928,7.713946281739436) -- (9.023787515628381,7.7168938971210554) -- (9.13883736203361,7.7198299208544485) -- (9.253830516881207,7.7227504890568754) -- (9.36875280728176,7.725651737845592) -- (9.483590060345863,7.7285298033378576) -- (9.598328103184109,7.731380821650931) -- (9.712952762907086,7.734200928902069) -- (9.827449866625388,7.73698626120853) -- (9.941805241449607,7.7397329546875735) -- (10.056004714490333,7.7424371454564564) -- (10.17003411285816,7.7450949696324365) -- (10.283879263663676,7.747702563332774) -- (10.397525994017476,7.750256062674724) -- (10.51096013103015,7.752751603775548) -- (10.62416750181229,7.755185322752501) -- (10.737133933474489,7.757553355722843) -- (10.849845253127334,7.759851838803832) -- (10.962287287881423,7.762076908112726) -- (11.074445864847343,7.764224699766783) -- (11.186306811135687,7.766291349883261) -- (11.297855953857047,7.7682729945794184) -- (11.409079120122016,7.770165769972514) -- (11.519962137041182,7.771965812179804) -- (11.630490831725139,7.773669257318549) -- (11.740651031284479,7.775272241506006) -- (11.85042856282979,7.776770900859432) -- (11.959809253471668,7.778161371496088) -- (12.068778930320704,7.779439789533229) -- (12.17732342048749,7.780602291088115) -- (12.285428551082614,7.7816450122780045) -- (12.39308014921667,7.782564089220154) -- (12.50026404200025,7.783355658031823) -- (12.606966056543946,7.784015854830269) -- (12.713172019958346,7.7845408157327505) -- (12.818867759354049,7.784926676856526) -- (12.92403910184164,7.785169574318853) -- (13.02867187453171,7.785265644236989) -- (13.132751904534857,7.785211022728194) -- (13.236265018961667,7.785001845909725) -- (13.339197044922734,7.78463424989884) -- (13.44153380952865,7.7841043708127975) -- (13.543261139890005,7.783408344768856) -- (13.644364863117392,7.782542307884273) -- (13.744830806321401,7.781502396276307) -- (13.844644796612625,7.7802847460622155) -- (13.943792661101657,7.778885493359258) -- (14.042260226899085,7.777300774284692) -- (14.140033321115503,7.775526724955775) -- (14.237097770861503,7.773559481489765) -- (14.333439403247674,7.771395180003922) -- (14.429044045384611,7.769029956615503) -- (14.523897524382903,7.766459947441765) -- (14.617985667353146,7.763681288599969) -- (14.711294301405925,7.76069011620737) -- (14.803809253651835,7.757482566381229) -- (14.895516351201469,7.754054775238802) -- (14.986401421165418,7.750402878897348) -- (15.07645029065427,7.746523013474125) -- (15.165648786778622,7.7424113150863905) -- (15.25398273664906,7.738063919851404) -- (15.341437967376182,7.733476963886424) -- (15.428000306070576,7.728646583308707) -- (15.513655579842832,7.7235689142355115) -- (15.598389615803544,7.718240092784096) -- (15.682188241063304,7.712656255071719) -- (15.765037282732703,7.7068135372156386) -- (15.84692256792233,7.700708075333112) -- (15.927829923742781,7.694336005541398) -- (16.007745177304646,7.687693463957755) -- (16.086654155718517,7.680776586699442) -- (16.164542686094983,7.673581509883714) -- (16.24139659554464,7.666104369627833) -- (16.317201711178075,7.658341302049054) -- (16.391943860105883,7.650288443264638) -- (16.465608869438654,7.641941929391841) -- (16.538182566286977,7.633297896547922) -- (16.609650777761452,7.624352480850138) -- (16.679999330972663,7.6151018184157495) -- (16.749214053031203,7.605542045362013) -- (16.817280771047667,7.595669297806186) -- (16.884185312132644,7.585479711865529) -- (16.949913503396722,7.574969423657298) -- (17.077784144904562,7.552971284907149) -- (17.200779312455925,7.529643970493805) -- (17.318785622935536,7.504956569355332) -- (17.431689693228133,7.478878170429794) -- (17.53937814021845,7.451377862655256) -- (17.641737580791215,7.422424734969784) -- (17.73865463183116,7.391987876311443) -- (17.83001591022303,7.360036375618299) -- (17.915708032851544,7.326539321828417) -- (17.995617616601443,7.291465803879861) -- (18.069631278357456,7.254784910710697) -- (18.19951730342676,7.176477354462806) -- (18.304459043137328,7.091369364589266) -- (18.383549432567012,6.999213652594598) -- (18.435881406793683,6.899762929983326) -- (18.460547900895207,6.792769908259969) -- (18.456641849949442,6.67798729892905) -- (18.42325618903426,6.555167813495093) -- (18.395225047049966,6.490666964584636) -- (18.359483853227516,6.424064163462617) -- (18.31591922445164,6.355328499067097) -- (18.26441777760708,6.284429060336144) -- (18.204866129578555,6.211334936207823) -- (18.13715089725081,6.136015215620199) -- (18.06115869750857,6.058438987511336) -- (17.976779929468975,5.978576057600689) -- (17.88429962370715,5.896471019312555) -- (17.784737490945346,5.81230769734872) -- (17.679193383035667,5.726285104179295) -- (17.5687671518301,5.638602252274373) -- (17.45455864918084,5.5494581541041) -- (17.337667726939785,5.4590518221385675) -- (17.219194236959183,5.36758226884791) -- (17.10023803109111,5.275248506702239) -- (16.98189896118764,5.18224954817166) -- (16.86527687910079,5.088784405726299) -- (16.751471636682723,4.995052091836262) -- (16.6415830857854,4.901251618971683) -- (16.536711078261064,4.8075819996026805) -- (16.43795546596168,4.714242246199332) -- (16.346416100739432,4.621431371231806) -- (16.263192834446286,4.529348387170179) -- (16.189385518934344,4.438192306484572) -- (16.12609400605585,4.348162141645124) -- (16.074418147662712,4.259456905121921) -- (16.035457795607144,4.172275609385116) -- (16.010312801741136,4.086817266904793) -- (16.00008301791665,4.00328089015108) -- (16.00558643115042,3.921827709130092) -- (16.026370504643097,3.842448646969771) -- (16.061624879530115,3.7650866626761683) -- (16.110539176898527,3.689684712568024) -- (16.172303017835333,3.6161857529640784) -- (16.2461060234277,3.5445327401830795) -- (16.331137814762513,3.4746686305437464) -- (16.426588012927112,3.406536380364841) -- (16.531646239008495,3.3400789459650966) -- (16.64550211409366,3.275239283663261) -- (16.76734525926966,3.21196034977806) -- (16.89636529562378,3.1501851006282493) -- (16.96331340458218,3.119843531678711) -- (17.03175184424279,3.0898564925325616) -- (17.101579317241402,3.060216852729617) -- (17.172694526213974,3.0309174818097375) -- (17.24499617379641,3.0019512493127607) -- (17.31838296262444,2.973311024778525) -- (17.392753595334028,2.944989677746868) -- (17.468006774561132,2.916980077757657) -- (17.54404120294177,2.8892750943507153) -- (17.62075558311139,2.8618675970658742) -- (17.698048617706235,2.8347504554430074) -- (17.775819009361925,2.8079165390219316) -- (17.853965460714704,2.781358717342499) -- (17.932386674400135,2.7550698599445553) -- (18.01098135305412,2.729042836367931) -- (18.08964819931279,2.703270516152486) -- (18.16828591581188,2.677745768838051) -- (18.246793205187174,2.652461463964471) -- (18.325068770074807,2.627410471071599) -- (18.403011313110397,2.6025856596992583) -- (18.480519536929904,2.5779798993872873) -- (18.55749214416946,2.5535860596755526) -- (18.633827837464622,2.5293970101038923) -- (18.709425319451356,2.5054056202121373) -- (18.78418329276562,2.481604759540147) -- (18.85800046004337,2.457987297627753) -- (18.930775523920346,2.4345461040147853) -- (19.002407187032333,2.411274048241104) -- (19.07279415201552,2.3881639998465545) -- (19.14183512150575,2.3652088283709602) -- (19.209428798138532,2.342401403354181) -- (19.275473884550053,2.3197345943360403) -- (19.402513097252864,2.2747943024551134) -- (19.5221423807009,2.230330909046913) -- (19.633551355981552,2.186287370430165) -- (19.735929644181766,2.14260664292361) -- (19.82846686638848,2.0992316828460105) -- (19.910352643689038,2.056105446516071) -- (19.980776597170433,2.0131708902525816) -- (20.08490839730058,1.9276943494070853) -- (20.14110503465446,1.8426831984301) -- (20.15189548618349,1.7581332929836258) -- (20.11981418897409,1.674040762708671) -- (20.047395580111868,1.5904017372462436) -- (19.93717409668301,1.507212346237349) -- (19.868679288595388,1.465785054165793) -- (19.791684175773128,1.4244687193229968) -- (19.70650556285233,1.3832628579140902) -- (19.613460254468407,1.3421669861441963) -- (19.512865055257407,1.3011806202184424) -- (19.40503676985469,1.2603032763419515) -- (19.29029220289607,1.2195344707198537) -- (19.168948159017702,1.1788737195572718) -- (19.041321442854894,1.1383205390593343) -- (18.975251084115655,1.1180841366237075) -- (18.907728859043686,1.0978744454311649) -- (18.838794368218373,1.0776914050073456) -- (18.768487212219384,1.0575349548778892) -- (18.69684699162616,1.0374050345684402) -- (18.6239133070182,1.0173015836046364) -- (18.549725758974887,0.9972245415121206) -- (18.474323948075778,0.9771738478165295) -- (18.3977474749002,0.9571494420435087) -- (18.320035940027594,0.9371512637186967) -- (18.241228944037516,0.9171792523677329) -- (18.16136608750952,0.897233347516261) -- (18.08048697102305,0.8773134886899197) -- (17.998631195157316,0.8574196154143507) -- (17.915838360492216,0.8375516672151946) -- (17.83214806760708,0.8177095836180897) -- (17.747599917081004,0.7978933041486833) -- (17.662233509493888,0.7781027683326087) -- (17.576088445424944,0.7583379156955097) -- (17.489204325453727,0.7385986857630265) -- (17.401620750159907,0.7188850180608028) -- (17.313377320122584,0.6991968521144738) -- (17.224513635921483,0.6795341274496867) -- (17.135069298136102,0.6598967835920764) -- (17.045083907345713,0.6402847600672885) -- (16.954597064129928,0.620697996400958) -- (16.863648369068017,0.6011364321187322) -- (16.77227742273959,0.5816000067462479) -- (16.68052382572415,0.562088659809147) -- (16.588427178601194,0.5426023308330681) -- (16.496027081950217,0.5231409593436567) -- (16.40336313635038,0.5037044848665495) -- (16.310474942381518,0.48429284692738506) -- (16.217402100622905,0.4649059850518107) -- (16.12418421165404,0.4455438387654631) -- (16.03086087605442,0.4262063475939861) -- (15.937471694403428,0.4068934510630129) -- (15.844056267280564,0.3876050886981943) -- (15.750654195265554,0.3683412000251636) -- (15.65730507893744,0.34910172456956445) -- (15.56404851887595,0.3298866018570372) -- (15.470924115660466,0.3106957714132239) -- (15.377971469870602,0.291529172763763) -- (15.285230182085513,0.27238674543429653) -- (15.192739852884927,0.25326842895046475) -- (15.10054008284834,0.2341741628379097) -- (15.008670472555025,0.21510388662226987) -- (14.917170622584479,0.19605753982918728) -- (14.826080133516427,0.17703506198430397) -- (14.735438605929915,0.15803639261325841) -- (14.64528564040478,0.13906147124169088) -- (14.555660837520293,0.12011023739524518) -- (14.466603797856067,0.10118263059955979) -- (14.378154121991372,0.08227859038027674) -- (14.29035141050582,0.06339805626303452) -- (14.203235263978911,0.044540967773476936) -- (14.116845282990027,0.025707264437242472) -- (14.031221068118555,0.006896885779973161) -- (13.946400269273255,-0.011889075238713076) -- (13.862394851326059,-0.030634338057524246) -- (13.77919865785529,-0.049311906968995345) -- (13.69680515245389,-0.06789456157966356) -- (13.61520779871455,-0.08635508149549764) -- (13.534400060230098,-0.10466624632280741) -- (13.454375400593392,-0.1228008356676753) -- (13.375127283397148,-0.14073162913607007) -- (13.296649172234254,-0.15843140633458574) -- (13.218934530697425,-0.17587294686902055) -- (13.141976822379519,-0.19302903034591168) -- (13.065769510873338,-0.2098724363711142) -- (12.990306059771683,-0.2263759445509379) -- (12.915579932667299,-0.24251233449157894) -- (12.841584593153044,-0.2582543857990629) -- (12.768313504821663,-0.27357487807987013) -- (12.695760131265956,-0.2884465909396283) -- (12.623917936078868,-0.3028423039852157) -- (12.552780382853058,-0.3167347968222032) -- (12.482340935181327,-0.3300968490571279) -- (12.412593056656476,-0.3429012402959586) -- (12.343530210871364,-0.3551207501448914) -- (12.275145861418792,-0.3667281582104067) -- (12.207433471891534,-0.3776962440983027) -- (12.140386505882475,-0.3879977874148892) -- (12.07399842698419,-0.3976055677663055) -- (11.943172784891658,-0.414630957998412) -- (11.814904254356463,-0.4285546516440206) -- (11.689140544120846,-0.4391588855522457) -- (11.565829362927332,-0.4462258965730541) -- (11.444918419518359,-0.44953792155538963) -- (11.326355422636198,-0.44887719734893494) -- (11.210088081023429,-0.444025960802918) -- (11.09606410342235,-0.4347664487664531) -- (10.984231198575287,-0.42088089808939344) -- (10.874537075224794,-0.4021515456207396) -- (10.766929442113195,-0.37836062821017435) -- (10.66135600798296,-0.3492903827069256) -- (10.557764481576442,-0.31472304596044864) -- (10.456102571636023,-0.27444085482019887) -- (10.356317986904116,-0.22822604613529052) -- (10.25835843612316,-0.17586085675529262) -- (10.162171628035594,-0.11712752352977418) -- (10.067705271383716,-0.05180828330773579) -- (9.974906625772348,0.020311396666784276) -- (9.883678873581914,0.09912902658425082) -- (9.793844624978703,0.18396260622603222) -- (9.705217898585317,0.2740683412530416) -- (9.61761271302447,0.36870243732846575) -- (9.530843086918765,0.4671211001145821) -- (9.4447230388908,0.5685805352723037) -- (9.35906658756329,0.6723369484648174) -- (9.273687751558725,0.7776465453534911) -- (9.188400549500045,0.8837655316028759) -- (9.103019000009851,0.9899501128729753) -- (9.01735712171029,1.0954564948269763) -- (8.931228933224531,1.1995408831262466) -- (8.844448453175175,1.301459483434428) -- (8.756829700184369,1.4004685014124334) -- (8.668186692875622,1.4958241427225403) -- (8.57833344987074,1.586782613028845) -- (8.48708398979312,1.6726001179922605) -- (8.394252331264795,1.7525328632727906) -- (8.299652492908365,1.8258370545377147) -- (8.203098493346886,1.891768897444308) -- (8.104404351202845,1.949584597658486) -- (8.00338408509873,1.9985403608407069) -- (7.899895565486588,2.0380649846758843) -- (7.793990385280438,2.0683497196973804) -- (7.68577361499905,2.0897962935837313) -- (7.575350326918908,2.1028064409338185) -- (7.462825593316154,2.1077818963449317) -- (7.348304486466816,2.105124394412087) -- (7.231892078646922,2.0952356697334835) -- (7.113693442132842,2.078517456906866) -- (6.99381364920049,2.0553714905277047) -- (6.8723577721260085,2.02619950519329) -- (6.749430883185653,1.9914032355022755) -- (6.625138054655679,1.9513844160492226) -- (6.499584358812001,1.9065447814323306) -- (6.372874867930875,1.857286066249344) -- (6.245114654288216,1.8040100050957335) -- (6.116408790160506,1.747118332568789) -- (5.98686234782366,1.6870127832671642) -- (5.856580399553707,1.6240950917858754) -- (5.791196445780145,1.5917072347892827) -- (5.725668017627129,1.5587669927222123) -- (5.6600082491287935,1.5253245824087571) -- (5.594230274319614,1.4914302206739194) -- (5.528347227234747,1.4571341243422467) -- (5.462372241907872,1.4224865102378317) -- (5.396318452374146,1.3875375951861315) -- (5.33019899266759,1.352337596011239) -- (5.26402699682302,1.3169367295372467) -- (5.1978155988750245,1.281385212589612) -- (5.1315779328578515,1.245733261992882) -- (5.065327132806317,1.2100310945711499)(8.678678924141101,7.708037885037359);
\draw [-to,line width=0.7pt,color=qqqqff] (7.194134128906732,5.038583211588977) -- (8.375943209142456,4.3814566109406705);
\draw [color=xfqqff](9.9,7.1) node[anchor=north west] {$\boldsymbol{F}_2 \mathbf{n}$};
\draw [color=qqqqff](8.4,4.45) node[anchor=north west] {$\mathbf{n}$};
\draw [color=xfqqff](3.5,5.7) node[anchor=north west] {$\boldsymbol{F}_1 \mathbf{n}$};
\draw [color=qqqqff](6.4,3.5) node[anchor=north west] {$\Xi = V_1 \cap V_2$};
\draw [-to,line width=0.7pt,color=xfqqff] (7.194134128906732,5.038583211588977) -- (5.462856882642039,5.172765134855208);
\draw [-to,line width=0.7pt,color=xfqqff] (7.194134128906732,5.038583211588977) -- (9.901186883420365,6.3984127501328985);
\draw [color=qqwwzz](-0.4229615748641753,3.5438831189007107) node[anchor=north west] {$V_1$};
\draw [color=wqwqwq](12.08682116585635,4.224331552496385) node[anchor=north west] {$V_2$};
\draw [-to,line width=0.7pt] (-6.,-3.) -- (-4.,-3.);
\draw [-to,line width=0.7pt] (-6.,-3.) -- (-6.,-1.);
\draw (-4,-3) node[anchor=west] {$x$};
\draw (-6,-1) node[anchor=south] {$y$};
\begin{scriptsize}
\draw [fill=black] (-6.,-3.) circle (2.5pt);
\end{scriptsize}
\end{tikzpicture}

%% file: figs/beam_dom.tex
\definecolor{xfqqff}{rgb}{0.4980392156862745,0,1}
\definecolor{qqwwzz}{rgb}{0,0.4,0.6}
\begin{tikzpicture}[scale=0.5, line cap=round,line join=round,>=triangle 45,x=1cm,y=1cm]
\fill[line width=0.7pt,color=qqwwzz,fill=qqwwzz,fill opacity=0.1] (6,7) -- (12,9) -- (24,5) -- (18,3) -- cycle;
\fill[line width=0.7pt,fill=black,fill opacity=0.2] (6,7) -- (12,9) -- (12,8) -- (6,6) -- cycle;
\fill[line width=0.7pt,color=qqwwzz,fill=qqwwzz,fill opacity=0.1] (6,6) -- (18,2) -- (18,3) -- (6,7) -- cycle;
\fill[line width=0.7pt,color=xfqqff,fill=xfqqff,fill opacity=0.1] (18,2) -- (24,4) -- (24,5) -- (18,3) -- cycle;
\draw [line width=0.7pt,color=qqwwzz] (6,7)-- (12,9);
\draw [line width=0.7pt,color=qqwwzz] (12,9)-- (24,5);
\draw [line width=0.7pt,color=qqwwzz] (24,5)-- (18,3);
\draw [line width=0.7pt,color=qqwwzz] (18,3)-- (6,7);
\draw [line width=0.7pt,dashed] (12,8)-- (24,4);
\draw [line width=0.7pt] (6,7)-- (12,9);
\draw [line width=0.7pt,dashed] (12,9)-- (12,8);
\draw [line width=0.7pt,dashed] (12,8)-- (6,6);
\draw [line width=0.7pt] (6,6)-- (6,7);
\draw [line width=0.7pt,color=qqwwzz] (6,6)-- (18,2);
\draw [line width=0.7pt,color=qqwwzz] (18,2)-- (18,3);
\draw [line width=0.7pt,color=qqwwzz] (18,3)-- (6,7);
\draw [line width=0.7pt,color=qqwwzz] (6,7)-- (6,6);
\draw [line width=0.7pt,color=xfqqff] (18,2)-- (24,4);
\draw [line width=0.7pt,color=xfqqff] (24,4)-- (24,5);
\draw [line width=0.7pt,color=xfqqff] (24,5)-- (18,3);
\draw [line width=0.7pt,color=xfqqff] (18,3)-- (18,2);
\draw [color=qqwwzz](18.5,8) node[anchor=north west] {$\overline{V} = [0,20]\times[0,10]\times[0,1]$};
\draw [color=xfqqff](21.860768958286833,3.391436319275696) node[anchor=north west] {$A_N^{\boldsymbol{\varphi}} \cap \{\vb{x} \in \partial \Vol \; | \; x = 20\}$};
\draw (7.25,9.6) node[anchor=north west] {$A_D^{\boldsymbol{\varphi}},A_D^{\boldsymbol{\Coss}}$};
\end{tikzpicture}

%% file: figs/bending_err_convergence_do_nothing.tex
\begin{tikzpicture}[scale = 0.6]
    \definecolor{asl}{rgb}{0.4980392156862745,0.,1.} 
    \definecolor{asb}{rgb}{0.,0.4,0.6}       
    \definecolor{aso}{rgb}{0.85,0.45,0.0}   
    \definecolor{asg}{rgb}{0.4,0.55,0.1} 

    \begin{loglogaxis}[
        /pgf/number format/1000 sep={},
        axis lines = left,
        xlabel={degrees of freedom},
        ylabel={relative error},
        xmin=2e3, xmax=9e5,
        ymin=1e-4, ymax=2e-1,
        xtick={1e2,1e3,1e4,1e5,1e6},
        ytick={1e-4,1e-3,1e-2,1e-1},
        legend style={at={(0.01,0.01)},anchor=south west,
        font=\footnotesize,row sep=0.5pt,inner sep=1pt},
        legend cell align={left},
        ymajorgrids=true,
        grid style=dotted,
    ]

        \addplot[color=asg, mark=diamond] coordinates {
            (2430,   0.05946123)
            (8670,   0.01181466)
            (32670,  0.00234013)
            (126750, 0.0006586)
            (499230, 0.00025844)
        };
        \addlegendentry{${\muc}/{\mu}=10^0$}

        \addplot[color=asb, mark=diamond] coordinates {
            (2430,   0.10807718)
            (8670, 0.03269495)
            (32670, 0.0066026)
            (126750, 0.0013402)
            (499230, 0.00041579)
        };
        \addlegendentry{${\muc}/{\mu}=10^{1}$}

        \addplot[color=aso, mark=diamond] coordinates {
            (2430,  0.15449121)
            (8670, 0.07542921)
            (32670, 0.02425845)
            (126750, 0.00555615)
            (499230, 0.00128236)
        };
        \addlegendentry{${\muc}/{\mu}=10^{2}$}

        \addplot[color=purple, mark=diamond] coordinates {
            (2430,   0.17228798)
            (8670, 0.10487194)
            (32670, 0.04681834)
            (126750, 0.01413465)
            (499230, 0.00347016)
        };
        \addlegendentry{${\muc}/{\mu}=10^{3}$}
        
        \addplot[color=asl, mark=diamond] coordinates {
            (2430,   0.17944315)
            (8670,   0.11305257)
            (32670,  0.05450948)
            (126750, 0.0196018)
            (499230, 0.00575619)
        };
        \addlegendentry{${\muc}/{\mu}=10^{4}$}

        
        \addplot[color=asg, mark=diamond, mark options={solid}, dotted] coordinates {
            (2430,   0.05941903)
            (8670, 0.01179533)
            (32670, 0.00232895)
            (126750, 0.00065145)
            (499230, 0.0002554)
        };

        \addplot[color=asb, mark=diamond, mark options={solid}, dotted] coordinates {
            (2430,   0.10773613)
            (8670, 0.03250879)
            (32670, 0.00655011)
            (126750, 0.00132574)
            (499230, 0.00041132)
        };

        \addplot[color=aso, mark=diamond, mark options={solid}, dotted] coordinates {
            (2430,  0.15352687)
            (8670, 0.0742539)
            (32670, 0.02380206)
            (126750, 0.0054877)
            (499230, 0.0012728)
        };

        \addplot[color=purple, mark=diamond, mark options={solid}, dotted] coordinates {
            (2430,   0.1715967)
            (8670,  0.10373776)
            (32670, 0.04612949)
            (126750, 0.0139269)
            (499230, 0.00343488)
        };
        
        \addplot[color=asl, mark=diamond, mark options={solid}, dotted] coordinates {
            (2430,   0.17917058)
            (8670, 0.11282937)
            (32670, 0.05437231)
            (126750, 0.01950916)
            (499230, 0.00570693)
        };

        \addplot[dashed,color=black, mark=none]
    				coordinates {
    					(90000, 0.25e-2)
    					(400000, 0.0005625000000000001)
    	};   

        \addplot[dashed,color=black, mark=none]
    				coordinates {
    					(90000, 0.35e-1)
    					(400000, 0.010097646034997484)
    	}; 

        \addplot[dashed,color=black, mark=none]
    				coordinates {
    					(5000, 0.175)
    					(20000, 0.0875)
    	}; 
    \end{loglogaxis}
    \draw (4.65,1.8) node[anchor=south west]{$_{\mathcal{O}(h^{3})}$};

    \draw (4.65,3.875) node[anchor=south west]{$_{\mathcal{O}(h^{2.5})}$};

    \draw (1.65,5.15) node[anchor=south west]{$_{\mathcal{O}(h^{1.5})}$};
\end{tikzpicture}

%% file: figs/bending_err_convergence_interp_Ned2_no_polar.tex
\begin{tikzpicture}[scale = 0.6]
    \definecolor{asl}{rgb}{0.4980392156862745,0.,1.} 
    \definecolor{asb}{rgb}{0.,0.4,0.6}       
    \definecolor{aso}{rgb}{0.85,0.45,0.0}   
    \definecolor{asg}{rgb}{0.4,0.55,0.1} 

    \begin{loglogaxis}[
        /pgf/number format/1000 sep={},
        axis lines = left,
        xlabel={degrees of freedom},
        ylabel={relative error},
        xmin=2e3, xmax=9e5,
        ymin=1e-4, ymax=2e-1,
        xtick={1e2,1e3,1e4,1e5,1e6},
        ytick={1e-4,1e-3,1e-2,1e-1},
        legend style={at={(0.01,0.01)},anchor=south west,
        font=\footnotesize,row sep=0.5pt,inner sep=1pt},
        legend cell align={left},
        ymajorgrids=true,
        grid style=dotted,
    ]

        \addplot[color=asg, mark=diamond] coordinates {
            (2430, 0.054991433009017446)
            (8670, 0.011499958270500075)
            (32670, 9.554672755223167e-05)
            (126750, 0.0014152169791283422)
            (499230, 0.001929466281893856)
        };
        \addlegendentry{${\muc}/{\mu}=10^0$}

        \addplot[color=asb, mark=diamond] coordinates {
            (2430, 0.07002567401229164)
            (8670, 0.02183502426386064)
            (32670, 0.0033407316275040964)
            (126750, 0.0006570894330622068)
            (499230, 0.001764928277734274)
        };
        \addlegendentry{${\muc}/{\mu}=10^{1}$}

        \addplot[color=aso, mark=diamond] coordinates {
            (2430, 0.07975904147207828)
            (8670, 0.029733521370659436)
            (32670, 0.007072121489582208)
            (126750, 0.0005181530752424125)
            (499230, 0.0014960806971915704)
        };
        \addlegendentry{${\muc}/{\mu}=10^{2}$}

        \addplot[color=purple, mark=diamond] coordinates {
            (2430, 0.0860834551685169)
            (8670, 0.03308644346695745)
            (32670, 0.008534925624462402)
            (126750, 0.0009646215006099791)
            (499230, 0.0013787505215284378)
        };
        \addlegendentry{${\muc}/{\mu}=10^{3}$}
        
        \addplot[color=asl, mark=diamond] coordinates {
            (2430, 0.09526754434591693)
            (8670, 0.03600943197497804)
            (32670, 0.009337927882680093)
            (126750, 0.0011305619767233933)
            (499230, 0.0013532237036167018)
        };
        \addlegendentry{${\muc}/{\mu}=10^{4}$}

        
        \addplot[color=asg, mark=diamond, mark options={solid}, dashed] coordinates {
            (2430, 0.054898546269953295)
            (8670, 0.011395502814468722)
            (32670, 1.5136953607770984e-05)
            (126750, 0.0014918196096570624)
            (499230, 0.001959706308753101)
        };

        \addplot[color=asb, mark=diamond, mark options={solid}, dashed] coordinates {
            (2430, 0.06993108379747119)
            (8670, 0.021724050504845722)
            (32670, 0.0032156228856264762)
            (126750, 0.0007553533931186989)
            (499230, 0.0018082596165093232)
        };

        \addplot[color=aso, mark=diamond, mark options={solid}, dashed] coordinates {
            (2430, 0.07965823818207103)
            (8670, 0.029613217176909595)
            (32670, 0.006931348469621979)
            (126750, 0.0004052568625316283)
            (499230, 0.0015485418734518452)
        };

        \addplot[color=purple, mark=diamond, mark options={solid}, dashed] coordinates {
            (2430, 0.08594312060005604)
            (8670, 0.032954960116715615)
            (32670, 0.008386538572031657)
            (126750, 0.0008470598872233187)
            (499230, 0.0014349286769052398)
        };
        
        \addplot[color=asl, mark=diamond, mark options={solid}, dashed] coordinates {
            (2430, 0.09509462043863386)
            (8670, 0.03587214199583815)
            (32670, 0.009187432169995796)
            (126750, 0.0010121875615607663)
            (499230, 0.0014101624327699193)
        };

        \addplot[dashed,color=black, mark=none]
    				coordinates {
    					(15000, 0.275e-1)
    					(95000, 0.0043421052631578945)
    	}; 
    \end{loglogaxis}

    \draw (2.95,3.4) node[anchor=south west]{$_{\mathcal{O}(h^{3})}$};
\end{tikzpicture}

%% file: figs/bending_err_convergence_interp_Ned2.tex
\begin{tikzpicture}[scale = 0.6]
    \definecolor{asl}{rgb}{0.4980392156862745,0.,1.} 
    \definecolor{asb}{rgb}{0.,0.4,0.6}       
    \definecolor{aso}{rgb}{0.85,0.45,0.0}   
    \definecolor{asg}{rgb}{0.4,0.55,0.1} 

    \begin{loglogaxis}[
        /pgf/number format/1000 sep={},
        axis lines = left,
        xlabel={degrees of freedom},
        ylabel={relative error},
        xmin=2e3, xmax=9e5,
        ymin=1e-4, ymax=2e-1,
        xtick={1e2,1e3,1e4,1e5,1e6},
        ytick={1e-4,1e-3,1e-2,1e-1},
        legend style={at={(0.01,0.01)},anchor=south west,
        font=\footnotesize,row sep=0.5pt,inner sep=1pt},
        legend cell align={left},
        ymajorgrids=true,
        grid style=dotted,
    ]

        \addplot[color=asg, mark=diamond] coordinates {
            (2430,   0.03454453)
            (8670, 0.00957251)
            (32670, 0.00218921)
            (126750, 0.00052937)
            (499230, 0.00013112)
        };
        \addlegendentry{${\muc}/{\mu}=10^0$}

        \addplot[color=asb, mark=diamond] coordinates {
            (2430,  0.04469284)
            (8670, 0.01640196)
            (32670, 0.0046823)
            (126750 , 0.00118065)
            (499230, 0.00027167)
        };
        \addlegendentry{${\muc}/{\mu}=10^{1}$}

        \addplot[color=aso, mark=diamond] coordinates {
            (2430,  0.05485501)
            (8670, 0.02306461)
            (32670, 0.00780509)
            (126750, 0.00221239)
            (499230, 0.00050258)
        };
        \addlegendentry{${\muc}/{\mu}=10^{2}$}

        \addplot[color=purple, mark=diamond] coordinates {
            (2430,   0.06274021)
            (8670, 0.02634012)
            (32670, 0.00920247)
            (126750, 0.00262052)
            (499230, 0.00056488)
        };
        \addlegendentry{${\muc}/{\mu}=10^{3}$}
        
        \addplot[color=asl, mark=diamond] coordinates {
            (2430,   0.0736834)
            (8670 , 0.02921091)
            (32670, 0.00997986)
            (126750, 0.00278942)
            (499230, 0.00056151)
        };
        \addlegendentry{${\muc}/{\mu}=10^{4}$}

        
        \addplot[color=asg, mark=diamond, mark options={solid}, dashed] coordinates {
            (2430,  0.03446729)
            (8670, 0.00948008)
            (32670, 0.00209476)
            (126750, 0.00047065)
            (499230, 0.00010878)
        };

        \addplot[color=asb, mark=diamond, mark options={solid}, dashed] coordinates {
            (2430,   0.04461782)
            (8670, 0.0163005)
            (32670, 0.00456031)
            (126750, 0.00109444)
            (499230, 0.00023577)
        };

        \addplot[color=aso, mark=diamond, mark options={solid}, dashed] coordinates {
            (2430,  0.05477565)
            (8670, 0.02295474)
            (32670, 0.00766923)
            (126750, 0.0021095)
            (499230, 0.00045687)
        };

        \addplot[color=purple, mark=diamond, mark options={solid}, dashed] coordinates {
            (2430,   0.06263216)
            (8670, 0.0262179)
            (32670 , 0.00906376)
            (126750, 0.00251465)
            (499230, 0.00051569)
        };
        
        \addplot[color=asl, mark=diamond, mark options={solid}, dashed] coordinates {
            (2430,   0.07355954)
            (8670, 0.02908428)
            (32670, 0.00984057)
            (126750, 0.00268329)
            (499230, 0.00051167)
        };

        \addplot[dashed,color=black, mark=none]
    				coordinates {
    					(75000, 0.065e-1)
    					(450000, 0.0010833333333333333)
    	}; 

    \end{loglogaxis}

    \draw (4.65,2.4) node[anchor=south west]{$_{\mathcal{O}(h^{3})}$};
\end{tikzpicture}

%% file: figs/huge_bending_convergence_muc1e0mu.tex
\begin{tikzpicture}[scale = 0.6]
    \definecolor{asl}{rgb}{0.4980392156862745,0.,1.} 
    \definecolor{asb}{rgb}{0.,0.4,0.6}       
    \definecolor{aso}{rgb}{0.85,0.45,0.0}   
    \definecolor{asg}{rgb}{0.4,0.55,0.1} 

    \begin{loglogaxis}[
        /pgf/number format/1000 sep={},
        axis lines = left,
        xlabel={degrees of freedom},
        ylabel={relative error},
        xmin=2e3, xmax=5e5,
        ymin=1e-4, ymax=4e-1,
        xtick={1e2,1e3,1e4,1e5,1e6},
        ytick={1e-4,1e-3,1e-2,1e-1},
        legend style={at={(0.02,0.02)},anchor=south west,
        font=\footnotesize,row sep=0.5pt,inner sep=1pt},
        legend cell align={left},
        ymajorgrids=true,
        grid style=dotted,
    ]
        \addplot[color=asl, mark=diamond] coordinates {
            (2430, 0.08391936907645424)
            (15606, 0.011161423906816107)
            (84942, 0.0016289647110634018)
            (430950, 0.00028662804631112863)
        };
        \addlegendentry{No interpolation}

        \addplot[color=aso, mark=diamond] coordinates {
            (2430, 0.17179559704883882)
            (15606, 0.021531738279104286)
            (84942, 0.0017934717176589941)
            (430950, 0.0001057873856988432)
        };
        \addlegendentry{Interpolation}

        \addplot[color=asb, mark=diamond] coordinates {
            (2430, 0.06629121219672165)
            (15606, 0.007954814673854288)
            (84942, 0.0011688077905974373)
            (430950, 0.00018256610127953286)
        };

        \addplot[color=asl, mark=diamond, mark options={solid}, dotted] coordinates {
            (2430, 0.08389558444167)
            (15606, 0.011157957731595548)
            (84942, 0.0016278801799527814)
            (430950, 0.0002858164746161008)  
        };

        \addplot[color=aso, mark=diamond, mark options={solid}, dotted] coordinates {
            (2430, 0.1717583968103551)
            (15606, 0.021509246073415273)
            (84942, 0.0017805104142499713)
            (430950, 9.766637352465901e-05)
        };
        
        \addplot[color=asb, mark=diamond, mark options={solid}, dotted] coordinates {
            (2430, 0.06624293459285921)
            (15606, 0.007932795785821371)
            (84942, 0.001156361236450671)
            (430950, 0.0001745016780453373)
        };
        \addlegendentry{Interp. and pol. decomp.}

        \addplot[dashed,color=black, mark=none]
    				coordinates {
    					(30000, 0.1e-1)
    					(350000, 0.0008571428571428571)
    	}; 
    \end{loglogaxis}

    \draw (4,2.5) node[anchor=south west]{$_{\mathcal{O}(h^{3})}$};
\end{tikzpicture}

%% file: figs/huge_bending_convergence_muc1e4mu.tex
\begin{tikzpicture}[scale = 0.6]
    \definecolor{asl}{rgb}{0.4980392156862745,0.,1.} 
    \definecolor{asb}{rgb}{0.,0.4,0.6}       
    \definecolor{aso}{rgb}{0.85,0.45,0.0}   
    \definecolor{asg}{rgb}{0.4,0.55,0.1} 

    \begin{loglogaxis}[
        /pgf/number format/1000 sep={},
        axis lines = left,
        xlabel={degrees of freedom},
        ylabel={relative error},
        xmin=2e3, xmax=5e5,
        ymin=1e-4, ymax=4e-1,
        xtick={1e2,1e3,1e4,1e5,1e6},
        ytick={1e-4,1e-3,1e-2,1e-1},
        legend style={at={(0.02,0.02)},anchor=south west,
        font=\footnotesize,row sep=0.5pt,inner sep=1pt},
        legend cell align={left},
        ymajorgrids=true,
        grid style=dotted,
    ]
        \addplot[color=asl, mark=diamond] coordinates {
            (2430, 0.34396773202945286)
            (15606, 0.09603720702851098)
            (84942, 0.03456242661971009)
            (430950, 0.012899849312091772)
        };
        \addlegendentry{No interpolation}

        \addplot[color=aso, mark=diamond] coordinates {
            (2430, 0.29956215834340805)
            (15606, 0.06406473445338196)
            (84942, 0.012480457790453892)
            (430950, 0.0021218404474742463)
        };
        \addlegendentry{Interpolation}

        \addplot[color=asb, mark=diamond] coordinates {
            (2430, 0.23901997445658937)
            (15606, 0.050167248868576306)
            (84942, 0.010951732047557025)
            (430950, 0.0020642088183211614)
        };
        \addlegendentry{Interp. and pol. decomp.}
        
        \addplot[color=asl, mark=diamond, mark options={solid}, dotted] coordinates {
            (2430, 0.3368805181014554)
            (15606, 0.09386304750934964)
            (84942, 0.03426352620550782)
            (430950, 0.012768337597705471)
            
        };

        \addplot[color=aso, mark=diamond, mark options={solid}, dotted] coordinates {
            (2430, 0.2994448335477453)
            (15606, 0.0640233337595912)
            (84942, 0.012458698309394783)
            (430950, 0.0021059533224272726)
        };
        
        \addplot[color=asb, mark=diamond, mark options={solid}, dotted] coordinates {
            (2430, 0.2389166239631516)
            (15606, 0.050132280774104185)
            (84942, 0.01093178658341685)
            (430950, 0.0020479056040956555)
        };
        
        \addplot[dashed,color=black, mark=none]
    				coordinates {
    					(15000, 0.1e-1)
    					(95000, 0.0015789473684210526)
    	};

        \addplot[dashed,color=black, mark=none]
    				coordinates {
    					( 10000, 2e-1)
    					(160000, 0.03149802624737184)
    	}; 
    \end{loglogaxis}

    \draw (3.4,2.4) node[anchor=south west]{$_{\mathcal{O}(h^{3})}$};

    \draw (3.6,4.4) node[anchor=south west]{$_{\mathcal{O}(h^{2})}$};
\end{tikzpicture}

%% file: figs/torsion_beam.tex
\definecolor{xfqqff}{rgb}{0.4980392156862745,0,1}
\definecolor{qqwwzz}{rgb}{0,0.4,0.6}
\begin{tikzpicture}[scale=0.6, line cap=round,line join=round,>=triangle 45,x=1cm,y=1cm]
\fill[line width=0.7pt,fill=black,fill opacity=0.2] (4,5) -- (4,3) -- (2,2) -- (2,4) -- cycle;
\fill[line width=0.7pt,color=qqwwzz,fill=qqwwzz,fill opacity=0.1] (2,2) -- (18,-2) -- (18,0) -- (2,4) -- cycle;
\fill[line width=0.7pt,color=qqwwzz,fill=qqwwzz,fill opacity=0.1] (18,0) -- (20,1) -- (4,5) -- (2,4) -- cycle;
\fill[line width=0.7pt,color=xfqqff,fill=xfqqff,fill opacity=0.1] (18,-2) -- (20,-1) -- (20,1) -- (18,0) -- cycle;
\draw [line width=0.7pt,dashed] (4,5)-- (4,3);
\draw [line width=0.7pt,dashed] (4,3)-- (2,2);
\draw [line width=0.7pt] (2,2)-- (2,4);
\draw [line width=0.7pt] (2,4)-- (4,5);
\draw [line width=0.7pt,color=qqwwzz] (2,2)-- (18,-2);
\draw [line width=0.7pt,color=qqwwzz] (18,-2)-- (18,0);
\draw [line width=0.7pt,color=qqwwzz] (18,0)-- (2,4);
\draw [line width=0.7pt,color=qqwwzz] (2,4)-- (2,2);
\draw [line width=0.7pt,color=qqwwzz] (18,0)-- (20,1);
\draw [line width=0.7pt,color=qqwwzz] (20,1)-- (4,5);
\draw [line width=0.7pt,color=qqwwzz] (4,5)-- (2,4);
\draw [line width=0.7pt,color=qqwwzz] (2,4)-- (18,0);
\draw [line width=0.7pt,color=xfqqff] (18,-2)-- (20,-1);
\draw [line width=0.7pt,color=xfqqff] (20,-1)-- (20,1);
\draw [line width=0.7pt,color=xfqqff] (20,1)-- (18,0);
\draw [line width=0.7pt,color=xfqqff] (18,0)-- (18,-2);
\draw [line width=0.7pt,dashed] (4,3)-- (20,-1);
\draw [color=xfqqff](20.15,0.8) node[anchor=north west] {$A_{N}^{\boldsymbol{\varphi}} \cap \{\vb{x} \in \partial \Vol \; | \; x = 8\}$};
\draw [color=xfqqff](20.15,0) node[anchor=north west] {$A_{N}^{\Coss} \cap \{\vb{x} \in \partial \Vol \; | \; x = 8\}$};

\draw [color=qqwwzz](10.65,4.5) node[anchor=north west] {$V = [0,8]\times[0,1]^2$};

\draw (0.5,5.5) node[anchor=north west] {$A_{D}^{\boldsymbol{\varphi}}, A_{D}^{\boldsymbol{R}}$};
\end{tikzpicture}

%% file: figs/torsion_convergence_muc1e0mu.tex
\begin{tikzpicture}[scale = 0.6]
    \definecolor{asl}{rgb}{0.4980392156862745,0.,1.} 
    \definecolor{asb}{rgb}{0.,0.4,0.6}       
    \definecolor{aso}{rgb}{0.85,0.45,0.0}   
    \definecolor{asg}{rgb}{0.4,0.55,0.1} 

    \begin{loglogaxis}[
        /pgf/number format/1000 sep={},
        axis lines = left,
        xlabel={degrees of freedom},
        ylabel={relative error},
        xmin=5e2, xmax=3e5,
        ymin=5e-5, ymax=6e-1,
        xtick={1e2,1e3,1e4,1e5,1e6},
        ytick={1e-4,1e-3,1e-2,1e-1,1e0},
        legend style={at={(0.02,0.02)},anchor=south west,
        font=\footnotesize,row sep=0.5pt,inner sep=1pt},
        legend cell align={left},
        ymajorgrids=true,
        grid style=dotted,
    ]

        \addplot[color=asl, mark=diamond] coordinates {
            (702, 0.029931028210876003)
            (4326, 0.005725251211751369)
            (30150, 0.0011957491091411099)
            (224646, 0.00026042401597612357)
            
        };
        \addlegendentry{No interpolation}

        \addplot[color=aso, mark=diamond] coordinates {
            (702, 0.061731801915583744)
            (4326, 0.007618570107516762)
            (30150, 0.0005961227628400716)
            (224646, 0.00011095563355964966)
        };
        \addlegendentry{Interpolation}

        \addplot[color=asb, mark=diamond] coordinates {
            (702, 0.036058972244004485)
            (4326, 0.006187124659561197)
            (30150, 0.0009837390115012144)
            (224646, 8.352281650492649e-05)
        };
        \addlegendentry{Interp. and pol. decomp.}

        \addplot[color=asl, mark=diamond, mark options={solid}, dotted] coordinates {
            (702, 0.029757470354748645)
            (4326, 0.005690681737702975)
            (30150, 0.001183305875111422)
            (224646, 0.0002550079525108435)

        };

        \addplot[color=aso, mark=diamond, mark options={solid}, dotted] coordinates {
            (702, 0.061590837579890485)
            (4326, 0.0075378999922182215)
            (30150, 0.0005606178287047742)
            (224646, 0.00012485950424251816)
        };
        
        \addplot[color=asb, mark=diamond, mark options={solid}, dotted] coordinates {
            (702, 0.03584838180419741)
            (4326, 0.0061105311674108575)
            (30150, 0.0009466104174490914)
            (224646, 6.943857191245268e-05)
        };

        \addplot[dashed,color=black, mark=none]
    				coordinates {
    					(1.8e4, 5e-3)
    					(1.8e5, 5e-4)
    	};
    \end{loglogaxis}

    \draw (5,1.8) node[anchor=south west]{$_{\mathcal{O}(h^{3})}$};
\end{tikzpicture}

%% file: figs/torsion_convergence_muc1e4mu.tex
\begin{tikzpicture}[scale = 0.6]
    \definecolor{asl}{rgb}{0.4980392156862745,0.,1.} 
    \definecolor{asb}{rgb}{0.,0.4,0.6}       
    \definecolor{aso}{rgb}{0.85,0.45,0.0}   
    \definecolor{asg}{rgb}{0.4,0.55,0.1} 

    \begin{loglogaxis}[
        /pgf/number format/1000 sep={},
        axis lines = left,
        xlabel={degrees of freedom},
        ylabel={relative error},
        xmin=5e2, xmax=3e5,
        ymin=5e-5, ymax=6e-1,
        xtick={1e2,1e3,1e4,1e5,1e6},
        ytick={1e-4,1e-3,1e-2,1e-1,1e0},
        legend style={at={(0.02,0.02)},anchor=south west,
        font=\footnotesize,row sep=0.5pt,inner sep=1pt},
        legend cell align={left},
        ymajorgrids=true,
        grid style=dotted,
    ]

        \addplot[color=asl, mark=diamond] coordinates {
            (702, 0.5112179066705261)
            (4326, 0.20546649253550825)
            (30150, 0.06028439041316002)
            (224646, 0.011274982566187591)
        };
        \addlegendentry{No interpolation}

        \addplot[color=aso, mark=diamond] coordinates {
            (702, 0.18754779542422587)
            (4326, 0.04436628108557691)
            (30150, 0.008748206589478452)
            (224646, 0.0005981556766323524)
        };
        \addlegendentry{Interpolation}

        \addplot[color=asb, mark=diamond] coordinates {
            (702, 0.14653821845584958)
            (4326, 0.03824419632510072)
            (30150, 0.008764211831162405)
            (224646, 0.0009313661448272678)
        };
        \addlegendentry{Interp. and pol. decomp.};

        \addplot[color=asl, mark=diamond, mark options={solid}, dotted] coordinates {
            (702, 0.5056558504025845)
            (4326, 0.20267810923679844)
            (30150, 0.05973587584880012)
            (224646, 0.011201760926975071)
            
        };

        \addplot[color=aso, mark=diamond, mark options={solid}, dotted] coordinates {
            (702, 0.18670531121678707)
            (4326, 0.04321399454836299)
            (30150, 0.00810888954330362)
            (224646, 0.00034286210884578)
        };
        
        \addplot[color=asb, mark=diamond, mark options={solid}, dotted] coordinates {
            (702, 0.14527440345637072)
            (4326, 0.03709096451069659)
            (30150, 0.008143121141635283)
            (224646, 0.000676893824457753)
        };

        \addplot[dashed,color=black, mark=none]
    				coordinates {
    					(1.5e4, 0.5e-2)
    					(9.5e4, 0.0007894736842105263)
    	};
        
        \addplot[dashed,color=black, mark=none]
    				coordinates {
    					(4.8e4, 5e-2)
    					(1.8e5, 0.016619219144522898)
    	};

        \addplot[dashed,color=black, mark=none]
    				coordinates {
    					(1e3, 6e-1)
    					(5e3, 0.2683281572999748)
    	};

    \end{loglogaxis}

    \draw (1.3,5.3) node[anchor=south west]{$_{\mathcal{O}(h^{1.5})}$};

    \draw (4.9,3.9) node[anchor=south west]{$_{\mathcal{O}(h^{2.5})}$};
    
    \draw (3.4,1.4) node[anchor=south west]{$_{\mathcal{O}(h^{3})}$};
\end{tikzpicture}

%% file: figs/thin_plate.tex
\definecolor{zzzzzz}{rgb}{0.4980392156862745,0,1}
\definecolor{qqwwzz}{rgb}{0,0.4,0.6}
\begin{tikzpicture}[scale=0.5, line cap=round,line join=round,>=triangle 45,x=1cm,y=1cm]
\fill[line width=0.7pt,color=qqwwzz,fill=qqwwzz,fill opacity=0.1] (6,7+0.5) -- (12,9+0.5) -- (24,5+0.5) -- (18,3+0.5) -- cycle;
\fill[line width=0.7pt,fill=black,fill opacity=0.2] (6,7+0.5) -- (12,9+0.5) -- (12,8) -- (6,6) -- cycle;
\fill[line width=0.7pt,color=qqwwzz,fill=qqwwzz,fill opacity=0.1] (6,6) -- (18,2) -- (18,3+0.5) -- (6,7+0.5) -- cycle;
\fill[line width=0.7pt,color=qqwwzz,fill=qqwwzz,fill opacity=0.1] (18,2) -- (24,4) -- (24,5+0.5) -- (18,3+0.5) -- cycle;
\draw [line width=0.7pt,color=qqwwzz] (6,7+0.5)-- (12,9+0.5);
\draw [line width=0.7pt,color=qqwwzz] (12,9+0.5)-- (24,5+0.5);
\draw [line width=0.7pt,color=qqwwzz] (24,5+0.5)-- (18,3+0.5);
\draw [line width=0.7pt,color=qqwwzz] (18,3+0.5)-- (6,7+0.5);
\draw [line width=0.7pt,dashed] (12,8)-- (24,4);
\draw [line width=0.7pt] (6,7+0.5)-- (12,9+0.5);
\draw [line width=0.7pt,dashed] (12,9+0.5)-- (12,8);
\draw [line width=0.7pt,dashed] (12,8)-- (6,6);
\draw [line width=0.7pt] (6,6)-- (6,7+0.5);
\draw [line width=0.7pt,color=qqwwzz] (6,6)-- (18,2);
\draw [line width=0.7pt,color=qqwwzz] (18,2)-- (18,3+0.5);
\draw [line width=0.7pt,color=qqwwzz] (18,3+0.5)-- (6,7+0.5);
\draw [line width=0.7pt,color=qqwwzz] (6,7+0.5)-- (6,6);
\draw [line width=0.7pt,color=qqwwzz] (18,2)-- (24,4);
\draw [line width=0.7pt,color=qqwwzz] (24,4)-- (24,5+0.5);
\draw [line width=0.7pt,color=qqwwzz] (24,5+0.5)-- (18,3+0.5);
\draw [line width=0.7pt,color=qqwwzz] (18,3+0.5)-- (18,2);
\draw [color=qqwwzz](18.5,8.4) node[anchor=north west] {$\overline{V} = [0,2]\times[0,1]\times[0,0.2]$};
\draw (7.25,10.1) node[anchor=north west] {$A_D^{\boldsymbol{\varphi}},A_D^{\boldsymbol{\Coss}}$};

\draw [color=zzzzzz](15,5.8) node[] {$\vb{m}$};

\draw [dashed,-to,shift={(15,5.8)},line width=0.7pt,color=zzzzzz]  plot[domain=-3.5:0.75,variable=\t]({0.75*cos(\t r)},{0.75*sin(\t r)});
\end{tikzpicture}

%% file: figs/spring.tex
\definecolor{xfqqff}{rgb}{0.4980392156862745,0,1}
\definecolor{qqwwzz}{rgb}{0,0.4,0.6}

\pgfmathsetmacro{\l}{100}
\pgfmathsetmacro{\b}{10}
\pgfmathsetmacro{\t}{4}
\pgfmathsetmacro{\sam}{60}

\pgfplotsset{
  colormap={blueviolet}{
    rgb255=(0,0,80)
    rgb255=(40,0,120)
    rgb255=(90,0,180)
    rgb255=(160,80,255)
  }
}

\pgfplotsset{
  colormap={allblack}{rgb255=(0,0,0) rgb255=(0,0,0)}
}

\pgfplotsset{
  colormap={qqwwzz}{rgb=(0,0.4,0.6) rgb=(0,0.4,0.6)}
}

\pgfplotsset{
  colormap={xfqqff}{rgb=(0.4980392156862745,0,1) rgb=(0.4980392156862745,0,1)}
}

\begin{tikzpicture}
\begin{axis}[
  grid=none,
  width=0.9\linewidth,
  height=0.3\linewidth,
  view={10}{50},
  axis lines=none,
  xlabel={}, ylabel={}, zlabel={},
  xtick=\empty,
ytick=\empty,
ztick=\empty,
  colormap name=blueviolet,
]
\addplot3[
  fill opacity=0.1,
  draw opacity=0.05,
  surf,
  domain=0:1,       
  domain y=0:1,     
  samples=\sam,
  samples y=20,
]
({\l*2*x},
 {\b*(-1+2*y)*(1.5+cos(deg(2*pi*(-1+2*x))))},
 {\t*(-5*sin(deg(2*pi*(-1+2*x))))});

 \addplot3[
  fill opacity=0.1,
  draw opacity=0.05,
  surf,
  domain=0:1,       
  domain y=0:1,     
  samples=\sam,
  samples y=20,
]
({\l*2*x},
 {\b*(-1+2*y)*(1.5+cos(deg(2*pi*(-1+2*x))))},
 {\t*(1-5*sin(deg(2*pi*(-1+2*x))))});

 \addplot3[
  fill opacity=0.1,
  draw opacity=0.05,
  surf,
  domain=0:1,       
  domain y=0:1,     
  samples=\sam,
  samples y=2,
]
({\l*2*x},
 {\b*(-1)*(1.5+cos(deg(2*pi*(-1+2*x))))},
 {\t*(y-5*sin(deg(2*pi*(-1+2*x))))});

 \addplot3[
  fill opacity=0.1,
  draw opacity=0.05,
  surf,
  domain=0:1,       
  domain y=0:1,     
  samples=\sam,
  samples y=2,
]
({\l*2*x},
 {\b*(-1+2)*(1.5+cos(deg(2*pi*(-1+2*x))))},
 {\t*(y-5*sin(deg(2*pi*(-1+2*x))))});

 \addplot3[
  point meta=none,
  colormap name=allblack,
  opacity=0.2,
  draw=none,
  surf,
  domain=0:1,       
  domain y=0:1,     
  samples=2,
  samples y=20,
]
({0},
 {\b*(-1+2*y)*(1.5+cos(deg(2*pi*(-1))))},
 {\t*(x-5*sin(deg(2*pi*(-1))))});

 \addplot3[
  point meta=none,
  colormap name=allblack,
  fill opacity=0.1,
  draw opacity=0.05,
  draw=none,
  surf,
  domain=0:1,       
  domain y=0:1,     
  samples=2,
  samples y=20,
]
({0},
 {\b*(-1+2*y)*(1.5+cos(deg(2*pi*(-1))))},
 {\t*(x-5*sin(deg(2*pi*(-1))))});

 \addplot3[
  point meta=none,
  colormap name=qqwwzz,
  fill opacity=0.5,
  draw opacity=0.05,
  draw=none,
  surf,
  domain=0:1,       
  domain y=0:1,     
  samples=2,
  samples y=20,
]
({\l*2},
 {\b*(-1+2*y)*(1.5+cos(deg(2*pi*(-1))))},
 {\t*(x-5*sin(deg(2*pi*(-1))))});
\end{axis}

\draw [color=xfqqff](8,3.1) node[anchor=north west] {
    $\overline{\Vol} \subset \R^3$
    };
    
\draw [color=qqwwzz](12.7,1.5) node[anchor=north west] {
    $\Area_D^{\Coss},$
    };
    
\draw [color=qqwwzz](12.7,1.9) node[anchor=north west] {
    $ \,\Area_N^{\defmap}\cap \{\vb{x} \in \partial \Vol \; | \; x = 200\}$
    };
    
\draw [color=black](-0.15,2.8) node[anchor=north west] {
    $\Area_D^{\Coss}, \,\Area_D^{\defmap}$
    };


    
    
    

\end{tikzpicture}